\documentclass[10pt,twoside]{article}
\usepackage{amssymb}
\usepackage[mathscr]{eucal}
\usepackage{eufrak}
\usepackage{amsmath}
\usepackage{mathrsfs}

\newcommand{\cD}{\mathcal{D}}

\def\rf#1{{\rm(\ref{#1})}}
\def\chiu{\hfill$\displaystyle\vspace{4pt}
\underset\Box\null$}
\def\Pr{{\bf Proof. }}
\def\O{\Omega}

\def\R{\Bbb R}
\def\N{\Bbb N}
\def\o{\"{o}}
\def\à{\`{a}}
\def\è{\`{e}}
\def\ì{\`{i}}
\def\ù{\`{u}}
\def\ò{\`{o}}
\def\é{\'{e}}

\def\dy{\displaystyle}
\def\ve{\varepsilon}

\def\pa{\partial}
\def\be{\begin{equation}}
\def\ba{\begin{array}}
\def\ea{\end{array}}
\def\ee{\end{equation}}
\def\vs1{\vspace{1ex}}
\def\vp{\varphi}
\def\ov{\overline}
\def\po{\partial\Omega}
\font\sc=cmcsc10
\title{A high regularity result of solutions to\\ modified $p$-Stokes equations}
\author{\sc F. Crispo and P. Maremonti
\thanks{
Dipartimento di Matematica e Fisica, Seconda
Universit\`{a} degli Studi di
 Napoli, via Vivaldi 43, 81100 Caserta,
 Italy.
francesca.crispo@unina2.it;
paolo.maremonti@unina2.it}}
\date{}
\begin{document}
\maketitle \noindent{\bf Abstract} - {\small This paper is
concerned with a special elliptic
system, which can be seen as a
perturbed $p$-Laplacean system,
$p\in(1,2)$, and, for its
``shape'', it is close to the
$p$-Stokes system. Since our ``stress tensor'' is
given by means of $\nabla u $ and
not by its symmetric part, then our
system is not a $p$-Stokes system.  Hence, the system is called {\it modified} $p$-Stokes system.  We look
for the high regularity of the
solutions $(u,\pi)$, that is $D^2u,
\nabla\pi \in L^q,q\in(1,\infty)$.
In particular, we get $\nabla
u,\pi\in C^{0,\alpha}$. As far as
we know, such a result of high
regularity is the first concerning
the coupling of unknowns $(u,\pi)$.
However, our result also holds for the $p$-Laplacean, and it is the first high regularity result in unbounded domains. }
 \par\noindent
 \vskip -0.7true cm\noindent
\newcommand{\red}{\protect\bf}
\renewcommand\refname{\centerline
{\red {\normalsize \bf References}}}
\newtheorem{ass}
{\bf Assumption}[section]
\newtheorem{defi}
{\bf Definition}[section]
\newtheorem{tho}
{\bf Theorem}[section]
\newtheorem{rem}
{\sc Remark}[section]
\newtheorem{lemma}
{\bf Lemma}[section]
\newtheorem{coro}
{\bf Corollary}[section]
\newtheorem{prop}
{\bf Proposition}[section]
\renewcommand{\theequation}{\thesection.\arabic{equation}}
\setcounter{section}{0}
\numberwithin{equation}{section}
\section{ Introduction}\label{intro}
In this paper we study the regularity of solutions to the following problem\be\label{stokes}
\nabla \cdot {\mathbb S}-\nabla\pi= f\,,\quad \nabla \cdot u=0\ \mbox{ in } \R^n\,,\ n\geq 3\,,\ee
where $u=(u_1,\cdots,u_n)$ is a vector field,  $\pi$ is a scalar field,   $S(\nabla u)$ is a special tensor of the kind
\be\label{MS2}{\mathbb S}(\nabla u)=(\mu+|\nabla u|^2)^\frac{p-2}{2} \,\nabla u \,,\ee
 $p\in
(1,2)$, $\mu\geq 0$.
The first contributions to the study of \eqref{stokes} with ${\mathbb S}$ given by \eqref{MS2} with $\mu=0$
are due to
J.L. Lions, in the sixtees, see
 \cite{lions1} and \cite{lions}.  The non-degenerate counterpart , i.e. $p>2$, of \eqref{stokes}  is studied by
 O.A. Ladyzhenskaya, see \cite{L1}, \cite{L2} and \cite{L}.
 The tensor ${\mathbb S}(\nabla u)$, of power-law type,  is ``close'' to the well known extra-stress tensor $\widetilde {\mathbb S}$ of non-Newtonian fluids
\be\label{VST}\widetilde {\mathbb S}(\cD u)=(\mu+|\cD u|^2)^\frac{p-2}{2}  \cD u \,,\ee
where  $\cD u=\frac{\nabla u+\nabla u^T}{2}$, which gives rise to the so called $p$-Stokes problem
\be\label{pS}
\nabla \cdot \widetilde {\mathbb S}(\nabla u)
-\nabla \pi=
f\,,\quad \nabla \cdot u=0 \ \mbox{ in }
\R^n\,.\ee
 Hence, we call, on the contrary, system \eqref{stokes} {\it modified} $p$-Stokes system.
 The difference between ${\mathbb S}$ and $\widetilde {\mathbb S}$ makes our results
 not interesting from a fluid dynamic point of
 view. However it is important to stress that the
 chief task of the paper is not to prove the
 regularity of solutions to the $p$-Stokes
 problem. Our task is to understand if the
 couple of the vector field $u$ and of the
 scalar field $\pi$ can be a smooth solution of
 \eqref{stokes} at least in the special case
 \eqref{MS2}. In other words, we
 investigate if, in spite of the presence of the
 scalar field $\pi$, the regularity of the field
 $u$ can be compared with the one of solutions
 to elliptic systems.
  \par The literature concerning the high regularity of solutions to the corresponding $p$-Laplacean equations  and systems in bounded domains is very rich. We refer to the
paper \cite{Ming} for a complete survey of the present status of regularity for the $p$-Laplacean and for more general opertors.
For results addressing the particular issue of
integrability of the second derivatives of solutions to the
$p$-Laplacean system we mention the papers \cite{acerbif},  \cite{BDVCRNproc}, \cite{BDVCRNLA}, \cite{CM3}, \cite{CM2}, \cite{IM}, \cite{LB},  \cite{Tolk}.  In connection with the $p$-Laplacean, our results enrich the existing literature, furnishing the first high regularity result in unbounded domains and, for $\mu>0$, the first existence result too.
\par
In contrast, as far as we know the literature concerning  the high
  regularity of solutions to the $p$-Stokes problem, in particular considering the singular case $\mu=0$,
   is not satisfactory, in the sense that the regularity of solutions is partial. In this regard, for the singular case we are just aware of the papers \cite{BDR}, \cite{Eb} and \cite{NW}, where at most the $L^{\frac{3p}{p+1}}$-integrability ($\frac{3p}{p+1}<2$) of second  derivatives are obtained in the special cases, respectively, of a space-periodic domain, of a bounded domain with slip boundary conditions, and in the interior. The corresponding evolutionary problem is studied in \cite{BDR} and \cite{MRR}.
   \par To better explain our results, we begin with the following
\begin{defi}[High regular
solution]\label{DS}{\sl Given a distribution
$f$,  by high regular solution to system \eqref{stokes} we mean a pair
$(u,\pi)$ such that
\begin{itemize}\item[i)] for some $q>n$,
$D^2u,\nabla \pi\in L^q(\R^n),\,\nabla u\in
L^p(\R^n)$, $\pi\in L^{p'}(\R^n),\ \frac
1p+\frac{1}{p'}=1\,,$ \item[ii)] $\nabla\cdot
u=0,$\item[iii)] $-({\mathbb S}(\nabla
u),\nabla\varphi)=(\nabla\pi,
\varphi)+(f,\varphi),\mbox{ for all }\varphi\in
C_0^\infty(\R^n).$\end{itemize}}
\end{defi}
\par
 Let $r\in (1,+\infty)$ and let $r'$ be its conjugate exponent. We set
 \be\label{Cri}
M(r)=1-(2-p)H(r')(5+H(r)),\ee
\be\label{cseg}\ov M(2):=2p-3-(2-p)(1+H(2))>0,\ee
where $H(s)$ is the $L^s$-singular transform norm of Calder\'on-Zygmund type. Hence $H(s)\leq C(s-1)$ if $s\geq 2$, $H(s)\leq \frac{C}{s-1}$ if $s\in (1,2]$, and $C=C(n)$ is a numerical constant depending on the space dimension $n$ (see \cite{stein}, chap. 2).
\par The aim of this paper is to show the following results.
\begin{tho}\label{mainTT}{\sl  Let ${\mathbb S} $ be as in \eqref{MS2}, with $p\in (1,2)$ and $\mu\geq 0$.
 Let $q\in
(n,+\infty)$, $q_1=\frac{np}{n+p}$, and
assume that $M(q_1)$, $M(q)$ and $\ov M(2)$ are positive constants.
Let
$f\in L^{q}(\R^n)\cap
L^{q_1}(\R^n)$.
Then there exists a high
regular  solution $(u,\pi)$
of system \eqref{stokes},
with  $\mathbb S$ given by \eqref{MS2}.  Further
the following estimates hold
  \be\label{mainETT1}   \|D^2 u\|_q+ \|D^2 u\|_{q_1}\leq \dy c\, \mu^{\frac{ 2-p}{2} } (\|f\|_{q_1}+\|f\|_q) +c
  \, (\|f\|_{q_1}+\|f\|_q) ^\frac{1}{p-1}\,,   \ee
 \be\label{mainETT2} \big\|{\nabla \pi}\big\|_{q_1}\leq c\, \|f\|_{q_1} \,, \mbox{ and }
 \big\|{\nabla\pi }\big\|_{q}\leq c\, \|f\|_{q}\,. \ee
Moreover, the solution
$(u,\pi)$ is  unique in the class of weak
solutions with $\nabla u\in L^p(\R^n)$. In particular $\nabla u, \pi \in
C^{0,\alpha}(\R^n)$, $\alpha=1-\frac nq$. } \end{tho}
\begin{tho}\label{maindual}{\sl Let ${\mathbb S} $ be as in \eqref{MS2}, with $p\in (1,2)$ and $\mu=0$.
Let $q\in (n,+\infty)$,  $q_1=\frac{np}{n+p}$,  $q_2=\frac{np}{np-n+p}$, and
assume that $M(q_1)$, $M(q)$ and $\ov M(2)$ are positive constants.  Let
$f\in L^{q}(\R^n)\cap L^{q_2}(\R^n)$.
Then there exists a solution
$(u,\pi)$, in the sense of
Definition\,\ref{DS} of system
\eqref{stokes}. Further
the following estimates hold
\be\label{EDUAL}\|\nabla u\|_p\leq c\|f\|_{q_2}^{\frac{1}{p-1}}\,,\ee
  \be\label{Edual1}
     \big\|{D^2u}\big\|_{q}\leq
 c\|f\|_q(1+ \|f\|_{q_2}^\frac{(1-a)(2-p)}{p-1}+\|f\|_{q}^\frac{a(2-p)}{p-a}
\|f\|_{q_2}^{\frac{(2-p)(1-a)}{p-a}\frac{p}{p-1}})\,,\ee
\be\label{EDUALP}\|\pi\|_{p'}\leq c\|f\|_{q_2}\,,\ee
 \be\label{Edual2}  \big\|{\nabla\pi }\big\|_{q}\leq c\, \|f\|_{q}\,, \ee
with $a=\frac{nq}{nq+q-n}$.  Moreover,
the solution $(u,\pi)$ is  unique
in the class of weak solutions with $\nabla u\in L^p(\R^n)$.
In particular $\nabla u, \pi \in
C^{0,\alpha}(\R^n)$, $\alpha=1-\frac nq$.
}
\end{tho}
\par As far as we know, the regularity results
contained in the above theorems
 are the first high regularity results, in the sense of $D^2u\in L^q(\R^n)$, and also of $C^{1,\alpha}$-regularity
 related to a system connecting $u$ and $\pi$, of $p$-Stokes kind.
 \par The choice of the exponent $q_1$ in Theorem\,\ref{mainTT} has been made, since the domain is unbounded, in order to get
 integrability of lower order derivatives,  precisely in order to get $\nabla u \in L^p(\R^n)$.  This enables us to prove the uniqueness of the solution not only in the existence class, but also if compared with weak
 solutions. Note that $q_1>1$ if and only if $p>\frac{n}{n-1}$, which excludes the value $n=2$. \par
 Theorem\,\ref{mainTT} together with the existence and uniqueness
of a weak solution for $\mu=0$ enable us to show
the existence result for high regular solution
stated in Theorem \ref{maindual}, under the
assumption $f\in L^{q}(\R^n)\cap L^{q_2}(\R^n)$.
Theorem\,\ref{maindual} will be a main tool, in
a forthcoming paper, for the proof of an
existence theorem for the {\it modified}
$p$-Navier-Stokes equation.  \par
 The requests on the exponent $p$ and on the constant $M(r)$
  translate a condition of proximity of $p$ to $2$,
 which is a sufficient condition in
 order to get the following kind of estimate
\be\label{stp}
 \big\|\frac{D^2 u}{(\mu+|\nabla u|^2)}{\atop^{\frac{2-p}{2}}}\big\|_{L^r(\R^n)}\leq c\, \|f\|_{L^r(\R^n)}\,.\ee
 Further, the validity of the result for ${\mathbb S}$ given by 
 \eqref{MS2}, and not for $ \widetilde {\mathbb S}$ given by \eqref{VST}, relies on our use of pointwise estimates of the kind $\frac{|D_{x_i}u_j|}{(\mu+|\nabla u|)}\leq 1\,,$ that clearly do not hold, in general, with  $(\mu+|\cD u|)$ in place of $(\mu+|\nabla u|)$.
 On the other hand, we could obtain the results for $\widetilde {\mathbb S}$
 if we were able to show the following crucial kind of estimates
 $$\int_\O\frac{|\nabla u|^r |D^2 u|^r}{(\mu+|\cD u|^2)}{\atop^{\!\frac{3-p}{2}r}}\, dx
 \leq c\int_\O\frac{|\cD u|^r |D^2 u|^r}
 {(\mu+|\cD u|^2)}{\atop^{\!\frac{3-p}{2}r}}\, dx$$
for $r=2$ and bounded $\O$, and for $r\not=2$ and $\Omega=\R^n$.
 \par We like to observe that our results still hold, with some minor changes in the restrictions on the exponent $p$,
 if we replace the tensor \eqref{MS2} with the following one
\be\label{MS1}{\mathbb S}(\nabla u) =(\mu+
|\nabla u|^2)^\frac{p-2}{2}\,\cD u\,.\ee
However, as kindly remarked by C.R. Grisanti,
the operator \eqref{MS1} is not a monotone
operator. Hence, we do not know a way to prove
the
 uniqueness of the corresponding solution, neither in the class of solution with $\nabla u\in L^p(\R^n)$, neither in the {\it a priori} smaller class of high regular solutions.
This is the reason way we desist from developing here the study of the corresponding system.
 \par
Now we would like to introduce our technique, since it appears original and suitable for problems of $p$-Stokes kind.
In doing this, we would like to highlight the obstacles which prevent us to use methods and ideas
from the classical Stokes theory.
 Firstly we observe that the usual regularity theory approaches the problems starting from weak solutions
 of suitable integral equations. Instead in our approach, that we have partially introduced in the previous paper \cite{CM3}
 concerning the elliptic problem, we directly produce ``smooth'' solutions (smooth in the sense
 of Definition \ref{DS}), and, if it is the case, by uniqueness  we deduce the regularity in the class of weak solution.
Generally speaking, when one employs the Faedo-Galerkin
 method to exhibit  a solution (see \cite{lions}), in some sense one is
 forced to aim at a twofold objective: the former is to drop the pressure field, the latter is to get a uniform bound of a suitable ``energy norm''.
 The first objective is a consequence of the classical Helmohltz-Weyl orthogonality. The second is a consequence of the coercive properties of the chief operator.
If we look for a weak solution, the Helmohltz-Weyl orthogonality is between the fields $u$ and $\nabla\pi$, while the metric concerns
$|\nabla u|_p$. If, for instance,  we look for regularity in $W^{2,2}(\O)$,
 then for the classical Stokes problem we need the Helmohltz-Weyl orthogonality between $P\Delta u$ and $\nabla\pi$. Then, the ``energy metric'' is $\|P\Delta u\|$, which implies, by a suitable estimate, the same bound on the second derivatives, and so on. If we reason in a similar way to obtain the high regularity (limited to the second derivatives) of solutions to the  $p$-Stokes problem or to the
modified $p$-Stokes problem, then the Helmohltz-Weyl orthogonality between $P\Delta u$ and $\nabla\pi$ clearly continues to hold, but
 we should be able to evaluate the quantities $$
\big(\nabla\cdot \widetilde {\mathbb S}(\nabla u), P\Delta u\big)=(f,P\Delta u)\quad \mbox{or}\quad
\big(\nabla\cdot {\mathbb S}(\nabla u), P\Delta u\big)=(f,P\Delta u),$$
where $(\cdot,\cdot)$ is the duality in $L^2$. 
This is the impasse that we meet if we formally reproduce the approach inherited
 from the classical analytic theory of the  Stokes problem. Arguing   in a different
 way, in place of multiplying by $P\Delta u$ and use the orthogonality between $P\Delta u$ and $\nabla \pi$,
 one could try to consider the orthogonality between $P\,(\nabla\cdot \widetilde {\mathbb S}(\nabla u))$ (respectively,
 $P\,(\nabla\cdot{\mathbb S}(\nabla u))$  and $\nabla \pi$. However the corresponding
 ``energy metric'' would be of the kind
 $$\|P\,(\nabla\cdot \widetilde {\mathbb S}(\nabla u))\|\,\ \ (\mbox{respectively } \
 \|P\,(\nabla\cdot {\mathbb S}(\nabla u))\|)\,, $$ from
 which we do not know how to make estimates on the second derivatives of $u$.
\par   So, in this paper, we introduce a new technique for the high regularity, where the modified $p$-Stokes problem is regarded as a suitable perturbed elliptic problem.
 We look for a Faedo-Galerkin approximation which does not preserve
 the null divergence and we bypass the difficulty of the Helmohltz-Weyl orthogonality by means of a suitable ``pressure'' function whose definition is just the one of the pressure field in the case
of ``regular'' solutions of problem \eqref{stokes}.
The advantage to handle a Faedo-Galerkin approximation not preserving
 the null divergence is that now we can employ as special Galerkin basis the eigenfunctions of the Laplace operator, and all the estimates on the second derivatives are in the $L^q$ space and not in the spaces of the hydrodynamic, which intrinsically contain the Helmohltz-Weyl orthogonality.
Summarizing, and roughly speaking, we gain a solution of the problem
\be\label{AM2}\nabla \cdot {\mathbb S}(\nabla u)-\nabla \widetilde \pi=
f\,,\quad\mbox{ in }  \R^n\,,\ \mu>0 \mbox{ in } \mathbb S\,,\ee
with
\be\label{repr}
\widetilde\pi:=
(2-p)\int_{\R^n}D_{y_i}  {\mathcal E} (x-y)
\frac{(D_{y_i}u_j\, \nabla u)(y)}{(\mu+|\nabla u(y)|^2)}{\atop^\frac{4-p}{2}}\cdot D_{y_j}\nabla u(y)\,dy\,.\ee
This representation is in agreement with the results known for the ordinary Stokes problem in $\R^n$,
since formally $\widetilde\pi$ becomes a constant when $p=2$.
A solution of problem \eqref{AM2} appears as no divergence free. Actually $\nabla\cdot u=0$ becomes a compatibility condition between
 equations \eqref{AM2}, the representation formula \eqref{repr} and $\nabla\cdot u\to 0$ as $|x|\to \infty$. Indeed,
 since in our approximation the operator is non-singular ($\mu>0$), more regularity is possible. As a consequence,
 from the definition of $\widetilde\pi$ and from equation \eqref{AM2}, we formally deduce a new equation
$$
\Delta(\nabla\cdot u)+ \frac{(p-2)}{2}\frac{
\nabla(\nabla\cdot u)\cdot \nabla|\nabla u|^2}{(\mu+|\nabla u|^2)}=0\,, \mbox{ in } \R^n,\ \nabla\cdot u\ \to 0 \mbox{ as } |x|\to \infty$$
which, from the maximum principle, ensures that $\nabla\cdot u=0$. So, the solution of \eqref{AM2} is
divergence free. Hence $u$ is solution of problem \eqref{stokes}.

\section{Notation and plan of the paper}\label{notations}
Throughout the paper we will assume $p\in (1,2)$. \par
For $\sigma>0$ let us denote by $B_\sigma=B(O, \sigma)$ the $n$-dimensional open ball
of radius $\sigma$  centered at the origin.\par
 By $\Omega$ we mean a domain in $\R^n$, $n\geq 3$. If there is no danger of confusion, we replace $\int_\O\,dx$ with $\int\, dx$. \par
We define an infinitely differentiable
function $\chi: [0,+\infty) \to [0,1]$ satisfying the conditions
$\chi(x)=1$ for $x\leq 1$, $\chi(x)=0$ for $x\geq 2$. If $\theta$ is a positive
constant and $x$ is a point of $\R^n$ , we let $\chi^\theta(x)=\chi(\frac{|x|}{\theta})$.  \par
By
$L^r(\O)$ and $W^{m,r}(\O)$, $m$ nonnegative integer and
$r\in[1,\infty]$, we denote the usual Lebesgue and Sobolev spaces,
with norms $\|\cdot\|_{L^r(\O)}$ and $\|\cdot\|_{W^{m,r}(\O)}$, respectively.
The $L^2$-norm, $L^r$-norm and $W^{m,r}$-norm on $\O$ will be simply
denoted, respectively,  by $\|\cdot\|$, $\|\cdot \|_{r}$ and $\|\cdot\|_{m,r}$, when no danger of confusion is possible.
 For
$m\geq 0$, $r \geq 1 $ we set
$\displaystyle{\widehat{W}}^{m,r}(\Omega):=\{u\!\in
\!L^1_{loc} (\Omega):D^{\alpha}u\in
L^r(\Omega),|\alpha| = m\}$,
where $D^{\alpha}u$
denotes weak  derivatives of
$u$ of order $|\alpha|$. When $|\alpha|=2$, by $D^2 u$
 we can also mean $D^2_{x_ix_j} u$, $\nabla\nabla u$ and $D_{x_i}\nabla u$. By
 $(\widehat W^{1,r}(\R^n))'$
we denote the space of bounded linear functionals defined on
$\widehat W^{1,r}(\R^n)$ such that $\|f\|_{-1,r'}:=\sup_{u\in \widehat W^{1,r}(\R^n), \|\nabla u\|_r=1} |f(u)|<\infty$, where
we denote by $r'$ the conjugate exponent of $r$, i.e. $\frac 1r+\frac{1}{r'}=1$.\par
Finally, we introduce spaces of solenoidal functions.  We set
${\mathscr C}_0(\O):=\{\vp\in C_0^{\infty}(\O)\!: \nabla\cdot
\vp=0\}$,  $J^r(\O):=\mbox
{completion of }{\mathscr C}_0(\O)\mbox{ in } \|\cdot\|_r$-norm,
$\widehat J^{1,r}(\O):=\mbox{completion of }{\mathscr C}_0(\O)\mbox{ in }$
$\|\nabla\|_r$-norm. \par
We use the symbols $\rightharpoonup$ and $\to$   to denote weak and strong convergences, respectively. \par
As defined in the Introduction, by $H(r)$ we  denote the $L^r$-singular transform norm of Calder\'on-Zygmund kind (see \cite{stein}, Chap. II).\par We shall use the lower case letter
 $c$ to denote a positive constant
  whose numerical value (and dependence on some parameters) is unessential for our aims.
  As well as, we can find in the same line $k>1$ and $k\,c\leq c$.
  \par Let $\mu>0$. For any  $s\in [0,+\infty)$, set
\be\label{amuq}
a_s(\mu,v):=\left(\mu+|\nabla
v|^2\right)^{s}\,.\ee
If $s=\frac{2-p}{2}$ we simply set
\be\label{amup}
a(\mu,v):=a_{\frac{2-p}{2}}(\mu, v)\,.\ee
Further
\be\label{amudelta}
a_s(\mu,v,\delta):=\left(\mu+|J_\delta(\nabla
v)|^2\right)^s\,,\ee where $J_{\delta}$ is the Friedrich's mollifier.
Similarly, by $A^s(\mu,v)$
 we denote the fourth-order tensor
 \be\label{fotin}
 A^s_{ijhk}(\mu,v):=\frac{(\nabla v)_{ij}(\nabla v)_{hk}}{a_s(\mu, v)}, \ee
 where $(\nabla v)_{ij}=D_{x_j} v_i$.
Sometimes we will avoid the index notation and write $ A^s$ as
\be\label{fot}
A^s(\mu,v)=\frac{\nabla v\otimes\nabla v}{a_s(\mu, v)}. \ee
We use the summation convention, according to which
pairs of identical indices imply summation. Finally, given a fourth-order tensor $B$, a third-order tensor $T$ and a vector field $v$,
by $B\cdot T$ we mean $B_{ijhk}T_{jhk}$, by $B\cdot v=v_i B_{ijhk}$ and by $T\cdot v=v_i T_{ijh}$. \par
The proof of Theorem \ref{mainTT} is based on suitable approximations of the solution of system \eqref{stokes}. Each approximation is constructed as a solution of a suitable system.
 Some approximating systems are introduced in order to prove existence and regularity, other approximations are related to the fact that we work in the unbounded domain $\R^n$. Below we introduce all the approximating systems we will use in our construction of the regular solution.
\par Firstly we assume that $f\in J^q(\R^n)\cap J^{q_1}(\R^n)$ and $\mu>0$. For $\rho>0$, let $\chi^\rho(x)$ be a smooth cut-off function.
We introduce the following auxiliary system
\be\label{stokesS1}\ve\Delta u+\frac{\Delta u}{a_{2-p}(\mu,u)}
+(p\!-\!2) \frac{\nabla u\otimes \nabla u}{a_{\frac{4-p}{2}}(\mu, u)}\cdot J_\eta(\frac{\nabla\nabla u\,\chi^\rho}{a(\mu, u)})=\!
\frac{\nabla \Pi(u,\chi^\rho)}{a(\mu,u)}+ \frac{f}{a(\mu,u)}, \mbox{ in }\R^n,
\ee
where $\mu$, $\rho$, $\eta$ and $\ve$ are positive constants,
 and
\be\label{pi}
\Pi(u,\chi^\rho):=(2-p)\int_{\R^n}D_{y_i}  {\mathcal E} (x-y)
\frac{D_{y_i}u_j (y)\nabla u(y)}{a_1(\mu, u(y))}\cdot J_\eta(\frac{D_{y_j}\nabla u(y)\,\chi^\rho(y)}{a(\mu, u(y))})dy\,,\ee
$ {\mathcal E} (x-y)=\frac{|x-y|}{(2-n)\omega_n}{^{^{\!\!\!\!\!2-n}}}$ being the fundamental solution of the Laplace equation $\Delta w=F$.
Roughly speaking the idea is to calculate the divergence in \eqref{stokes}, when $\mu>0$, with the focused introduction of a regularizer, a cut-off function, the diffusion term $\ve\Delta u$ and the ``perturbation term'' $\nabla\Pi(u, \chi^\rho)$ in place of the ``pressure function'' $\nabla\pi$. We explicitly stress that the system is further multiplied by $a(\mu, u)^{-1}$. This is strategic to obtain the estimate \eqref{aip} for the Faedo-Galerkin approximations, which implies other crucial estimates. This last artifice can be removed as soon as the solutions of the
approximating systems are sufficiently smooth.
 Due to the presence of the diffusion term
$\ve\Delta u$, by using the Faedo-Galerkin approximation method on domains $B_\sigma$ invading $\R^n$ and then known results applied
on an approximating linear elliptic system in each $B_\sigma$, we prove  that there exists a solution which satisfies $D^2 u\in L^{q_1 }(\R^n)\cap L^q(\R^n)$, for $q_1=\frac{np}{n+p}$ and some $q>n$.  This, in particular ensures that $\nabla u\in L^p(\R^n)$.
The existence and regularity for solutions on the sequence of invading domains are obtained in sec. \ref{secconvex}, while
the existence and regularity results for solutions of system \eqref{stokesS1} in the whole $\R^n$ are proved in sec. \ref{secnoep}.  Clearly, the solution of this system depends on
 the parameters $\mu$, $\rho$, $\eta$ and $\ve$ and the estimates in the norms of the previous spaces are not uniform in these parameters. Therefore $u=u^{\mu, \rho, \eta, \ve}$.
 We will let these four parameters tend to zero in four different sections, one after the other, proving that suitable norms are bounded with respect to each parameter. For the reader's convenience, in each section we will stress the dependence of the sequence just on the parameter which is
  going to zero.  In sec. \ref{secep} (Proposition \ref{Vteo4}) we prove estimates uniform, at the same time, in $\ve$, $\eta$ and $\rho$.
 This enables us to pass to the limit firstly as $\ve$ goes to zero, and show that the family of solution
$\{u^{\mu, \rho, \eta, \ve}\}$ of system \eqref{stokesS1} tends to a solution, say $u=u^{\mu, \rho, \eta}$, of the following system
\be\label{stokesS4}\frac{\Delta u}{a_{2-p}(\mu,u)}
+(p-2) \frac{\nabla u\otimes \nabla u}{a_{\frac{4-p}{2}}(\mu, u)}\cdot J_\eta(\frac{\nabla\nabla u\,\chi^\rho}{a(\mu, u)})=
\frac{\nabla \Pi(u,\chi^\rho)}{a(\mu,u)}+ \frac{f}{a(\mu,u)}\,,\ \mbox{ in } \R^n\,,
\ee
with $\Pi(u,\chi^\rho)$ given by \eqref{pi}. Before passing to the limit on $\eta$, in sec. \ref{secet} we show that, for any $\eta>0$, the solution of system \eqref{stokesS4} admits third-order derivatives in $L^q(\R^n)$, for some $q>n$, and the estimates are uniform in $\eta$ and $\rho$. This is an important step in order to prove that the solution of system \eqref{stokes} is divergence free. Therefore, still in section sec.\ref{secet} we pass to the limit as $\eta$ goes to zero, and prove that the limit function of the family of solutions
 $\{u^{\mu, \rho, \eta}\}$ tends to a solution, say $u=u^{\mu, \rho}$, which solves the system
 \be\label{stokesS5}\frac{\Delta u}{a(\mu,u)}
+(p-2) \frac{\nabla u\otimes \nabla u}{a_{\frac{4-p}{2}}(\mu, u)}\cdot (\nabla\nabla u\,\chi^\rho)=
\nabla \Pi(u,\chi^\rho)+ f\,,\ \mbox{ in } \R^n\,,
\ee
with
\be\label{pia}
\Pi(u,\chi^\rho):=(2-p)\int_{\R^n}D_{y_i}  {\mathcal E} (x-y)
\frac{(D_{y_i}u_j\, \nabla u)(y)}{a_{\frac{4-p}{2}}(\mu, u(y))}\cdot D_{y_j}\nabla u(y)\,\chi^\rho(y)\,dy\,,\ee
and has the properties stated in Proposition\,\ref{Vteo7} ($D^3u\in L^q(\R^n)$).  Since all the estimates are also uniform in $\rho$, we can pass
to the limit as $\rho$ tends to infinity, and show that the family of solution
$\{u^{\mu, \rho}\}$ tends to a solution, say $u=u^{\mu}$ of the following system
\be\label{stokesS8}\frac{\Delta u}{a(\mu,u)}
+(p-2) \frac{(\nabla u\otimes \nabla u)\cdot \nabla\nabla u}{a_{\frac{4-p}{2}}(\mu, u)}=
\nabla \Pi(u)+ f\,,\ \mbox{ in } \R^n\,,
\ee
where
\be\label{piu}
\Pi(u):=(2-p)\int_{\R^n}D_{y_i}  {\mathcal E} (x-y)
\frac{(D_{y_i}u_j\, \nabla u)(y)}{a_{\frac{4-p}{2}}(\mu, u(y))}\cdot D_{y_j}\nabla u(y)\,dy\,,\ee
We
show that the solutions of system \eqref{stokesS8} have a divergence satisfying
a suitable elliptic system with bounded coefficients. By virtue of a well known maximum principle, it follows that $\nabla\cdot u=0$.   Clearly for each fixed $\mu>0$, the solution $u^\mu$ of \eqref{stokesS8} is also a solution of \eqref{stokes}.
The final steps consist in showing that the one parameter family of solutions $\{u^\mu\}$ to system \eqref{stokesS8}
converges to a solution of system \eqref{stokes} as $\mu\to 0$, in extending the results to $f\in L^q(\R^n)\cap L^{q_1}(\R^n)$, with nonnull divergence,
and in showing that the solution is unique. These steps are proved in sec. \ref{proof}. Sec. \ref{proof1} is concerned with the proof of Theorem\,\ref{maindual}.

\section{\large Some preliminary results}\label{preliminary}
For the reader's convenience, we recall below some well known results and introduce
 some basic estimates.
\par Since we will use the Faedo-Galerkin approximation scheme, we recall the classical ``fixed point'' theorem,
for which we refer to \cite{lions}, Lemma I.4.3.
 \begin{lemma}\label{Lionslemma}{\sl Let $P$ be a continuous function of $\R^k$, $k\geq 1$, into itself such that, for some $R> 0$,
$P(\xi)\cdot \xi\geq 0$,  for all $\xi\in \R^k$,  with $|\xi |=R$. Then there exists a $\xi_0\in \R^k$, with $|\xi_0| \leq R$, such that $P(\xi_0) = 0$.}
\end{lemma}
 \par
 Let us consider the following system
 \be\label{pstokes}
\nabla \cdot {\mathbb S}(\nabla u)-\nabla\pi= f\,,\quad \nabla \cdot u=0\ \mbox{ in } \R^3\,.
 \ee
\begin{defi}\label{DWS}
Given a distribution $f$, by a weak solution of
system \eqref{pstokes} we mean a field $u\in
\widehat J^{1,p}(\R^n)$,  such that
 $$\int {\mathbb S}(\nabla u)\cdot \nabla \varphi
\,dx=\int f\cdot \vp\, dx\,, \ \forall\,\vp\in
{\mathscr C}_0(\R^n)\,.$$
\end{defi}
\begin{lemma}\label{exists}{\sl Let $\mu=0$.
Let $f\in (\widehat W^{1,p}(\R^n))'$. Then there
exists a unique weak solution of \eqref{pstokes},
and the following estimate holds
\be\label{grad}\|\nabla u\|_p\leq \|f\|_{-1,p'}^{\frac{1}{p-1}}\,.\ee}
 \end{lemma}
 \Pr For the result we refer for instance  to \cite{lions} (Chap. 2, sec. 5), whose proof is performed in a bounded domain, but it can be easily extended to the case of $\R^n$. \chiu
 \vskip0.1cm
Further we recall the following regularity theorem, as given in \cite{GiaqMart}, Theorem 7.3.
\begin{lemma}\label{lemmaGM}{\sl
Let $v$ be a $W^{2,2}$-solution of the linear system
$$A_{ij}^{\alpha\beta}D_\alpha D_\beta v^j=f_i,$$
with $A_{ij}^{\alpha\beta}\in
C(\overline\Omega)$ satisfying the
Legendre-Hadamard condition. If $f$
belongs to $L^q(\Omega)$, for some
$q\geq 2$,  then $D^2v\in
L^q(\Omega)$,  with
$$\|D^2 v\|_q \leq C(q,n,L,\omega)\, (\|f\|_q+\|D^2 v\|)\,,$$
where $L$ is the  constant of the
Legendre-Hadamard condition and $\omega$ is the modulus of continuity of $A$.}
\end{lemma}
We recall the following classical inequality, for which we refer, for instance, to \cite{DER}, Lemma 6.3.
\begin{lemma}\label{giu}
Let $p\in (1, +\infty)$, and let $\mu\geq 0$.
There exists a constant $C=C(n,p)$, independent of $\mu$, such that for any $A, B \in \R^{n\times n}$,
$$ ((\mu+ |A|^2)^\frac{p-2}{2}A- (\mu+ |B|^2)^\frac{p-2}{2}B)\cdot(A-B)\geq C\,(\mu + |A|^2 + |B|^2)^\frac{p-2}{2}|A-B|^2\, .$$
\end{lemma}
In Lemma \ref{density} we give a known density result and, for completeness, we also sketch the proof.
This result will be used in the sequel for the whole space $\R^n$. However
we perform the proof  in the more general case of an exterior domain. One of the authors, in \cite{Mar88}, already gave the same result,
 but the proof cantains an oversight. The argument we use relies on standard procedures
for proving density results, see, for instance, \cite{Glibro}, \cite{marem} and  \cite{marem1}.
\begin{lemma}\label{density}{\sl
Let $\O\subseteq\R^n$ be an exterior domain. Then ${\mathscr C}_0(\O)$ is dense in $J^r(\O)\cap J^s(\O)$, for any
$1< r<s<+\infty$.  }
\end{lemma}
\Pr Let $\Omega_{2R}=\Omega\cap B_{2R}$, $R>d=diam (\O^c)$, where $\O^c$ is the complementary set of $\O$ in $\R^n$.
 Let $\chi^R$ be a $C_{0}^{\infty}(\R^n)$
 cut-off function as defined in the notation. For $v\in J^r(\O)\cap J^s(\O)$, denote by $w_R$ the solution of the Bogovskii
 problem \be\label{bog}\ba{rll}\dy \nabla\cdot w_R=&\!\!\!\!\!\dy-v\cdot \nabla \chi^R,&\ \mbox{ in }\ \O_{2R},\\
 w_R=&\dy0, &\ \mbox{on } \ \partial \O_{2R}\,. \ea
 \ee
 Due to the assumptions on $v$ and $\chi^R$,
  the compatibility condition $$\int_{\O_{2R}} v \cdot \nabla\chi^R\,dx=0$$
  is satisfied. Thanks to a well known result due to Bogovskii \cite{bogov} (see also \cite{Glibro}, chap. 3), there exists a solution $w_R\in W_0^{1, r}(\O_{2R})\cap W_0^{1,s}(\O_{2R})$ of the above
  problem, such that
\be\label{perP}\|\nabla w_R\|_{L^r(\O_{2R})}\leq c
\| v\cdot \nabla \chi^R \|_{L^r(\O_{2R})}\leq \frac cR \|v\|_{L^r(R\leq |x|\leq 2R)},\ee
 \be\label{Poi1}
 \|w_R\|_{L^r(\O_{2R})}\leq cR  \|\nabla w_R\|_{L^r(\O_{2R})}\leq c\|v\|_{L^r(R\leq |x|\leq 2R)}\,,\ee
  \be\label{perP3}\|\nabla w_R\|_{L^s(\O_{2R})}\leq c \| v\cdot \nabla \chi^R \|_{L^s(\O_{2R})}\leq \frac cR \|v\|_{L^s(R\leq |x|\leq 2R)},\ee
  \be\label{Poi2}
   \|w_R\|_{L^s(\O_{2R})}\leq cR  \|\nabla w_R\|_{L^s(\O_{2R})}\leq c\|v\|_{L^s(R\leq |x|\leq 2R)}\,,\ee
  with the  constants $c$  independent of $R$, where we have used
  Poincar\'e's inequality in \eqref{Poi1} and \eqref{Poi2}.
 Let us define the function $$v_R=v\chi^R+w_R.$$
  Clearly, from \eqref{bog}--\eqref{Poi2}, $ v_R\in J^r(\O_{2R})\cap J^s(\O_{2R})$. Let be $\ve>0$.
  There exists a $R$ such that
   $\|v\|_{L^r(|x|\geq R)}<\ve$ and
$\|v\|_{L^s(|x|\geq R)}<\ve$.
  Since $ v_R\in J^r(\O_{2R})\cap J^s(\O_{2R})$, there is a
   function $\psi\in {\mathscr C}_0(\O_{2R})\subset {\mathscr C}_0(\O)$
  such that
 $$\|v_R-\psi\|_{L^s(\O_{2R})} <\frac{\ve}{R}{\atop^{n\frac{s-r}{sr}}}$$
 and, from H\o lder's inequality with exponents $s$ and $\frac{sr}{s-r} $,
 one also gets
 $$ \|v_R-\psi\|_{L^r(\O_{2R})} \leq R^{n\frac{s-r}{sr}} \|v_R-\psi\|_{L^s(\O_{2R})}
 <\ve\,.$$
 Let us extend the function $w_R$ and the corresponding function $v_R$ to zero outside $\Omega_{2R}$.
 It is now easy to verify that the function $\psi$ approximates $v$ in the norms $L^{r}(\O)$ and $L^s(\O)$.
 Indeed,
$$\ba{rl}\dy \vs1
 \|v-\psi\|_{s}&\dy\leq \|v-v_R\|_s+\|v_R-\psi\|_{s}\leq \|v(1-\chi^R)\|_{s}+\|w_R\|_s+\frac{\ve}{R}{\atop^{n\frac{s-r}{sr}}}
\\& \dy \leq \|v\|_{L^s(|x|\geq R)}+\|w_R\|_{L^s(\O_{2R})}+\frac{\ve}{R}{\atop^{n\frac{s-r}{sr}}}< 3\frac{\ve}{R}{\atop^{n\frac{s-r}{sr}}}\,,\ea$$
 and
$$\ba{rl}\dy \vs1
 \|v-\psi\|_{r}&\dy\leq \|v-v_R\|_r+\|v_R-\psi\|_{r}\leq \|v(1-\chi^R)\|_{r}+\|w_R\|_r+\ve
\\& \dy \leq \|v\|_{L^r(|x|\geq R)}+\|w_R\|_{L^r(\O_{2R})}+\ve< 3\ve.\ea$$
 \chiu\vskip0.1cm
The following result is completely similar to Lemma 2.1 in \cite{CM2}, for which we refer for further details and for an extension to more general domains. It is a generalization to $p\not=2$ of a well known
estimate, that can be found, for instance in \cite{LU} (chap. 3, sec. 8).
\begin{lemma}\label{LL1}{\sl
Let $\O$ be a bounded convex domain of class $C^2$. Let $\mu>0$ and $p\in (\frac 32,2]$.  Assume that $v \in
W^{2,2}(\O)\cap W_0^{1,2}(\O)$. Then,
$$ \big\|(\mu+|\nabla v|^2)^\frac{p-2}{2}\nabla\nabla v\big\|\!\leq \frac{1}{2p-3}\,\big\| (\mu+|\nabla v|^2)^\frac{p-2}{2}\Delta v\big\|\,.
$$ }
 \end{lemma}
  \Pr  Following \cite{LU}, we just prove the result for sufficiently smooth functions. It can be extended to functions
 in $W^{2,2}(\O)\cap W_0^{1,2}(\O)$ by density arguments.
 So, let $v$ be a function which is continuously differentiable three
times and vanishes on $\po$. Integrations by parts give
\be\label{entr1}
\ba{ll}\vs1 \dy \int_\O(\mu+|\nabla v|^2)^{p-2}|\Delta v|^2\, dx=\!\dy-\int_\O (\mu+|\nabla v|^2)^{p-2} D_{x_k} \Delta v \cdot D_{x_k} v\,dx\\
\hfill \vs1 \! \dy- 2(p-2)\int_\O (\mu+|\nabla v|^2)^{p-3}(\Delta v\cdot D_{x_k} v) \,
(\nabla v\cdot D_{x_k} \nabla v) \,dx \dy\\ \vs1\hfill\dy
+\int_{\po} (\mu+|\nabla v|^2)^{p-2} \Delta
v \cdot \frac{\partial v}{\partial n}\,d\sigma\\ \hfill= \vs1 \dy
\!\int_\O\!\!(\mu+|\nabla v|^2)^{p-2}D^2_{x_jx_k} v\cdot D^2_{x_jx_k} v
\,dx \\\vs1\hfill \dy + 2(p-2)\!\int_\O\!\!(\mu+|\nabla v|^2)^{p-3}\left(D_{x_j} \nabla v\cdot \nabla v\right)^2 dx
\\\vs1\dy\hfill
 -2(p-2)\int_\O (\mu+|\nabla v|^2)^{p-3}(\Delta v\cdot D_{x_k} v)\,(\nabla
v\cdot D_{x_k}\nabla v) \,dx \dy
\\\dy\hfill
+\int_{\po} (\mu+|\nabla
v|^2)^{p-2}\left[\Delta v\cdot \frac{\partial v}{\partial
n}-\frac{\partial^2 v}{\partial x_k\partial n}\cdot\frac{\partial
v}{\pa x_k} \right] \,d\sigma \,.\ea\ee
By
using the arguments of \cite{L}, based on a localization technique,
one can show that the last boundary integral, say $I_{\po}$, is non-negative since
$\O$ is convex.
\par From \eqref{entr1} and recalling position \eqref{amuq} one gets
 \be\label{L2a}\ba{rl}\vs1 \dy (1\!-\!2(2\!-p))\!\!\int _\O\!
\frac{|\nabla\nabla v|^2}{a_{2-p}(\mu,v)}\, dx\leq& \!\!\!\dy\!\!\int_\O\! \frac{|\Delta v|^2}{a_{2-p}(\mu,v)}\,dx\\ \vs1
\dy \!\!\!&\dy\!\!\!\! \!+
2(p-2)\!\int_\O\!\frac{(\Delta v\cdot D_{x_k} v)(\nabla
v\cdot D_{x_k}\!\nabla v)}{a_{3-p}(\mu,v)} \,dx \dy
- I_{\po}\\
\dy \leq &\!\!\!\dy \!\int_\O \frac{|\Delta v|^2}{a_{2-p}(\mu,v)} \,dx+ 2(2-p)\!\int_\O \frac{|\nabla\nabla v|\, |\Delta v|} {a_{2-p}(\mu,v)}\,dx \,.\ea\ee
By applying H\"{o}lder's and  Young's
inequalities to the last integral one readily has
 \be\label{L3a}\ba{ll}\dy
\left[2p-3-2\ve(2-p)^2\right]\!\int _\O \frac{|\nabla\nabla v|^2}{a_{2-p}(\mu,v)} \, dx  \leq
\left(\!1+\frac{1}{2\ve}\right)\! \int_\O\frac{|\Delta v|^2}{a_{2-p}(\mu,v)}\,dx\,,\ea \ee hence \be\label{L3} \int
_\O \frac{|\nabla\nabla v|^2}{a_{2-p}(\mu,v)} \, dx \leq C(\ve)
\int_\O\frac{|\Delta v|^2}{a_{2-p}(\mu,v)}\,dx\,,\ee
with $$C(\ve):=\mbox{\large${\frac{1+2\ve}{2\ve[2p-3-2\ve(2-p)^2]}}$}\,.$$
By an
easy computation, one verifies that the minimum of $C(\ve)$ equals
$1/(2p-3)^2$ and it is attained for $\ve=\frac{2p-3}{2(2-p)}$. Hence the result follows.
 \chiu\vskip0.1cm
\begin{rem}{\rm
The following Lemma \ref{LL2} generalizes the above inequalities to any $L^q$-space, $q>1$, when $\O=\R^n$. Note that Lemma \ref{LL1}
still holds in $\R^n$, for any $p\in (\frac 32, 2]$ and with a better constant, but its proof is strictly connected to the fact that $q=2$.  }
\end{rem}
The following well known result is a main tool in the proof of Lemma\,\ref{LL2} and Lemma\,\ref{LL3}.
\begin{lemma}\label{perLL2}
{\sl Let $h\in C_0^{\infty}(\R^n)$. Then the solution of
\be\label{ric}
\Delta w=h, \mbox{ in } \R^n, \ee
is smooth. In particular, for any $s\in (1,\infty)$ and $|\alpha|\geq0$
 \be\label{Poi}\|D^{2+\alpha} w \|_{s}\leq H(s) \|D^{\alpha}h\|_{s}.\ee}
\end{lemma}

 \begin{lemma}\label{LL2}{\sl
Let $\mu>0$, $p\in (1,\infty)$ and $q\in (1,+\infty)$.  Assume that $(\mu+|\nabla v|^2)^\frac{p-2}{2} D^2v\in L^{q}(\R^n)$,
 and $\nabla v\in L^r(\R^n)$, for some $r\in (1,+\infty)$ if $p\leq 2$, for some $r\in [p',+\infty)$ if $p>2$.
Then
\be\label{ELL2}\big\|(\mu+|\nabla v|^2)^{\frac{p-2}{2}}\nabla\nabla v\big\|_q\,
\leq \frac{H(q')}{1-4 H(q')|2-p|}\big\|(\mu+|\nabla v|^2)^{\frac{p-2}{2}}\Delta v\big\|_q,\ee
provided that $4H(q')|2-p|<1$\,, where $q'$ is the conjugate exponent of $q$.}
 \end{lemma}
 \Pr Let $h$ and $w$ be as in Lemma\,\ref{perLL2}.  Multiplying $ (\mu+|\nabla v|^2)^{\frac{p-2}{2}}D^2_{x_ix_j}v $ by $h$ and integrating by parts in $\R^n$, we obtain
 $$\ba{ll}\vs1 \dy \int (\mu+|\nabla v|^2)^{\frac{p-2}{2}}D^2_{x_ix_j} v \cdot h\, dx  =\dy\!\int(\mu+|\nabla v|^2)^{\frac{p-2}{2}}
 D^2_{x_ix_j} v \cdot \Delta w\, dx\\
 =\dy\!-\!\!\int\! (\mu+|\nabla v|^2)^{\frac{p-2}{2}}
 D_{x_i}v\cdot D_{x_j}\Delta w\,dx+\frac{2-p}{2}\!
\int\! (\mu+|\nabla v|^2)^{\frac{p-4}{2}}
D_{x_i}v\cdot\Delta w\, D_{x_j} |\nabla v|^2\,dx\\ \hfill \dy :=\dy\!-\!\int (\mu+|\nabla v|^2)^{\frac{p-2}{2}}
D_{x_i}v\cdot D_{x_j}\Delta w
\,dx+I_0,
\ea$$
where each term makes sense, thanks to our assumptions on $v$ and the
integrability of $D^2w$ and $D^3 w$. We then integrate three times more by parts, denoting at the $i-th$ integration, by $I_i$ the integral with the
derivatives of the term $a(\mu,v)$. Then using the previous identity we get
\be\label{sep1}
\ba{rl}\vs1 \dy \int (\mu+|\nabla v|^2)^{\frac{p-2}{2}}D^2_{x_ix_j} v \cdot h\, dx= \!\!\!&\dy
\!\!\!\!\int\! (\mu+|\nabla v|^2)^{\frac{p-2}{2}}
D^2_{x_ix_h}v\cdot D^2_{x_jx_h}w\,dx+\sum_{i=0}^1 I_i
\\
\vs1
\hfill =\!\!\!\dy &\dy\!\!\!-\!
\!\int \!\!(\mu+|\nabla v|^2)^{\frac{p-2}{2}}
D_{x_h}v\cdot D^3_{x_ix_jx_h}w\,dx+\!\sum_{i=0}^2 I_i
\\\vs1
\hfill\dy =\!\!\!&\dy\!\!\! \int (\mu+|\nabla v|^2)^{\frac{p-2}{2}}
\Delta v\cdot D^2_{x_ix_j}w\,dx+\sum_{i=0}^3\,I_i
 \,,\ea\ee
where
$$I_1:=\frac{p-2}{2}\int
(\mu+|\nabla v|^2)^{\frac{p-4}{2}}
D_{x_i}v\cdot D^2_{x_jx_h}w\, D_{x_h} |\nabla v|^2\,dx\,, $$
$$I_2:=\frac{2-p}{2}\int (\mu+|\nabla v|^2)^{\frac{p-4}{2}}
D_{x_h}v\cdot D^2_{x_jx_h}w\, D_{x_i} |\nabla v|^2\,dx\,.$$
$$I_3:=\frac{p-2}{2}\int
(\mu+|\nabla v|^2)^{\frac{p-4}{2}}
D_{x_h}v\cdot D^2_{x_ix_j}w\, D_{x_h} |\nabla v|^2\,dx\,.$$
Taking into account estimate \eqref{Poi}, the first term on the right-hand side of \eqref{sep1} is estimated as
\be\label{sep2}
  |\int (\mu+|\nabla v|^2)^{\frac{p-2}{2}}
  \Delta v\cdot D^2_{x_ix_j}w\,dx|\leq H(q')
\big\|(\mu+|\nabla v|^2)^{\frac{p-2}{2}}\Delta v\big\|_q
  \|h\|_{q'}\,.\ee
 Let us estimate the generic term $I_i$:
  \be\label{sep3}
  \ba{rl}\vs1 \dy |I_i(x)|\leq &\dy|2-p| \int
(\mu+|\nabla v|^2)^{\frac{p-4}{2}}
{|\nabla v|^2 |\nabla\nabla v| |D^2 w|}\,dx\\\vs1
\leq &\dy |2-p|\int (\mu+|\nabla v|^2)^{\frac{p-2}{2}}
|\nabla\nabla v| |D^2 w|\,dx\\  \leq &\dy H(q')|2-p|
\big\|(\mu+|\nabla v|^2)^{\frac{p-2}{2}} \nabla\nabla v\big\|_q
  \|h\|_{q'}\,.\ea\ee
  Hence, from \eqref{sep1}--\eqref{sep3} we obtain
  $$\ba{rl}\vs1 \dy | \int (\mu+|\nabla v|^2)^{\frac{p-2}{2}}
  D^2_{x_ix_j} v \cdot h\, dx|\leq&\dy\!\!\!   H(q')
\big\|(\mu+|\nabla v|^2)^{\frac{p-2}{2}}\Delta v\big\|_q
  \|h\|_{q'}\\& \dy + 4 H(q')|2-p|
\big\|(\mu+|\nabla v|^2)^{\frac{p-2}{2}} \nabla\nabla v\big\|_q
  \|h\|_{q'},\ea$$
  which ensures that
  $$\ba{ll}\vs1 \dy\big\|(\mu+|\nabla v|^2)^{\frac{p-2}{2}} D^2_{x_ix_j} v\big\|_q
 \! =\!\!\sup_{h\in C_0^\infty(\O) \atop
|h|_{q'}=1}\!\!|((\mu+|\nabla v|^2)^{\frac{p-2}{2}} D^2_{x_ix_j} v,h)|\\ \dy
\hskip2.5cm \leq\! H(q')
\big\|(\mu+|\nabla v|^2)^{\frac{p-2}{2}}\Delta v\big\|_q\dy+4 H(q')|2-p|
\big\|(\mu+|\nabla v|^2)^{\frac{p-2}{2}}\nabla\nabla v\big\|_q. \ea$$
By the assumptions on $p$ and $q$, a straightforward computation gives estimate \eqref{ELL2}. \chiu
     \begin{lemma}\label{LL3}{\sl
Let $\mu>0$, $p\in (1,\infty)$ and $q\in (1,+\infty)$.  Assume that $(\mu+|\nabla v|^2)^\frac{p-3}{4}D^2 v\in L^{2q}(\R^n)$,
 $(\mu+|\nabla v|^2)^\frac{p-2}{2}D^3v\in L^{q}(\R^n)$ and  $(\mu+|\nabla v|^2)^\frac{p-2}{2}D^2 v\in L^{r}(\R^n)$, for some $r\in (1,+\infty)$.  Then
\be\label{ELL3}\ba{rl} \vs1\dy
\big\|(\mu+|\nabla v|^2)^{\frac{p-2}{2}}\nabla\nabla\nabla v\big\|_q
\leq&\dy \!\!\!H(q')\big\|(\mu+|\nabla v|^2)^{\frac{p-2}{2}}
D\Delta v\big\|_q\\ &\dy \hfill +4 |2-p|H(q')\| (\mu+|\nabla v|^2)^{\frac{p-3}{4}}
\nabla\nabla v\|_{2q}^2 \,.\ea\ee
}
 \end{lemma}
 \Pr Let $h$ and $w$ be as in Lemma\,\ref{perLL2}.  Multiplying $(\mu+|\nabla v|^2)^{\frac{p-2}{2}}D^3_{x_ix_jx_k}v $
  by $h$ and integrating by parts in $\R^n$, we obtain
$$
\ba{ll}\vs1 \dy \int(\mu+|\nabla v|^2)^{\frac{p-2}{2}}
D^3_{x_ix_jx_k} v \cdot h\, dx  =\dy\int (\mu+|\nabla v|^2)^{\frac{p-2}{2}}
D^3_{x_ix_jx_k} v \cdot \Delta w\, dx\\
=\vs1 \dy\!-\!\int\!(\mu+|\nabla v|^2)^{\frac{p-2}{2}}
D^2_{x_ix_j}v\!\cdot\! D_{x_k}\!\Delta w\,dx\!+\frac{2-p}{2}\!
\int\! (\mu+|\nabla v|^2)^{\frac{p-4}{2}}D^2_{x_ix_j}v\!\cdot\!\Delta wD_{x_k}|\nabla v|^2dx\\ \dy
\hfill :=\dy -\int(\mu+|\nabla v|^2)^{\frac{p-2}{2}}
D^2_{x_ix_j}v\cdot D_{x_k}\Delta w\,dx+I_0
.\ea$$
We then integrate three times more by parts, denoting at the $i-th$ integration, by $I_i$ the integral with the
derivatives of the term $a(\mu,v)$. Then using the previous identity we get
\be\label{sep4}\ba{ll}\vs1 \dy \!\!\int \!(\mu+|\nabla v|^2)^{\frac{p-2}{2}}
D^3_{x_ix_jx_k} v \!\cdot\! h dx \! =\!  \dy
\!\int\!(\mu+|\nabla v|^2)^{\frac{p-2}{2}}
D^3_{x_ix_jx_h}v\!\cdot\! D^2_{x_kx_h}w\,dx\!+\!\sum_{i=0}^1\! I_i
\\\hfill
\vs1
\dy =-\!\int (\mu+|\nabla v|^2)^{\frac{p-2}{2}}D^2_{x_ix_h}v\cdot D^3_{x_jx_kx_h}w\,dx+\sum_{i=0}^2 I_i\\
\hfill =\dy \!\int(\mu+|\nabla v|^2)^{\frac{p-2}{2}}D_{x_i}\Delta v\cdot D^2_{x_jx_k}w\,dx+\sum_{i=0}^3 I_i\,,\ea\ee
where
$$I_1:= \frac{p-2}{2}\int(\mu+|\nabla v|^2)^{\frac{p-4}{2}}
D^2_{x_ix_j}v\cdot D^2_{x_kx_h}w\, D_{x_h} |\nabla v|^2\,dx\,,$$
$$I_2:= \frac{2-p}{2}\!\int(\mu+|\nabla v|^2)^{\frac{p-4}{2}}
D^2_{x_ix_h}v\cdot D^2_{x_kx_h}w\, D_{x_j} |\nabla v|^2\,dx
\,,$$
$$ I_3:=\frac{p-2}{2}\int(\mu+|\nabla v|^2)^{\frac{p-4}{2}}
D^2_{x_ix_h}v\cdot D^2_{x_jx_k} w\, D_{x_h} |\nabla v|^2\,dx
\,.$$
By using estimate \eqref{Poi}, the first term on the right-hand side of \eqref{sep4} is estimated as
\be\label{sep5}
 |\int_\O (\mu+|\nabla v|^2)^{\frac{p-2}{2}}D_{x_i}\Delta v\cdot D^2_{x_jx_k}w\,dx|\leq H(q')
\big\|(\mu+|\nabla v|^2)^{\frac{p-2}{2}}D_{x_i}\Delta v\big\|_q
  \|h\|_{q'}\,.\ee
 The generic term $I_i$ can be estimated as follows
 \be\label{sep6}
 \ba{rl}\vs1 \dy |I_i(x)|\leq &\dy|2-p| \int
(\mu+|\nabla v|^2)^{\frac{p-4}{2}}
|\nabla v| |\nabla\nabla v|^2 |\nabla\nabla w|\,dx\\\vs1 \leq &\dy|2-p|
\int(\mu+|\nabla v|^2)^{\frac{p-3}{2}}
|\nabla\nabla v|^2 |\nabla\nabla w| \,dx\\  \leq &\dy
 |2-p|H(q')\|(\mu+|\nabla v|^2)^{\frac{p-3}{4}}
 \nabla\nabla v\|_{2q}^2
  \|h\|_{q'}\,.\ea\ee
  Hence, from \eqref{sep4}--\eqref{sep6}  we obtain
  $$\ba{ll}\vs1 \dy \left|\int(\mu+|\nabla v|^2)^{\frac{p-2}{2}}D^3_{x_ix_jx_k} v \cdot h\, dx\right|\leq   H(q')
\big\|(\mu+|\nabla v|^2)^{\frac{p-2}{2}}D_{x_i}\Delta v\big\|_q\|h\|_{q'}
\\ \dy \hfill  + 4|2-p|H(q')\|(\mu+|\nabla v|^2)^{\frac{p-3}{4}} \nabla\nabla v\|_{2q}^2  \|h\|_{q'},\ea$$
  which easily gives \eqref{ELL3}.
\chiu
\vskip0.1cm In our construction of the solution to system \eqref{stokes},
we will consider family of solutions to suitable approximating problems, and then pass to the limit in the different parameters,
one after the other.
Since the adopted convergence procedures will be similar for some parameters, we collect below some useful lemmas, that will be applied for these convergences.
\begin{lemma}\label{converg1}{\sl Let $r\in (1,+\infty)$, and $ \ov r\in [1,n)$.
Let $\psi^\nu$ be a third-order tensor and let $\{\psi^\nu\}$ be a bounded sequence in $L^r(E)$, $E\subseteq \R^n$ open set with the cone property.
 Let $\{v^\nu\}$ be a sequence of functions with $\{D^2 v^\nu\}$ bounded
 in $L^{\overline r}(E)$, and $\{\nabla v^\nu\}\in L^{\frac{n\ov r}{n-\ov r}}(E)$.
Then, for any $s\in [0,+\infty)$, there exist a subsequence of $\{\frac{\psi^\nu}{a_{s}(\mu, v^\nu)}\}$, a function $v$ and a third-order tensor $\psi$ such that
$$ \psi^\nu\rightharpoonup \psi\textrm{ in }  L^r(E)\,,$$
$$  \frac{\psi^\nu}{a_{s}(\mu, v^\nu)} \rightharpoonup  \frac{\psi}{a_{s}(\mu, v)} \textrm{ in }  L^r(E)\,,$$
where $a_s(\mu,\cdot)$ is defined in \eqref{amuq}.
}
\end{lemma}

\Pr The assumption of boundedness
of $\{\psi^\nu\}$ in $L^r(E)$ and
the bound $\frac{1}{a_s(\mu,
u^\nu)}\leq \mu^{-s}$ in $E$ imply
the existence of a tensor field
$\ov\psi$ of the third-order and of
a subsequence of
$\{\frac{\psi^\nu}{a_s(\mu,
v^\nu)}\}$, that we do not relabel,
such that \be\label{1extr}
\frac{\psi^\nu}{a_s(\mu,
v^\nu)}\rightharpoonup \ov\psi
\textrm{ in }  L^{r}(E). \ee Let us
work on this subsequence. The same
assumption of  boundedness of
$\{\psi^\nu\}$ in $L^r(E)$, implies
the existence of a tensor field
$\psi$ of the third-order and of
another subsequence, that we still
denote by $\{\psi^\nu\}$, such that
\be\label{2extr}
\psi^\nu\rightharpoonup
\psi\textrm{  in }  L^{r}(E). \ee
Let us consider the corresponding
subsequence of $\{a_{s}(\mu,
v^\nu)\}$, and let us fix an
arbitrary bounded open set
$K\subset \R^n$, such that $K\cap
E$ satisfies the cone property.
From the Rellich compactness
theorem, the $L^{\overline
r}$-uniform bound of $\{D^2v^\nu\}$
implies the existence, for any
$\overline{\overline r}\in
[1,\frac{n{\overline
r}}{n-{\overline r}})$, of a
subsequence of $\{\nabla v^\nu\}$
such that $\nabla v^\nu$ strongly
converges to $\nabla v$ in
$L^{\overline{\overline r}}(K\cap
E)$. This ensures that, along
another subsequence, as
$\nu\to\infty$, one has
\be\label{3extr} a_{s}(\mu, v^\nu)
\to a_{s}(\mu, v), \ \mbox{ a.e. in
} K\cap E\,.\ee Therefore, by
extracting three times from the
 extract in \eqref{1extr}, we have obtained a subsequence
 $\{\frac{\psi^\nu}{a_{s}(\mu, v^\nu)}\}$ such that both \eqref{2extr} and \eqref{3extr} are satisfied.
Let $\vp$ be a vector function in $L^{r'}(K\cap E)$, and set
\be\label{weakter}\ba{ll}\dy\vs1
 (\frac{\psi^\nu}{a_{s}(\mu, v^\nu)}- \frac{\psi}{a_{s}(\mu, v)}, \vp)=&
 \dy (\frac{\psi^\nu}{a_{s}(\mu, v^\nu)}-
 \frac{\psi^\nu}{a_{s}(\mu, v)}, \vp)\\ &\dy +
 \, (\frac{\psi^\nu}{a_{s}(\mu, v)}- \frac{\psi}{a_{s}(\mu, v)}, \vp):=I_1^\nu+I_2^\nu\,.\ea\ee
Clearly  \be\label{onge}| \frac{1}{a_{s}(\mu, v^\nu)}- \frac{1}{a_{s}(\mu, v)}| \leq 2\, \mu^{-s}\,, \mbox{ a.e. in } K\cap E\,.\ee
By applying H\o lder's inequality,
 $$|I_1^\nu|\leq \|\psi^\nu\|_r \|(\frac{1}{a_{s}(\mu, v^\nu)}-
 \frac{1}{a_{s}(\mu, v)})\,\vp\|_{r'}\,.$$
Thanks to estimate \eqref{onge}, we can apply the
 Lebesgue dominated convergence theorem, which, together with the $L^r$-bound on $\{\psi^\nu\}$, ensures that the
 sequence $I_1^\nu$ goes to zero as $\nu\to\infty$. The sequence $I_2^\nu$ goes to zero by the weak convergence of
 $\psi^\nu$ to $\psi$ in $L^r(K\cap E)$, observing that $\frac{\vp}{a_{s}(\mu, v)}\in L^{r'}(K\cap E)$. Therefore $\frac{\psi^\nu}{a_{s}(\mu, u^\nu)}$ weakly converges to $\frac{\psi}{a_{s}(\mu, v)}$ in $L^r(K\cap E)$.
 On the other hand, the function $\ov\psi$ is  weak limit also in $L^r(K\cap E)$.
 By uniqueness of the weak limit we get $\ov\psi=\frac{\psi}{a_{s}(\mu, v)}$ a.e. in $K\cap E$. The arbitrariness of $K$ gives the result.
 \chiu
\begin{lemma}\label{converg2}{\sl Let the assumptions of Lemma \ref{converg1} be satisfied,
and let $h\in L^\infty(E)$, with compact support.
 Then, for any $s \in [0,+\infty)$, $t\in [1,+\infty)$, there exist a subsequence of \linebreak $\{A^t(\mu, v^\nu)\cdot J_\eta(\frac{\psi^\nu\, h}{a_{s}(\mu, v^\nu)}\}$, a function $v$ and a third-order tensor $\psi$, such that
$$ A^t(\mu, v^\nu)\cdot J_\eta(\frac{\psi^\nu\, h}{a_{s}(\mu, v^\nu)})
  \rightharpoonup
 A^t(\mu, v)\cdot J_\eta(\frac{\psi\, h}{a_{s}(\mu, v)})
 \textrm{  in }  L^r(E)
 \,,$$
where $A^t(\mu, \cdot )$ is defined in \eqref{fot} and $J_\eta$ denotes the Friedrich's mollifier.\par
The same result holds if we replace the sequence $\{A^t_{ijhk}(\mu, v^\nu) J_\eta(\frac{\psi^\nu_{jhk}\, h}{a_{s}(\mu, v^\nu)}\}$ with the sequence
$\{A^t_{jihk}(\mu, v^\nu)J_\eta(\frac{\psi^\nu_{jhk}\, h}{a_{s}(\mu, v^\nu)}\}$.

}
\end{lemma}
\Pr From Lemma \ref{converg1}, there exist a subsequence of $\{\frac{\psi^\nu}{a_{s}(\mu, v^\nu)}\}$, a third-order tensor $\psi$ and
a function $v$ such that,
for any $s\in [0,+\infty)$,
\be\label{extr1}
  \frac{\psi^\nu\, h}{a_{s}(\mu, v^\nu)} \rightharpoonup  \frac{\psi\, h}{a_{s}(\mu, v)} \textrm{  in }  L^r(E)\,.\ee
Let us consider the corresponding
subsequence of $\{A^t(\mu, v^\nu)\}$.
Denote by $K_\eta$ the compact support of $J_\eta(\frac{\psi^\nu\, h}{a_{s}(\mu, v^\nu)})$, which, for
$\eta$ sufficiently small, is contained in $E$. From Rellich's compactness theorem, for any $\overline{\overline r}\in [1,\frac{n{\overline r}}{n-{\overline r}})$ there exists a subsequence of $\{\nabla v^\nu\}$ such that
$$\nabla v^\nu \to \nabla v\ \mbox{ in }\ L^{\overline{\overline r}}(K_\eta)\,,$$
hence, taking another subsequence,
 \be\label{extr3}A^t(\mu, v^\nu)\to A^t(\mu, v) , \ \mbox{ a.e. in } \ K_\eta\,.\ee
 Let us fix the last subsequence
 of $\{ A^t(\mu, v^\nu)\cdot J_\eta(\frac{\psi^\nu\, h}{a_{s}(\mu, v^\nu)})\}$. Let $\vp\in L^{r'}(E)$ and let us set
\be\label{forpr}\ba{ll}\vs1\dy
\big(\,( A^t(\mu, v^\nu)-
A^t(\mu, v)\,)\cdot J_\eta(\frac{\psi^\nu h}{a_s(\mu, v^\nu)}), \vp\big)\\\hskip3cm
\dy +\big (A^t(\mu, v)\cdot J_\eta\big(\frac{\psi^\nu h}{a_s(\mu, v^\nu)}-\frac{\psi h}{a_s(\mu, v)}\big), \vp\big)
:=I_1^\nu+I_2^\nu\,.
\ea\ee
We have $$\big\|J_\eta(\frac{\psi^\nu\, h}{a_s(\mu, v^\nu)})\big\|_r\leq \big\|
 \frac{\psi^\nu\, h}{a_s(\mu, v^\nu)}\big\|_r\leq C(\mu) \|\psi^\nu\|_r\,, \forall \nu\in \N$$
and $$\|A^t(\mu, v^\nu) )-A^t(\mu, v) ) \|_\infty\leq 2 \mu^{1-t}\,, \forall \nu\in \N.$$
Since convergence \eqref{extr3} holds, by applying the H\o lder inequality and then the Lebesgue dominated convergence theorem,
 it follows that $I_1^\nu\to 0$.
As far as the sequence of integrals $I_2^\nu$ is concerned, using Fubini's theorem we get
$$I_2^\nu= \big(\big(\frac{\psi^\nu \,h}{a_s(\mu, v^\nu)}-\frac{\psi\,h}{a_s(\mu, v)}\big), J_\eta(A^{t}(\mu, v)
\cdot \vp\big)\big)\,, $$
which tends to zero from the weak convergence \eqref{extr1},
since $J_\eta(A^t(\mu, v) \cdot \vp\big)\in L^{p'}(E) $, uniformly in $\nu$.
\chiu
\begin{lemma}\label{converg3}{\sl Ler $r\in (1,+\infty)$, and $ {\overline r}\in [1,n)$.
Let $\psi^\nu$ be a third-order tensor and let $\{\psi^\nu\}$ be a bounded sequence in $L^r(E)$, $E\subseteq \R^n$ open set with the cone property. Let $\{v^\nu\}$ be a sequence of functions such that $\{D^2 v^\nu\}$ is  bounded
 in $L^{\overline r}(E)$ and $\{\nabla v^\nu\}\in L^{\frac{n\ov r}{n-\ov r}}(E)$.
Then, for any $s\in [0,+\infty)$, there exist a subsequence of $\{A^s_{ijhk}(\mu, v^\nu)\psi^\nu_{jhk}\}$, a third-order tensor  $\psi$ and a function $v$ such that
$$ \psi^\nu \rightharpoonup \psi \textrm{  in }  L^r(E)\,,$$
$$ A^s(\mu, v^\nu)\cdot \psi^\nu \rightharpoonup A^s(\mu, v)\cdot \psi \textrm{ in }  L^r(E)\,,$$
where $A^s(\mu,\cdot)$ is defined in \eqref{fot}.
\par The same result holds if we replace the sequence $\{A^s_{ijhk}(\mu, v^\nu)\psi^\nu_{jhk}\}$ with the sequence
$\{A^s_{jihk}(\mu, v^\nu)\psi^\nu_{jhk}\}$.
}
\end{lemma}
\Pr The proof is the same as in Lemma \ref{converg1}, with the sequence $\{a_s(\mu,v^\nu)\}$ replaced by the sequence
$\{ A^s(\mu, v^\nu)\}$, since they share the same properties of boundedness and almost everywhere convergence on
compact sets, along subsequences.\chiu

\begin{lemma}\label{GNT}{\sl Let $\O$ be a Lipschitz  bounded domain. Assume that $v=0$ on $\pa\O$ and
$D^2 v\in L^{\frac{2n}{n+2}}(\O)$. Then
\be\label{nd}
\|\nabla v\|\leq C_s\|D^2 v\|_{\frac{2n}{n+2}}\,,\ee
where $C_s$ is the Sobolev constant.}
\end{lemma}
\Pr Since the boundedness of $\O$ implies that $\nabla v\in L^2(\O)$, by Sobolev's embeddings
\be\label{sob}
\|v\|_{\frac{2n}{n-2}}\leq C_s\|\nabla v\|\,,\ee
where  $C_s$ is the Sobolev's constant, which is independent on $\O$.
Integrating by parts, by using that $v\in W_0^{1,2}(\O)$ and then employing estimate \eqref{sob}, we get
$$\|\nabla v\|^2=-(\Delta v, v)\leq \|D^2 v\|_{\frac{2n}{n+2}}\|v\|_{\frac{2n}{n-2}}\leq C_s \|D^2 v\|_{\frac{2n}{n+2}}\|\nabla v\|\,,$$
and the result is proved. \chiu

\section{\large An auxiliary problem in a smooth convex and bounded domain $\O$}\label{secconvex}
Let $\Omega\subset \R^n$ be a bounded convex domain, whose boundary $\pa\O$ is $C^2$-smooth.
Let us
 introduce the following auxiliary
system \be\label{stokesSo}\ba{rl}\dy \ve\Delta u+\frac{\Delta u}{a_{2-p}(\mu,\!u)}
+(p\!-\!2)\frac{\nabla u\!\otimes\!\nabla u}{a_{\frac{4-p}{2}}(\mu,\!u)}\!\cdot\! J_\eta(\frac{\nabla\nabla u\,\chi^\rho}{a(\mu,\!u)})\!\!\!\!\!&\dy=\!
\frac{\nabla \Pi(u,\chi^\rho)}{a(\mu,\!u)}\!+\! \frac{f}{a(\mu,\!u)},\mbox{in\,}\O,\\\dy u&=0,\ \mbox{ on } \pa\O\,,
\ea\ee
where $\mu>0$. Here we assume that the support of the cut-off function $\chi^\rho$, introduced in sec. \ref{notations}, is included in $\O$, and, recalling the definition of
 $\Pi(u,\chi^\rho)$  given by \eqref{pi}, a suitable extension to $\R^n$ of $\nabla u$ is assumed.
 Since $\nabla\Pi(u,\chi^\rho)$ is expressed by means of a Calder\'on-Zygmund singular kernel, its $L^2$-norm can be estimated as
\be\label{CZ2}
\|\nabla \Pi(u,\chi^\rho)\|\leq (2-p)H(2)
\big\|\frac{\nabla u_j}{a_1(\mu, u)} \cdot \big(\nabla u\cdot J_\eta(\frac{D_{y_j}\nabla u\,\chi^\rho}{a(\mu, u))})\big)\|
.\ee
Throughout this section the norm $\|\cdot\|_{r}$ will be
always the $L^r$-norm on the domain $\O$.
Note that we aim at estimates which are uniform with respect to the size of the domain $\O$, in order to apply, in the next section, the method of expanding domains. \begin{prop}\label{Vteo2}{\sl
Let $\ov M(2)>0$,  $\ov M(2)$ defined in \eqref{cseg}. Let $f\in L^2(\O)$. There exists a solution $u\in W_0^{1,2}(\O)\cap
W^{2,2}(\O)$ of problem \eqref{stokesSo}, with
\begin{equation}\label{main1}\|\, D^2 u \,\|
\leq C(\ve)\,\|f\|\,,\end{equation}
\begin{equation}\label{main1p}\|\, \Pi(u,\chi^\rho)\,\|_{\frac{2n}{n-2}}
\leq c\,\|f\|\,,\end{equation}
\begin{equation}\label{maingp}\|\,\nabla\Pi(u,\chi^\rho) \,\|\leq c\,\|f\|\,.\end{equation}
}
\end{prop}
\Pr
 Let $\{a_j\}$ be the eigenfunctions of
the Laplace operator $-\Delta$, and denote by
$\lambda_j$ the corresponding (positive)
eigenvalues:
$$\ba{ll}\vs1 \dy-\Delta a_j\!\!\!\!&\dy=\lambda_j a_j,\ \mbox{ in }
\O,\\\dy \hskip0.6cm a_j\!\!\!\!&=\dy0,\hskip0.5cm  \mbox{ on }
\pa \O\,. \ea$$ Recall that $\{a_j\}$ is a
complete system in
$W_0^{1,2}(\O)\cap W^{2,2}(\O)$,
orthonormal in $L^{2}(\O)$.   We
construct the Faedo-Galerkin
approximations related to system
\eqref{stokesSo}, such that, for each
$k\in \N$, \be\label{G1}u^k(x)=
\sum_{j=1}^k c_{jk}\,a_j(x)\,,\ee
and
\be\label{sis1}\ba{ll}\vs1\dy Q_j(c_k):=&\dy\hskip-0.3cm \ve (\Delta
u^k,a_j)+(\frac{\Delta u^k}{a_{2-p}(\mu,u^k)},a_j)
+(p-2)(A^{\frac{4-p}{2}}(\mu, u^k) \cdot J_\eta(\frac{\nabla\nabla u^k\,\chi^\rho}{a(\mu, u^k)}), a_j)
\\
&\dy -(\frac{\nabla \Pi (u^k, \chi^\rho)}{a(\mu, u^k)}, a_j)- (\frac{f}{a(\mu, u^k)},a_j)=0
\,,\ea\ee
for any $j=1,\cdots,k$, with
$c_k=(c_{1k},\cdots, c_{kk})$. This is a system of $k$ equations in the
unknowns $c_{jk}$. Set \be\label{siss}\ba{ll}\vs1 P_j(c_k):=&\dy\hskip-0.3cm- \ve (\Delta
u^k,\lambda_j a_j)-(\frac{\Delta u^k}{a_{2-p}(\mu,u^k)},\lambda_ja_j)\\ \vs1&\dy+(2-p)(A^\frac{4-p}{2}(\mu, u^k)
\cdot J_\eta(\frac{\nabla\nabla u^k\,\chi^\rho}{a(\mu, u^k)}), \lambda_ja_j)
\\ &\dy +(\frac{\nabla \Pi (u^k ,\chi^\rho)}{a(\mu, u^k)}, \lambda_ja_j)+(\frac{f}{a(\mu, u^k)},\lambda_ja_j)\ea\ee $j=1,\cdots,
k$, and observe that \be\label{position}Q_j(c_k)=-\,\frac{1}
{\lambda_j} P_j(c_k).\ee Let us verify that $P(c_k)\cdot c_k\geq 0$,
for a suitable $c_k$.  Thanks to our choice of the basis, taking the scalar product  of $P(c_k)$ with $c_k$, the first two terms on the right-hand side of \eqref{siss}  give
\be\label{s1} \ve \|\Delta u^k\|^2+\big\|\frac{\Delta u^k}{a(\mu,u^k)}\big\|^2\,.\ee
 Let us estimate the term coming from the product of $c_{jk}$ with the third term in \eqref{siss}.  Firstly we observe that, since $\frac{\nabla\nabla u^k\,\chi^\rho}{a(\mu, u^k)}\in L^2(\O)$, then
\be\label{geta}
\big\|J_\eta(\frac{\nabla\nabla u^k\,\chi^\rho}{a(\mu, u^k)})\big\|
\leq \big\|\frac{\nabla\nabla u^k\,\chi^\rho}{a(\mu, u^k)}\big\|\,.\ee Hence, from
 $|\frac{\nabla u^k\otimes \nabla u^k}{a_{1}(\mu, u^k)}|\leq 1$ and  $|\chi^\rho| \leq 1$, by applying
 Lemma \ref{LL1} (note that $\ov M(2)>0$ implies $p>\frac 32$),
we  get
$$\ba{ll}\dy\vs1\sum_{j=1}^k (\frac{\nabla u^k\otimes \nabla u^k}{a_{1}(\mu, u^k)}
\cdot J_\eta(\frac{\nabla\nabla u^k\,\chi^\rho}{a(\mu, u^k)}), \frac{c_{jk}\lambda_ja_j}{a(\mu, u^k)})=
-(\frac{\nabla u^k\otimes \nabla u^k}{a_{1}(\mu, u^k)}
\cdot J_\eta(\frac{\nabla\nabla u^k\,\chi^\rho}{a(\mu, u^k)}), \frac{\Delta u^k}{a(\mu, u^k)})\\ \dy\hfill
\geq
-\big\|\frac{\nabla\nabla u^k\,\chi^\rho}{a(\mu, u^k)} \big\| \big\|\frac{\Delta u^k}{a(\mu, u^k)} \big\|
\geq \frac{1}{3-2p} \big\|\frac{\Delta u^k}{a(\mu, u^k)} \big\|^2\,.
 \ea$$
The product of $c_{jk}$ with the fourth term in \eqref{siss} can be estimated similarly. As remarked in \eqref{CZ2},
$\nabla \Pi(u^k, \chi^\rho)$ is
expressed by means of a Calder\'on-Zygmund singular kernel. Hence estimate \eqref{CZ2}
and inequality \eqref{geta}, easily imply
\be\label{press}\|\nabla \Pi (u^k, \chi^\rho)\|\leq (2-p)H(2)
\big\|\frac{\nabla\nabla u^k\,\chi^\rho}{a(\mu, u^k)} \big\| \,.\ee
Then, arguing as for the previous term, we have
$$\ba{ll}\vs1\dy\sum_{j=1}^k (\frac{\nabla \Pi (u^k, \chi^\rho)}{a(\mu, u^k)}, c_{jk}\lambda_ja_j)=-
(\nabla \Pi (u^k, \chi^\rho), \frac{\Delta u^k}{a(\mu, u^k)})\geq
\frac{2-p}{(3-2p)}H(2) \big\|\frac{\Delta u^k}{a(\mu, u^k)} \big\|^2\,.\ea$$
For the last term in \eqref{siss}, we
just apply H\o lder's inequality and get
$$ \ba{ll}\dy |\sum_{j=1}^k (\frac{f}{a(\mu, u^k)},c_{jk}\lambda_ja_j)|= |-(f,\frac{\Delta u^k}{a(\mu, u^k)})|
\leq \|f\| \big\| \frac{\Delta u^k}{a(\mu, u^k)}\big\|\,.\ea$$
 Therefore, since $\ov M(2)>0$,
$$\ba{rl}\vs1\dy  P(c_k)\cdot c_k
\dy \geq&\dy \ve \|\Delta u^k\|^2+
 \big\| \frac{\Delta u^k}{a(\mu, u^k)} \big\|\,\big(
 \big\| \frac{\Delta u^k}{a(\mu, u^k)} \big\| - \frac{2-p}{2p-3} \big\| \frac{\Delta u^k}{a(\mu, u^k)} \big\|\\& \dy
 -\frac{2-p}{(2p-3)}\,H(2) \big\| \frac{\Delta u^k}{a(\mu, u^k)} \big\|- \|f\| \big)\geq 0\,,\ea$$ provided
that  \be\label{onD}
\big\| \frac{\Delta u^k}{a(\mu, u^k)} \big\|\geq \frac{2p-3}{\ov M(2)}\|\,f\|\,.\ee
Since the vector valued function $P$ is a continuous function, in order to apply Lemma \ref{Lionslemma}
it is sufficient to show that there exists a  $R$ such that, for $|c_k|=R$,
 condition \eqref{onD} is satisfied.
 Denoting by
 $\lambda_1$ the smallest eigenvalue,  we observe that $$-\sum_{j=1}^k (\Delta u^k, c_{jk}\lambda_j a_j)= \sum_{j=1}^k c_{jk}^2 \lambda_j^2\geq \lambda_1^2 |c_k|^2.$$
 Further
 $$\ba{ll}\dy a(\mu, u^k(x))\leq (\mu+(\sum_{j=1}^k |c_{jk}| |\nabla a_j(x)|)^2)^{\frac{2-p}{2}}
\leq  (\mu+(\sum_{j=1}^k |c_{jk}|\, \max_{\ov\O} |\nabla a_j|)^2)^{\frac{2-p}{2}}\\
\dy \hfill \leq  (\mu+ (\sum_{j=1}^k |c_{jk}|\, \max_{1\leq j\leq k} |\nabla a_j|_{C(\ov\O)})^2)^{\frac{2-p}{2}}\leq
  (\mu+ k^2 |c_k|^2\, (\max_{1\leq j\leq k} |\nabla a_j|_{C(\ov\O)})^2)^{\frac{2-p}{2}},\ea$$
where we have used the regularity of $\{a_j\}$. Therefore
one gets  $$\big\| \frac{\Delta u^k}{a(\mu, u^k)} \big\|\geq \dy\frac{\lambda_1\,|c_k|}{
\mu^\frac{2-p}{2}+(k\, \max_{1\leq j\leq k} |\nabla a_j|_{C(\ov\O)})^2)^{\frac{2-p}{2}} |c_k|^{2-p}} \,,
$$ which tends to infinity as $|c_k|$ increases. Hence there exists a $R>0$ such that for any $|c_k|\geq R$ estimate \eqref{onD} holds.
 By Lemma
\ref{Lionslemma}, this proves the
existence of a solution $\widetilde c_k$ of the
algebraic system $P_j(c_k)=0$ with $|\widetilde c_k|\leq R$. Recalling
inequality \eqref{position}, we also obtain a
solution of the kind \eqref{G1} for
system \eqref{sis1}. \par The above arguments also give an estimate on the second derivatives, uniformly
 with respect to $k\in\N$.
Actually, since $\widetilde c_k$ solves \rf{siss}, that
is $P(\widetilde c_k)=0$, we
get \be\label{aip}\ba{rl}\vs1\dy \ve \|\Delta u^k\|^2+
 \big\| \frac{\Delta u^k}{a(\mu, u^k)} \big\|^2 \leq  \frac{2-p}{2p-3}(1+H(2)) \big\| \frac{\Delta u^k}{a(\mu, u^k)} \big\|^2+
  \big\| \frac{\Delta u^k}{a(\mu, u^k)} \big\| \|f\|\,.\ea\ee Hence,
   \be\label{onDD}
   \big\| \frac{\Delta u^k}{a(\mu, u^k)} \big\|\leq \frac{2p-3}{\ov M(2)} \|f\|\,,\ee
    \be\label{onL}  \|\Delta u^k\|\leq \frac{c}{\sqrt \ve}\|f\|\,, \ee
   and, due to Lemma \ref{LL1},
     \be\label{onDn}
   \big\| \frac{\nabla\nabla u^k}{a(\mu, u^k)} \big\|\leq \frac{1}{\ov M(2)} \|f\|\,,\ee   uniformly in $k$.
   The next step is to prove the existence of a subsequence and of a limit function $u$, as $k\to\infty$,
   which satisfies the same bounds of the sequence $\{u^k\}$, and satisfies
   \be\label{limuk}
   (\ve\Delta u+\frac{\Delta u}{a_{2-p}(\mu,u)}
+(p-2) A^{\frac{4-p}{2}}(\mu, u)\cdot J_\eta(\frac{\nabla\nabla u\,\chi^\rho}{a(\mu, u)})-
\frac{\nabla \Pi(u, \chi^\rho)}{a(\mu,u)}- \frac{f}{a(\mu,u)},\vp)=0\,,\ee
for all $\vp\in C_0^\infty(\O)$.
   Thanks to  \eqref{onL} and Lemma\,\ref{LL1} with $p=2$, the sequence $\{D^2u^k\}$ is uniformly bounded in $L^2(\O)$. Hence we get that the sequence $\{u^k\}$ is uniformly bounded in $W_0^{1,2}(\O)\cap W^{2,2}(\O)$. Therefore we are in the hypotheses
      of Lemma \ref{converg1}, with  $r=\ov r=2$, $E=\O$, and in the hypotheses of
   Lemma\,\ref{converg2}, with  $h=\chi^\rho$. The quoted lemmas imply the existence of a field $u$
and of a subsequence, that we still denote by $\{u^k\}$, such that, for any $\vp\in L^2(\O)$,
$$ \lim_k\,(\Delta u^k-\Delta u, \vp)\,=0\,,$$
$$
\lim_k (\frac{\Delta u^k}{a_{2-p}(\mu, u^k)}- \frac{\Delta u}{a_{2-p}(\mu, u)}, \vp)=0\,,$$
and
\be\label{der2}\lim_k  \big(\,A^\frac{4-p}{2}(\mu, u^k)\cdot J_\eta(\frac{\nabla\nabla u^k\,\chi^\rho}{a(\mu, u^k)})-
A^\frac{4-p}{2}(\mu, u)\cdot J_\eta(\frac{\nabla\nabla u\,\chi^\rho}{a(\mu, u)}), \vp\big)=0\,.\ee
The boundedness of $\O$ and the lower semi-continuity ensure that
\eqref{main1} holds and, clearly,
$$\|u\|_{2,2}\leq c
\|f\|,$$
where the constant $c$ depends on $\ve$ and the size of $\O$.
 Finally, let us consider the corresponding subsequence of $\{\nabla \Pi^k\}$, where $\Pi^k:=\Pi(u^k, \chi^\rho)$.
From the Hardy-Littlewood-Sobolev
theorem (see \cite{stein}, Chap. 5,
Theorem 1) and estimate
\eqref{onDn}, the sequence
 $\{\Pi^k\}$ is bounded in $L^r(\O)$, $r\in [1,\frac{2n}{n-2})$. Moreover, by
using estimates \eqref{press} and \eqref{onDn}, it follows that the sequence $\{\nabla\Pi^k\}$ is bounded in $L^2(\O)$.
Therefore there exist a limit $\Pi$ and a subsequence such that
$$\Pi^k\rightarrow \Pi \textrm{ in } L^{r}(\O),\ \forall r\in[1,\frac{2n}{n-2}), $$
$$\Pi^k\rightharpoonup \Pi \textrm{ in } W^{1,2}(\O), $$
and, by applying  Lemma\,\ref{converg1}, for any $\vp\in L^2(\O)$,
$$(\frac{\nabla\Pi^k}{a(\mu, u^k)}, \vp) \to(\frac{\nabla\Pi}{a(\mu, u)}, \vp)\,.$$ We want to show that
$\Pi=\Pi(u, \chi^\rho)\, $.
It is enough to show that there exists a subsequence such that, for any $\vp\in C_0^\infty(\R^n)$,
$$\lim_{k\to \infty}\, (\Pi^k,\vp)=(\Pi(u, \chi^\rho),\vp)\,, $$
and the result will follow by uniqueness of the limit.
 We set
$$g(y)= \int_{\R^n}\nabla_x {\mathcal E} (x-y)\, \vp(x)\,dx\,.$$Then
$$\ba{rl}\vs1\dy
(\Pi^k,\vp)\dy=&\dy (p-2) \int_{\R^n}\frac{\nabla u_i^k   (y)D_{y_j}u_h^k(y) }{a_1(\mu, u^k(y))} J_\eta(\frac{D^2_{y_iy_j} u_h^k(y)\,\chi^\rho(y)}{a(\mu, u^k(y))}\big)\cdot g(y)\,dy\,.
\ea$$
The convergence to $(\Pi,\vp)$ of a subsequence follows from Lemma\,\ref{converg2},
observing that, from the Hardy-Littlewood-Sobolev theorem,
$g\in L^2(\R^n)$.
\par Summarizing, we have found a subsequence such that \eqref{limuk} holds, which ensures
that $u$ satisfies \eqref{stokesSo} a.e. in $\O$.
 The proof is therefore
completed.\chiu \vskip0.1cm
Let us
 introduce the following
   linear elliptic
system \be\label{stokesS2}\ba{rl}\vs1 \dy \ve\Delta w\!+\!\frac{\Delta w}{a_{2\!-\!p}(\mu,\!u,\!\delta)}
\!\!\!\!&\dy\!=\!(2\!-\!p) \frac{\nabla u\!\otimes\!\nabla u}{a_{\frac{4-p}{2}}(\mu,\!u)}\!\cdot\!J_\eta(\frac{\nabla\nabla u\chi^\rho}{a(\mu,\!u)})\!+\!\frac{\nabla \Pi(u,\!\chi^\rho)}{a(\mu,\!u)}\!+\!\frac{f}{a(\mu,\!u)},\mbox{in\,}\O,\\
w&\dy=0\,, \ \mbox{ on }\partial \O\,,\ea\ee where  $u$ is the solution of system
\eqref{stokesSo} given in Proposition \ref{Vteo2}, $a_{2-p}(\mu, u, \delta)$ is defined in \eqref{amudelta}, and $\Pi(u, \chi^\rho)$ is defined in \eqref{pi}.
The key tool in
proving  the following result on $u$ is to find suitable regularity properties for solutions of the above linear elliptic system are
\begin{prop}\label{proplq}{\sl
Let $\ov M(2)>0$. Let $q \in (n, +\infty)$, $q_1=\frac{np}{n+p}$ and assume $f\in L^q(\O)\cap L^{q_1}(\O)$.
Then, there exist constants $C$, independent of $\O$, such that the solution $u$ of Proposition\,\ref{Vteo2} satisfies the following estimates:
\be\label{passeu1}\|D^2 u\|_q \leq C(\mu,\eta,\rho, \ve)\, (\|f\|_{q}+ \|f\|) \,,\ee
\be\label{passeu2}\|D^2 u\|_{q_1} \leq C(\mu,\rho,\ve)\, (\|f\|_{q_1}+ \|f\|) \,.\ee
 }
\end{prop}
\Pr The proof will be achieved in three steps. Firstly we show that the
second derivatives of the solution of system \eqref{stokesS2} satisfy $L^q(\O)\cap L^{q_1}(\O)$ estimates,
uniformly in $\delta$.
 Then we show that, in the limit as $\delta$ tends to zero, the sequence of
solutions $\{w^\delta\}$ tends to a limit function $w$. In the final step we will prove that $w$ coincide with the solution $u$ of Proposition \ref{Vteo2}. \vskip0.2cm\noindent
{\it\underline{Step I:  $D^2 w\in L^q(\O)\cap L^{q_1}(\O)$, uniformly in $\delta>0$. }} -
We show that there exist  constants $C$, independent of $\delta$ and $meas\,\O$, such that
\be\label{passe1}\|D^2 w\|_q \leq C(\mu,\eta,\rho, \ve)\, (\|f\|_{q}+ \|f\|) \,,\ee
\be\label{passe2}\|D^2 w\|_{q_1} \leq C(\mu,\rho,\ve)\, (\|f\|_{q_1}+ \|f\|) \,.\ee
System \eqref{stokesS2} is a linear elliptic system
with $C^{\infty}(\R^n)$ coefficients. Observing that, from Proposition \ref{Vteo2}, the right-hand side of \eqref{stokesS2} belongs
to $L^2(\O)$ (see estimates \eqref{geta}, \eqref{press} and \eqref{onDn}),
 by classical existence and regularity results, one has the existence
of  $w\in W^{2,2}(\O)$ with
\be\label{con1}
\|D^2w\|\leq C(\mu, \eta, \ve) \|f\|\,.\ee
Firstly we estimate the $L^q$-norm of the right-hand side of \eqref{stokesS2}.
We get $$ \big\|A^\frac{4-p}{2}(\mu, u)\cdot J_\eta(\frac{\nabla\nabla u\,\chi^\rho}{a(\mu, u)})\big\|_q\leq
 \frac{1}{\mu}{\atop^{\frac{2-p}{2}}} \big\|  J_\eta(\frac{\nabla\nabla u\,\chi^\rho}{a(\mu, u)})\big\|_q\leq c(\mu,\eta,\rho)
 \big\| \nabla\nabla u\,\chi^\rho\big\|\,.$$
From \eqref{pi} we can apply to
$\nabla\Pi(u, \chi^\rho)$ the
Calder\'on-Zygmund theorem. Hence
$$\big\|\frac{\nabla \Pi(u, \chi^\rho)}{a(\mu,u)}\big\|_q\leq  \frac{(2-p)}{\mu^\frac{2-p}{2}}H (q)
 \big\|  J_\eta(\frac{\nabla\nabla u\,\chi^\rho}{a(\mu, u)})\big\|_q\leq c(\mu,\eta,\rho)H(q)(2-p) \big\| \nabla\nabla u\,\chi^\rho\big\|\,.$$
Finally $$\| \frac{f}{a(\mu,u)}\|_q\leq \frac{1}{\mu}{\atop^{\frac{2-p}{2}}} \|f\|_q\,.$$
By applying Lemma \ref{lemmaGM} we find that  the second derivatives of $w$ belong to $L^q(\O)$ and satisfy
$$\|D^2 w\|_q \leq C(\mu,\eta,\rho,\ve,\delta)\, (\|f\|_q+\|D^2 w\|)\,.$$
On the other, the knowledge of $D^2w\in L^q(\O)$, due to the
previous arguments, enables us to recover estimate \eqref{passe1}, which is
uniform in $\delta$, but not in the other parameters.
In order to show this,  we multiply equation \eqref{stokesS2} by $\Delta w |\Delta w|^{q-2}$, and integrate over $\O$. By H\o lder' s inequality
 and the previous estimates on the $L^{q}$-norms of the right-hand side of \eqref{stokesS2} we get
 $$\ba{rl}\dy\vs1  \ve \|\Delta w\|_{q}^{q}\!\leq\!\!\!&\dy   \!\!
(2-p) \big\|A^\frac{4-p}{2}(\mu, u)\cdot J_\eta(\frac{\nabla\nabla u\,\chi^\rho}{a(\mu, u)})\big\|_{q}
 \|\Delta w\|_{q}^{q-1}\\ \vs1 &\dy +\big\|\frac{\nabla \Pi(u,\chi^\rho)}{a(\mu,u)}\big\|_{q}  \|\Delta w\|_{q}^{q-1}\! +\dy
 \frac{1}{\mu}{\atop^{\frac{2-p}{2}}} \|f\|_{q} \|\Delta w\|_{q}^{q-1}\\
 \dy\! \leq\!\!\!&\dy\! \!\! \left(c(\mu,\eta,\rho)(2-p)(1+H(q))
 \big\| \nabla\nabla u\,\chi^\rho\big\|+ c(\mu)  \|f\|_{q}\right)\! \|\Delta w\|_{q}^{q-1},
    \ea$$
    which, from \eqref{main1}, and Lemma\,\ref{LL2} with $p=2$, immediately gives \eqref{passe1}.  In
order to get estimate \eqref{passe2}, at first we estimate the $L^{q_1}$-norm of the right-hand side of \eqref{stokesS2}.
Note that, since  $p\in (\frac 32, 2)$, then $q_1\in (1,2)$.
 Since $D^2u\in L^2(\O)$ and, by the assumption, $supp\, \chi^\rho\subset \O$, then
 $\nabla\nabla u\,\chi^\rho\in L^{q_1}(\R^n)$.
  Therefore
 $$\ba{ll}\dy  \big\|A^\frac{4-p}{2}(\mu, u)\cdot J_\eta(\frac{\nabla\nabla u\,\chi^\rho}{a(\mu, u)})\big\|_{q_1}&\dy
\leq
 \frac{1}{\mu}{\atop^{\frac{2-p}{2}}} \big\|  J_\eta(\frac{\nabla\nabla u\,\chi^\rho}{a(\mu, u)})\big\|_{q_1}\leq \frac{1}{\mu}{\atop^{\frac{2-p}{2}}}
 \big\| \nabla\nabla u\,\chi^\rho\big\|_{q_1}\\ &\dy \leq  c(\mu,\rho)
 \big\| \nabla\nabla u\,\chi^\rho\big\|\,,\ea$$
 and
$$\big\|\frac{\nabla \Pi(u, \chi^\rho)}{a(\mu,u)}\big\|_{q_1}\leq  \frac{(2-p)}{\mu^\frac{2-p}{2}}H(q_1)
 \big\|  J_\eta(\frac{\nabla\nabla u\,\chi^\rho}{a(\mu, u)})\big\|_{q_1}\leq c(\mu,\rho)H(q_1)(2-p) \big\| \nabla\nabla u\,\chi^\rho\big\|\,.$$
 Moreover $$\| \frac{f}{a(\mu,u)}\|_{q_1}\leq \frac{1}{\mu}{\atop^{\frac{2-p}{2}}} \|f\|_{q_1}\,.$$
 Let
 us multiply equation \eqref{stokesS2} by $\Delta w  |\Delta w|^{q_1-2}$  and integrate over $ \O$.
  By H\o lder's inequality
 and the previous estimates on the $L^{q_1}$-norms we get
 $$\ba{rl}\dy\vs1  \ve \|\Delta w\|_{q_1}^{q_1}\!\leq\!\!\!&\dy   \!\!
(2-p) \big\|A^\frac{4-p}{2}(\mu, u)\cdot J_\eta(\frac{\nabla\nabla u\,\chi^\rho}{a(\mu, u)})\big\|_{q_1}
 \|\Delta w\|_{q_1}^{q_1-1}\\ \vs1 &\dy +\big\|\frac{\nabla \Pi(u, \chi^\rho)}{a(\mu,u)}\big\|_{q_1}  \|\Delta w\|_{q_1}^{q_1-1}\! +\dy
 \frac{1}{\mu}{\atop^{\frac{2-p}{2}}} \|f\|_{q_1} \|\Delta w\|_{q_1}^{q_1-1}\\
 \dy\! \leq\!\!\!&\dy\! \!\! \left(c(\mu,\rho)(2-p)(1+H(q_1))
 \big\| \nabla\nabla u\,\chi^\rho\big\|+ c(\mu)  \|f\|_{q_1}\right)\! \|\Delta w\|_{q_1}^{q_1-1}\,,
    \ea$$ from which, using \eqref{main1} and Lemma\,\ref{LL2} with $p=2$,  \eqref{passe2} follows.
\vskip0.2cm\noindent
{\it\underline{Step  II: convergence of the sequence $\{w^\delta\}$ to a solution $w$ of system:}}
 \be\label{stokesS3} \ve\Delta w+\frac{\Delta w}{a_{2-p}(\mu,u)}
=(2-p) A^\frac{4-p}{2}\cdot J_\eta(\frac{\nabla\nabla u\,\chi^\rho}{a(\mu, u)})+
\frac{\nabla \Pi(u, \chi^\rho)}{a(\mu,u)}+ \frac{f}{a(\mu,u)}, \mbox{ in } \O.
\ee\vskip0.1cm
We prove that, in the limit as $\delta\to 0$,
   there exists a limit function $w$ which satisfies the same bounds of the sequence $\{w^\delta\}$ and satisfies system \eqref{stokesS3} a.e. in $\O$.
   Thanks to  \eqref{passe1} and \eqref{passe2}, the sequence $\{ w^\delta\}$ is uniformly bounded in $W^{2,q}(\O)\cap W^{2,q_1}(\O).$ From the
Rellich compactness theorem, there exist a field $w$
and a subsequence such that
\be\label{straw}w^\delta\rightharpoonup w\textrm{  in }  W^{2,q}(\O), \ee
$$w^\delta\rightarrow w \textrm{ in } W^{1,r}(\O),$$
 for any
$r\in[1,\infty)$,
 and the limit $w$ satisfies estimates \eqref{passe1} and \eqref{passe2}.
 Therefore, from \eqref{straw}, for any  $\mu>0$ and for any $\vp\in L^{q'}(\O)$,
\be\label{weakt1}
\lim_{\delta \to 0}((\Delta w^\delta- \Delta w),\frac{\vp}{a_{2-p}(\mu, u)})=0\,.\ee
 Corresponding to this subsequence in $\delta$, let us consider the sequence $\{J_\delta(\nabla u)\}$. This sequence strongly converges to $\nabla u$ in $L^2(\O)$, hence there exists a subsequence such that
 \be\label{jd} J_\delta(\nabla u) \to \nabla u,\ \mbox{ almost everywhere in } \O\,.\ee
 From \eqref{weakt1}, we have
$$\lim_{\delta \to 0}(\frac{\Delta w^\delta}{a_{2-p}(\mu,u ,\delta)}- \frac{\Delta w}{a_{2-p}(\mu, u)}, \vp)=
\lim_{\delta \to 0}(\Delta w^\delta(\frac{1}{a_{2-p}(\mu, u,\delta)}- \frac{1}{a_{2-p}(\mu, u)}, \vp) \,.$$
Applying H\o lder's inequality, thanks to the uniform boundedness of
$\|\Delta w^\delta\|_q$, property \eqref{jd} and then the application of the Lebesgue dominated convergence theorem, we obtain that the limit is equal to zero.
Therefore
$$(\ve\Delta w+\frac{\Delta w}{a_{2-p}(\mu,u)}
+(p-2)A^{\frac{4-p}{2}}(\mu, u)\cdot J_\eta(\frac{\nabla\nabla u\,\chi^\rho}{a(\mu, u)})-
\frac{\nabla \Pi(u, \chi^\rho)}{a(\mu,u)}- \frac{f}{a(\mu,u)},\vp)=0\,,$$
for all $\vp\in C_0^\infty(\O)$, which ensures
that $w$ satisfies \eqref{stokesS3} a.e. in $\O$.
\vskip0.2cm\noindent
{\it\underline{Step  III: $w=u$.}} -
 In Step II we have constructed a solution $w$ of problem
\eqref{stokesS3} belonging to
$W^{2, q}(\O)$.  By
taking the difference of
\eqref{stokesSo} and \eqref{stokesS3},
side by side, and by setting
$v=u-w$, with $u$ constructed in Proposition \ref{Vteo2}, we get that $v\in W_0^{1,2}(\O)\cap W^{2,2}(\O)$ satisfies
$$\ve \Delta v+\frac{\Delta v}{a_{2-p}(\mu,u)}=0\,, \mbox{ in } \O\,.$$
Multiplying the above equation by
$\Delta v$ and integrating in
$\O$, it follows that
$$\ve \|\Delta v\|+\| \frac{\Delta v}{a(\mu,u)}\| =0,$$hence $\|\,\Delta v\,\|=0$, and, since $v=0$ on $\pa\O$,
 $v\equiv 0$ holds. This, in particular, implies that $u\in W^{2,q}(\O)$,
$D^2u\in L^{q_1}(\O)$,
and $u$ satisfies estimates \eqref{passeu1} and \eqref{passeu2} as required. \chiu  \vskip0.1cm

 \section{ \large The solution of problem \eqref{stokesS1}}\label{secnoep}
In the previous section we have obtained an approximating sequence whose existence and regularity are independent of the particular bounded domain $\O$. Here our aim is the construction of a sufficiently regular approximating sequence, defined on the whole $\R^n$. \par In the proof of Proposition\,\ref{Vteo1} below we follow the arguments used in paper \cite{hey} (sec. 2)
for the construction of the solution of the Navier-Stokes equations in unbounded
domains. The difference is in the fact that in system \eqref{sm} below the convergence is ensured by lemmas \ref{converg1}--\ref{converg3} as in the proof of Proposition\,\ref{Vteo2}.\begin{prop}\label{Vteo1}{\sl
Let $\ov M(2)>0$. Let $q\in (n,+\infty)$, $q_1=\frac{np}{n+p}$, and assume that $f\in L^{q_1}(\R^n)\cap L^q(\R^n)$.
 Then there exists a solution $u\in \widehat W^{1,p}(\R^n)$ of problem \eqref{stokesS1}, such that $u\in L^{\frac{np}{n-p}}(\R^n)$, and
 \be\label{passeur1}\|D^2 u\|_{q} \leq C(\mu,\eta,\rho, \ve)\, (\|f\|_{q_1}+ \|f\|_{q}) \,,\ee
\be\label{passeur2}\|D^2 u\|_{q_1} \leq C(\mu,\eta,\rho,\ve)\, (\|f\|_{q_1}+ \|f\|_{q}) \,.\ee
}
\end{prop}
\Pr
Let $\{B_\sigma\}$ be a sequence of open balls such that
\be\label{ball}
\overline B_\sigma\subset  B_{\sigma+1}\ \mbox{  and }\ \R^n =\bigcup_{\sigma=1}^\infty B_\sigma\,.\ee
For any fixed $\rho>0$, there exists a positive integer $\sigma_\circ$ such that $\sigma_\circ>2\rho$. Therefore, for any $\sigma\geq \sigma_\circ$, $supp\,\chi^\rho\subset B_\sigma$.
  Due to Proposition\,\ref{Vteo2} and Proposition\,\ref{proplq}, in
each $B_\sigma$ we find a solution $u^\sigma\in W_0^{1,2}(B_\sigma)\cap W^{2,q}(B_\sigma)$ of
system \eqref{stokesSo}, which satisfies the estimates \eqref{main1}, \eqref{passeu1} and \eqref{passeu2},
uniformly in $\sigma\in \N$. Let $l\leq \sigma$. Since $D^2 u^\sigma\in L^{q_1}(B_\sigma)\cap L^q(B_\sigma)$, uniformly in $\sigma$, and since $p<2$, by interpolation
$D^2 u^\sigma$ are bounded in $L^{\frac{2n}{n+2}}(B_\sigma)$, and there exists a constant $C$, independent of $\sigma$, such that
 $$\|D^2 u^\sigma\|_{L^\frac{2n}{n+2}(B_l)} \leq C(\mu,\eta,\rho, \ve)\, (\|f\|_{q_1}+ \|f\|_{q}) \,.$$
Then, thanks to Lemma\,\ref{GNT}, we also have
 $$\|\nabla u^\sigma\|_{L^{2}(B_l)} \leq C(\mu,\eta,\rho,\ve)\, (\|f\|_{q_1}+ \|f\|_q) \,,$$
 and, by Sobolev's imbedding,
 $$\|u^\sigma\|_{L^{\frac{2n}{n-2}(B_l)}} \leq C(\mu,\eta,\rho,\ve)\, (\|f\|_{q_1}+ \|f\|_q) ,$$
with constants independent of $\sigma$.
 Collecting all the $\sigma$-uniform estimates thus obtained, we find that, for any fixed $l\leq \sigma$, $u^\sigma$  satisfies
 \be\label{aa}\|D^2 u^\sigma\|_{L^q(B_l)} \leq C(\mu,\eta,\rho, \ve)\, (\|f\|_{q_1}+ \|f\|_q) \,,\ee
\be\label{ab}\|D^2 u^\sigma\|_{L^{q_1}(B_l)} \leq C(\mu,\eta,\rho,\ve)\, (\|f\|_{q_1}+ \|f\|_q) \,,\ee
 \be\label{ac}\|\nabla u^\sigma\|_{L^{2}(B_l)} \leq C(\mu,\eta,\rho,\ve)\, (\|f\|_{q_1}+ \|f\|_q)\,, \ee
  \be\label{ad}\|u^\sigma\|_{L^{\frac{2n}{n-2}}(B_l)} \leq C(\mu,\eta,\rho,\ve)\, (\|f\|_{q_1}+ \|f\|_q)\,. \ee
Thus, for any fixed $l$, there exists a subsequence $\{u^\sigma\}$, again labeled in $\sigma$,  which
weakly converges in the norms
 listed in \eqref{aa}--\eqref{ad}. Passing to subsequences, the convergences are obtained for
 any $l$. Denoting by $\{u^\sigma\}$ the final subsequence and by $u$ its limit,
 from the lower semicontinuity the limit $u$ verifies the same estimates
 \eqref{aa}-- \eqref{ad}. Further, since $\nabla u\in \widehat W^{1,q_1}(\R^n)\cap L^2(\R^n)$,
 by using Sobolev's inequality, we also have that $\nabla u\in L^p(\R^n)$. Similarly, $u\in L^{\frac{np}{n-p}}(\R^n)$.
 \par Let us prove that
  $u$ satisfies system \eqref{stokesS1} in $\R^n$. Let us consider a generic but fixed  function
  $\vp\in C_0^\infty(\R^n)$ and let $K$ be its compact support. Then, there exists a $\ov \sigma\in \N$, such that
  $K\subset B_{\ov \sigma}$. Hence,
   for any $\sigma\geq \ov \sigma$, $\vp\in C_0^\infty(B_\sigma)$, and, since $u^\sigma$ is solution of system \eqref{stokesSo} in $B_\sigma$, we have
 \be\label{sm}\ba{ll}\dy\vs1  (\ve\Delta u^\sigma+\frac{\Delta u^\sigma}{a_{2-p}(\mu,u^\sigma)}
+(p-2)A^\frac{4-p}{2}(\mu, u^\sigma)
\cdot J_\eta(\frac{\nabla\nabla u^\sigma\,\chi^\rho}{a(\mu, u^\sigma)})\\ \dy \hskip6cm-
\frac{\nabla \Pi(u^\sigma, \chi^\rho)}{a(\mu,u^\sigma)}- \frac{f}{a(\mu,u^\sigma)},\vp)=0\,.\ea \ee
 We can  pass to the limit in the above integral identity, as $\sigma$ tends to $\infty$,
by using the same arguments employed in Proposition\,\ref{Vteo2} for the
 convergence of the Faedo-Galerkin approximation.
 Finally we find that
 \be\label{sm1}\ba{ll}\vs1\dy (\ve\Delta u+ \frac{\Delta u}{a_{2-p}(\mu,u)}
+(p-2)A^\frac{4-p}{2}(\mu, u)\cdot J_\eta(\frac{\nabla\nabla u\,\chi^\rho}{a(\mu, u)})-
\frac{\nabla \Pi(u, \chi^\rho)}{a(\mu,u)}- \frac{f}{a(\mu,u)},\vp)=0\,,\ea\ee
for any $\vp\in C_0^\infty(\R^n)$, which completes the proof.
\chiu
   \section{\large The solution of problem \eqref{stokesS4} ($\ve\to 0$)}\label{secep}
Starting from this section, our aim
is to let all the parameters, $\ve,
\eta, \rho, \mu$, tend to zero in
\eqref{stokesS1}, preserving the
regularity properties of the
approximating sequence. We begin
with the limit as $\ve$ tend to
zero, and show that the sequence of
solutions $\{u^{\mu, \rho, \eta,
\ve}\}$ of \eqref{stokesS1},
obtained in Proposition
\ref{Vteo1}, converges to a limit
function $u^{\mu, \rho, \eta}$
which satisfies system
\eqref{stokesS4}, and estimates
\eqref{passeue2}--\eqref{passeqq1}
below. The first step consists in
proving $\ve$-uniform estimates on
$\{u^{\mu, \rho, \eta, \ve}\}$.
This is the aim of Proposition
\ref{Vteo4}. As a matter of course,
these are also $\eta$-uniform
estimates. For the reader's
convenience we denote the family of
solutions $u^{\mu, \rho, \eta,
\ve}$ simply by $u$. From now on we
will work in the whole space
$\R^n$. Hence the symbol
$\|\cdot\|_r$ will always denote
the $L^r(\R^n)$-norm.
 Recall that  \be\label{Cqr}
M(r)=1-(2-p)H(r')(5+H(r)).\ee
Note that the assumption $M(r)>0$ is equivalent to the pair of assumptions
 \be\label{conL3.3}
 1-4H(r')(2-p)>0\ \mbox{ and } \ 1-\frac{(2-p)H(r')}{1-4H(r')(2-p)}(1+H(r))>0.\ee
\begin{prop}\label{Vteo4}{\sl
Let the assumption of Proposition\,\ref{Vteo1} be satisfied.  Further, assume that $M(q_1)$ and
$M(q)$ are positive. Then there exist some constants $C$,
independent of $\eta$, $\rho$ and $\ve$, such that, for $r=q_1$ and $r=q$,
     \be\label{passeue2}
     \big\|\frac{\nabla\nabla u}{a(\mu, u)}\big\|_{r}
     \leq C(r)\, \|f\|_r\,,  \ee
  \be\label{passepe2}
     \big\|\nabla\Pi(u,\chi^\rho)\big\|_r \leq C(r)\, \|f\|_r\,,   \ee
    \be\label{passeqq1} \big\|{D^2u}\big\|_{q_1}+\big\|{D^2u}\big\|_{q}\leq\dy
C(q_1,q)\, \mu^{\frac{ 2-p}{2} } (\|f\|_{q_1}+\|f\|_q) +C(q_1,q)\, (\|f\|_{q_1}+\|f\|_q) ^\frac{1}{p-1}\,.\ee
Finally,
\be\label{gp1}\|\nabla u\|_p+\|\Pi(u,\chi^\rho)\|_p\leq \dy
C(q_1,q)\, \mu^{\frac{ 2-p}{2} } (\|f\|_{q_1}+\|f\|_q) +C(q_1,q)\, (\|f\|_{q_1}+\|f\|_q) ^\frac{1}{p-1}\,.\ee
}
\end{prop}
 \Pr
 Let be either $r=q_1$ or $r=q$. Let us multiply equation \eqref{stokesS1} by $\frac{\Delta u |\Delta u|^{r-2}}{a(\mu, u)^{r-2}}$. By H\o lder' s inequality
 we get
\be\label{hold}\ba{ll}\dy\vs1  \ve\!\int\!\!\frac{|\Delta u|^r}{a(\mu,\!u)}{\atop^{r-2}}\,dx+ \big\|\frac{\Delta u}{a(\mu,\!u)}\big\|_r ^{r}\!\leq\!
(2-p) \big\|A^1(\mu,\!u)\cdot J_\eta(\frac{\nabla\nabla u\chi^\rho}{a(\mu,\!u)})\big\|_{r}
 \big\|\frac{\Delta u}{a(\mu, u)}\big\|_{r}^{r-1}\\ \hfill \dy +\|\nabla \Pi(u,\chi^\rho)\|_{r}  \big\|\frac{\Delta u}{a(\mu, u)}\big\|_{r}^{r-1}\! +
 \|f\|_{r}  \big\|\frac{\Delta u}{a(\mu, u)}\big\|_{r}^{r-1}.
    \ea\ee
Inequality \eqref{hold} is well-posed, via inequalities  \eqref{passeur1} and \eqref{passeur2}. Hence
for any $\ve>0$, $\eta>0$ and $\mu>0$,  $ \frac{\nabla\nabla u }{a(\mu, u)}\in L^r(\R^n)$. Then, since $ 1-4H(r')(2-p)>0$,
   applying Lemma\,\ref{LL2} with $p<2$, we get
\be\label{nn1} \ba{ll}\dy\vs1\big\|A^1(\mu, u)\cdot J_\eta(\frac{\nabla\nabla u\,\chi^\rho}{a(\mu, u)})\big\|_r\\
\dy \leq
 \big\|  J_\eta(\frac{\nabla\nabla u\,\chi^\rho}{a(\mu, u)})\big\|_r\leq   \big\|
 \frac{\nabla\nabla u }{a(\mu, u)}\big\|_r\leq \frac{H(r')}{1-4H(r')(2-p)}\big\|
 \frac{\Delta u }{a(\mu, u)}\big\|_r\,.\ea\ee
Since from \eqref{pi} we can apply to $\nabla\Pi(u, \chi^\rho)$ the Calder\'on-Zygmund theorem, using the above estimate we find
\be\label{nn2}\big\|\nabla \Pi(u, \chi^\rho)\big\|_r\leq (2-p) H(r)
 \big\|  J_\eta\big(\frac{\nabla\nabla u\,\chi^\rho}{a(\mu, u)}\big)\big\|_r\leq \frac{(2-p)H(r')H(r)}{1-4H(r')(2-p)}\big\|
 \frac{\Delta u }{a(\mu, u)}\big\|_r
 \,.\ee
 Therefore estimate \eqref{hold} gives
 $$
\big(1-\frac{(2-p)H(r')}{1-4H(r')(2-p)}(1+H(r))\big)\big\|\frac{\Delta u}{a(\mu, u)}\big\|_r\leq \|f\|_{r}\,,$$
which, provided that $M(r)>0$,
with a further application of Lemma\,\ref{LL2}, implies  \eqref{passeue2}.
 From \eqref{nn2} and \eqref{passeue2}, we also obtain \eqref{passepe2}.
\par Estimate \eqref{passeue2} is the main tool in deducing an estimate on $D^2u $ in the $L^q$-norm and $L^{q_1}$-norm,
uniformly in $\ve>0$ and $\eta>0$, hence proving \eqref{passeqq1}.
   Recalling that, from Proposition\,\ref{Vteo1},  $u\in \widehat W^{1,p}(\R^n)$,  applying
    Gagliardo-Nirenberg's inequality (see \cite{gagliardo} or \cite{nirenberg} or also \cite{CM1}),  we  have
      \be\label{GNinfty}
      \, \|\nabla u\|_\infty\leq  c\, \|D^2 u\|_q^a\|\nabla u\|_{p}^{1-a}\leq  c\, \|D^2 u\|_q^a\|D^2 u\|_{q_1}^{1-a}\,
      ,\ee with $a=\frac{nq}{n(q-p)+pq}$.
Let us multiply
system \eqref{stokesS1} by $a_{2-p}(\mu, u) \Delta u |\Delta u|^{r-2}$, $r=q_1$ or $r=q_1$, and integrate over $\R^n$.
By H\o lder' s inequality
 we get
$$\ba{ll}\dy\vs1  \ve \int |\Delta u|^r a_{2-p}(\mu,u)\,dx+ \big\|{\Delta u}\big\|_r ^{r}\\ \vs1\hfill \leq\dy \!
(2-p) \big\|A^{1}(\mu, u)\cdot J_\eta(\frac{\nabla\nabla u\,\chi^\rho}{a(\mu, u)})\big\|_{r}
 \big\|\Delta u\big\|_r^{r-1} \big\|a(\mu, u)\big\|_{\infty}\\\hfill+\|\nabla \Pi(u,\chi^\rho)\|_{r}
  \big\|\Delta u\big\|_r^{r-1}\|a(\mu, u)\big\|_{\infty} +
 \|f\|_{r}  \big\|\Delta u\|_r^{r-1}\|a(\mu, u)\|_{\infty},
    \ea$$
    which, by using firstly estimates \eqref{nn1} and \eqref{nn2} and then estimate \eqref{passeue2},   gives
      \be\label{holdq}\ba{rl}\dy\vs1  \big\|{\Delta u}\big\|_r\!\!\!\!\! &\leq\dy  \!
\frac{(2-p)H(r')} {1-4H(r')(2\!-\!p)}(1+ H(r)) \big\|
 \frac{\Delta u }{a(\mu,\!u)}\big\|_r
 \big\|a(\mu,\!u)\big\|_{\infty}\!\!+ \|f\|_{r} \|a(\mu,\!u)\|_{\infty}\\
\hfill \!\! &\leq \dy C(r)\|f\|_{r} \|a(\mu, u)\|_{\infty}\,.
 \ea\ee
  Finally, by using the estimate \be\label{qd}
 \|D^2 u\|_{r} \leq H (r)\|\Delta u\|_{r}\,,\ee
 which follows from Lemma\,\ref{LL2} with $p=2$,
    and summing \eqref{holdq} written once with $r=q_1$ and then with $r=q$,
    we find
       \be\label{holdqsum}\ba{rl}\dy  \big\|{D^2u}\big\|_{q_1}+\big\|{D^2u}\big\|_{q}\!\!\! &\leq\dy
C(q_1,q) (\|f\|_{q_1}+\|f\|_q) \|a(\mu, u)\|_{\infty}\,.
 \ea\ee
   From \eqref{GNinfty}, the $L^{\infty}$-norm of $a(\mu,u)$ can be estimated as follows
   \be\label{amuuinfty}\|a(\mu, u)\|_{\infty}\leq  \mu^{\frac{ 2-p}{2} }+ \, \|\nabla u\|_\infty^{2-p}\leq
       \mu^{\frac{ 2-p}{2} }+ c\, (\|D^2 u\|_q+\|D^2u \|_{q_1})^{2-p}\,. \ee
Using estimate \eqref{amuuinfty} in \eqref{holdqsum} and then Young's inequality,
    we find \eqref{passeqq1}. Estimate \eqref{gp1} for $\nabla u$ follows by Sobolev's embedding and estimate \eqref{passeqq1}, while the estimate on $\Pi(u,\chi^\rho)$ follows from the Hardy-Littlewood-Sobolev theorem and the same estimate \eqref{passeqq1}.
  \chiu
\vskip0.1cm
We are now in position to perform the limit as $\ve\to 0$. \begin{prop}\label{Vteo3}{\sl
Let $q\in (n,+\infty)$, $q_1=\frac{np}{n+p}$, and assume that $\ov M(2)$, $M(q_1)$ and $M(q)$ are positive constants.
Let $f\in L^{q_1}(\R^n)\cap L^q(\R^n)$.
 Then there exist a solution $u$  of system \eqref{stokesS4}, and  some constants $C$,
independent of $\eta$, $\rho$, such that, for $r=q_1$ and $r=q$,
 \be\label{passeunoee}
     \big\|\frac{\nabla\nabla u}{a(\mu, u)}\big\|_{r}\leq C(r)\, \|f\|_{r}\,,  \ee
     \be\label{spres1} \|\nabla \Pi(u,\chi^\rho )\|_r\leq C(r)\, \|f\|_{r} \,,\ee
        \be\label{passeqq1noee} \big\|{D^2u}\big\|_{q_1}+\big\|{D^2u}\big\|_{q}\leq\dy
C(q_1,q)\, \mu^{\frac{ 2-p}{2} } (\|f\|_{q_1}+\|f\|_q) +C(q_1,q)\, (\|f\|_{q_1}+\|f\|_q) ^\frac{1}{p-1}\,,\ee
where $\Pi(u, \chi^\rho)$ is given by \eqref{pi}.
Finally,
\be\label{gp2}\|\nabla u\|_p+\|\Pi(u,\chi^\rho)\|_p\leq \dy
C(q_1,q)\, \mu^{\frac{ 2-p}{2} } (\|f\|_{q_1}+\|f\|_q) +C(q_1,q)\, (\|f\|_{q_1}+\|f\|_q) ^\frac{1}{p-1}\,.\ee
 }
\end{prop}
\Pr By virtue of Proposition \ref{Vteo4}, for any $\ve>0$
there exists a solution $u^\ve\in \widehat W^{1,p}(\R^n)$ of system \eqref{stokesS1}, satisfying properties \eqref{passeue2}--\eqref{passeqq1}.  These bound imply the existence of a field $u$ and a subsequence such that
$$\nabla u^\ve\rightharpoonup \nabla u \textrm{ in } L^{p}(\R^n),$$
$$D^2 u^\ve\rightharpoonup D^2 u \textrm{ in } L^{q_1}(\R^n)\cap L^q(\R^n),$$
and, from the lower
semi-continuity, estimate \eqref{passeqq1noee} and the $L^p$ bound for $\nabla u$ given in \eqref{gp2} hold. Further,
we can apply Lemma \ref{converg1}, either with $r=q$, $\ov r=q_1$, $E=\R^n$, and with $r=\ov r=q_1$, and Lemma \ref{converg2} with $h=\chi^\rho$. Therefore, along a subsequence we find that,
in the limit as $\ve\to 0$,
\be\label{weakte}
\frac{\Delta u^\ve}{a_{2-p}(\mu, u^\ve)}\rightharpoonup\frac{\Delta u}{a_{2-p}(\mu, u)}\,,\textrm{  in }  L^{q_1}(\R^n)\cap L^{q}(\R^n),\ee
 \be\label{forpe}
A^\frac{4-p}{2}(\mu, u^\ve) J_\eta(\frac{\nabla\nabla u^\ve\,\chi^\rho}{a(\mu, u^\ve)})\rightharpoonup
A^\frac{4-p}{2}(\mu, u)\,)\cdot J_\eta(\frac{\nabla\nabla u\,\chi^\rho}{a(\mu, u)})\textrm{  in }  L^{q_1}(\R^n)\cap L^{q}(\R^n),\ee
and estimates \eqref{passeunoee} holds.
Finally, let us consider the corresponding subsequence of $\{\nabla \Pi^\ve\}$, where $\Pi^\ve:=\Pi(u^\ve, \chi^\rho)$ is defined in \eqref{pi}.
 From estimate \eqref{gp1}, the sequence
 $\{\Pi^\ve\}$ is bounded in $L^p(\R^n)$.
 Moreover, from estimate \eqref{passepe2}, it follows that the sequence $\{\nabla\Pi^\ve\}$  is uniformly bounded in $L^{q_1}(\R^n)\cap L^q(\R^n)$.
Therefore there exist a limit $\Pi$ and a subsequence such that, in the limit as $\ve$ tends to zero,
$$\Pi^\ve\rightharpoonup\Pi \textrm{  in } L^{p}(\R^n)\,,$$
$$\nabla\Pi^\ve\rightharpoonup \nabla\Pi \textrm{  in } L^{q_1}(\R^n)\cap L^q(\R^n), $$
and, thanks to Lemma\,\ref{converg1},
$$\frac{\nabla\Pi^\ve}{a(\mu, u^\ve)}\rightharpoonup \frac{\nabla\Pi}{a(\mu, u)}\textrm{  in } L^{q_1}(\R^n)\cap L^q(\R^n)\,.$$
The limit clearly satisfies estimates \eqref{spres1} and \eqref{gp2}. It  remains to show that
$\Pi=\Pi(u, \chi^\rho)\, $. It is enough to show that there exists a subsequence such that, for any $\vp\in C_0^\infty(\R^n)$,
$$\lim_{\ve\to 0}\, (\Pi^\ve,\vp)=(\Pi(u, \chi^\rho),\vp)\,. $$
We set
$$g(y)= \int_{\R^n}\nabla_x {\mathcal E} (x-y)\, \vp(x)\,dx\,.$$Then
$$\ba{rl}\vs1\dy
(\Pi^\ve,\vp)\dy=&\dy (p-2) \int_{\R^n}\frac{\nabla u_i^\ve   (y)D_{y_j}u_h^\ve (y)}{a_1(\mu, u^\ve(y))} J_\eta(\frac{D^2_{y_iy_j} u_h^\ve(y)\,\chi^\rho(y)}{a(\mu, u^\ve(y))}\big)\cdot g(y)\,dy\,,
\ea$$
and the convergence to $(\Pi(u, \chi^\rho),\vp)$ can be obtained from Lemma\,\ref{converg2}, observing that,
from the Hardy-Littlewood-Sobolev theorem,
the function $g$ belongs to $L^{q_1'}(\R^n)$ . 
\par Since clearly one has
 $$ \lim_{\ve\to 0}\,\ve(\Delta u^\ve, \vp)\,=0,$$
 we have obtained
$$(\frac{\Delta u}{a_{2-p}(\mu,u)}
+(p-2) A^{\frac{4-p}{2}}(\mu, u)\cdot J_\eta(\frac{\nabla\nabla u\,\chi^\rho}{a(\mu, u)}-
\frac{\nabla \Pi(u, \chi^\rho)}{a(\mu,u)}- \frac{f}{a(\mu,u)},\vp)=0\,,$$
for all $\vp\in C_0^\infty(\R^n)$, which completes the proof. 
\chiu

\section{\large  The solution of problem \eqref{stokesS5} ($\eta\to 0$)}\label{secet}
Here we investigate on further regularities of $u$, with bounds uniform in $\eta>0$ and $\rho>0$.
 Subsequently we pass to the limit as $\eta$ goes to zero and prove that the limit function solves  system
 \eqref{stokesS5} and has the properties stated in Proposition\,\ref{Vteo7}.
\begin{prop}\label{Vteo6}{\sl Let $q\in (n,+\infty)$, $q_1=\frac{np}{n+p}$, and assume that $\ov M(2)$, $M(q_1)$ and $M(q)$ are positive constants.
Let $f\in L^{q_1}(\R^n)\cap L^{q}(\R^n)$,
with $\nabla f\in L^q(\R^n)$, and let $u$ be the solution of system \eqref{stokesS4}
obtained in Proposition\,\ref{Vteo3}. Then,
there exists a constant $C=C(\mu, \|f\|_{q_1}, \|f\|_{q},\|\nabla f\|_{q})$, independent of $\eta$ and $\rho$, such that
     \be\label{trev}
     \big\|D^3 u\big\|_{q}\leq C\,.
      \ee
}
\end{prop}
 \Pr 
 \noindent {\it Step I: \underline{$\forall \eta>0, \ D^3 u\in L^q(\R^n)$.}}
 \vskip0.1cm Let us multiply both sides of system \eqref{stokesS4} by $a_{2-p}(\mu, u)$:
\be\label{pois}
\Delta u=(2-p) \frac{\nabla u\otimes \nabla u}{a_{\frac{p}{2}}(\mu, u)}\cdot J_\eta(\frac{\nabla\nabla u\,\chi^\rho}{a(\mu, u)})+
\nabla \Pi(u,\chi^\rho)\,{a(\mu,u)}+f\,{a(\mu,u)}\,,\ \mbox{ in } \R^n\,.\ee
If  the right-hand side of \eqref{pois} belongs to $W^{1,q}(\R^n)$, then, by classical results on the Poisson's equation, $D^2u\in L^q(\R^n)$
and $D^3u\in L^q(\R^n)$. From Proposition\,\ref{Vteo3}, the right-hand side belongs to $L^q(\R^n)$, and the estimates are uniform in $\eta$ and $\rho$.
 Let us verify that the first order derivatives of the right-hand side also belong to $L^q(\R^n)$.
 We have
$$\ba{ll}\vs1\dy D_{x_i}\big(A^{\frac{p}{2}}(\mu, u)\cdot J_\eta(\frac{\nabla\nabla u\,\chi^\rho}{a(\mu, u)})\big)=
D_{x_i}A^{\frac p2}(\mu, u)\cdot J_\eta(\frac{\nabla\nabla u\,\chi^\rho}{a(\mu, u)})\\
\vs1\hfill\dy +A^\frac p2 (\mu, u) \cdot D_{x_i} J_\eta(\frac{\nabla\nabla u\,\chi^\rho}{a(\mu, u)})\,.\ea$$
By interpolation, all the tems are bounded in $L^r(\R^n)$, for any $r\in [q_1, q]$.
Indeed, from Proposition \ref{Vteo3},  $D^2 u\in L^q(\R^n)\cap L^{q_1}(\R^n)$,  and  $\nabla u\in L^\infty(\R^n)$,
thanks to estimate
\eqref{GNinfty}. Moreover
$J_\eta(\frac{\nabla\nabla u\,\chi^\rho}{a(\mu, u)})\in L^\infty(\R^n) $ and
$D_{x_i} J_\eta(\frac{\nabla\nabla u\,\chi^\rho}{a(\mu, u)})\in L^r(\R^n)$
thanks to the properties of the mollifier. Further, recalling the definition of $\Pi(u,\chi^\rho)$ given in \eqref{pi} and taking its gradient, we can write
\be\label{rapf}
\nabla\Pi(u,\chi^\rho)=(p-2)\int_{\R^n}\nabla_x {\mathcal E} (x-y)D_{y_i}\, \big(
\frac{D_{y_i}u_j \nabla u}{a_1(\mu, u)}\cdot J_\eta(\frac{D_{y_j}\nabla u\,\chi^\rho}{a(\mu, u)})\big)\,dy\,.\ee
The integration by parts is possible since we have proved above that the density in \eqref{rapf} belongs to $L^r(\R^n)$,
for any $r\in [q_1,q]$. Note that, from the Hardy-Littlewood-Sobolev theorem, $\nabla\Pi(u,\chi^\rho)\in L^s(\R^n)$, with  $s=\frac{nr}{n-r}$ for any
$r\in [q_1,n)$.  Therefore  $\nabla\Pi(u,\chi^\rho)\in L^s(\R^n)$, for any $s\in [\frac{nq_1}{n-q_1}, +\infty)$. Finally,
the Calder\'on-Zygmund theorem ensures in particular that $\nabla\nabla \Pi(u,\chi^\rho)\in L^q(\R^n)$.
\vskip0.2cm\noindent {\it Step II: \underline{There exist a constant $C(\mu)$, independent of $\eta$ and $\rho$, such that}}
     \be\label{tre1}
     \big\|\frac{D^3 u}{a(\mu, u)}\big\|_{q}\leq C(\mu)\|D^3 u\|_q^\frac{n}{q}\|D^2 u\|_q^{2-\frac{n}{q}}+ C(\mu)\|D^2 u\|_q+\|\nabla f\|_q\,.
      \ee
    We remark that, by interpolation,
 \be\label{II}
 \|D^2 u\|_{2q}\leq c\|D^3 u\|_q^\frac{n}{2q}\|D^2 u\|_q^{1-\frac{n}{2q}}\,.\ee
 By multiplying system \eqref{stokesS4} by $a(\mu, u)$, we obtain
\be\label{stokesS6}\frac{\Delta u}{a(\mu,u)}
+ (2-p)
 \frac{\nabla u\otimes \nabla u}{a_1(\mu, u)}\cdot J_\eta(\frac{\nabla\nabla u\,\chi^\rho}{a(\mu, u)})
 =
\nabla \Pi(u,\chi^\rho)+ f\,,\ \mbox{ in } \R^n\,.
\ee
 Thanks to Step I, we can derive equation \eqref{stokesS6}, and, a.e. in $\R^n$, we have
\be\label{stokesS7} \frac{D\Delta u}{a(\mu,u)}
=\sum _{i=1}^{6}\, I_i\,,\ee
with $$I_1=\frac{(2-p)}{2}\frac{\Delta u\, D|\nabla u|^2}{a_{\frac{4-p}{2}}(\mu,u)}\,,$$
$$\ba{ll} \vs1\dy I_2= (p-2)
 \frac{D\nabla u\otimes \nabla u}{a_{1}(\mu, u)}\cdot J_\eta(\frac{\nabla\nabla u\,\chi^\rho}{a(\mu, u)})
\dy +(p-2)\frac{\nabla u\otimes D\nabla u}{a_1(\mu, u)}\cdot J_\eta(\frac{\nabla\nabla u\,\chi^\rho}{a(\mu, u)})
\\\dy\hfill +(2-p)\frac{\nabla u\otimes \nabla u \cdot D|\nabla u|^2}{a_2(\mu, u)}\cdot J_\eta(\frac{\nabla\nabla u\,\chi^\rho}{a(\mu, u)})\,,\ea$$
 $$ I_3=
(p-2) \frac{\nabla u\otimes \nabla u}{a_1(\mu, u)}\cdot J_\eta(\frac{D\nabla\nabla u\,\chi^\rho}{a(\mu, u)})\,,
$$
$$I_4=(p-2) \frac{\nabla u\otimes \nabla u}{a_1(\mu, u)}\cdot J_\eta(\frac{\nabla\nabla u\,D\chi^\rho}{a(\mu, u)})+
(p-2) \frac{\nabla u\otimes \nabla u}{a_1(\mu, u)}\cdot J_\eta(\frac{\nabla\nabla u\,\chi^\rho\, D|\nabla u|^2}{a_{\frac{4-p}{2}}(\mu, u)})\,,$$
$$I_5=
D\nabla \Pi(u,\chi^\rho)\,,$$
$$I_6=
D\,f\,.$$
Let us estimate the $L^q$-norm of these terms. We have
\be\label{I1}\|I_1\|_q\leq (2-p)C(\mu) \|D^2 u\|_{2q}^2\,,\ee
\be\label{I234}\ba{rl}\dy\vs1
 \|I_2\|_q+\|I_4\|_q
\leq&\dy \!\!(2-p)C(\mu) (\|D^2 u\|_{2q}\|D^2 u\,\chi^\rho\|_{2q}+ \|D^2 u\|_{q})\\
\leq &\dy \!\!(2-p)C(\mu) (\|D^2 u\|_{2q}^2+\|D^2u\|_q)\,.\ea\ee
Moreover, using the properties of the mollifier and then applying Lemma\,\ref{LL3} with $p<2$, we get
\be\label{I98}\ba{rl}\dy \vs1\|I_3\|_q &\dy \leq (2-p) \|J_\eta(\frac{D\nabla\nabla u\,\chi^\rho}{a(\mu, u)})\|_{q}
\leq (2-p) \|\frac{D\nabla\nabla u}{a(\mu, u)}\|_{q}
\\ \dy &\dy \leq (2-p)H(q') \|\frac{D\Delta u}{a(\mu, u)}\|_{q}+ 4(2-p)^2H(q')C(\mu)\|D^2u\|_{2q}^2\,.\ea\ee
These estimates, together with the representation formula \eqref{rapf} and the
Calder\'on Zygmund theorem, ensure that
\be\label{I89}\|I_5\|_q  \leq (2-p)H(q')H(q)\big( C(\mu)\|D^2 u\|_{2q}^2+
\|\frac{D\Delta u}{a(\mu, u)}\|_{q}\big)\,.\ee
Finally, due to the assumptions on $f$,
\be\label{I1011}\|I_{6}\|_q
 \leq \, \|\nabla \, f\|_{q}\,.\ee
Let us multiply equation \eqref{stokesS7} by $\frac{D\Delta u |D\Delta u|^{q-2}}{a{(\mu, u)^{q-1}}}$ and integrate on $\R^n$. By using
 H\o lder' s inequality, then the above estimates \eqref{I1}--\eqref{I1011},  and finally the interpolation inequality \eqref{II},  we get
\be\label{holdS}\ba{ll}\vs1\dy \big\|\frac{D\Delta u}{a(\mu,u)}\big\|_q^q\leq &\!\!\!\!\dy (C(\mu)\|D^3 u\|_q^\frac{n}{q}\|D^2 u\|_q^{2-\frac{n}{q}}+ C(\mu)\|D^2 u\|_q+\|\nabla f\|_q)
 \big\|\frac{D\Delta u}{a(\mu,u)}\big\|_q^{q-1}\\& \dy +(2-p)H(q')(1+H(q))
  \big\|\frac{D\Delta u}{a(\mu,u)}\big\|_q ^q\,.
\ea\ee
Since $M(q)>0$, one has
 $1-(2-p)H(q')(1+H(q))>0$. Therefore, we have found that
$$\big\|\frac{D\Delta u}{a(\mu,u)}\big\|_q\leq  C(\mu)\|D^3 u\|_q^\frac{n}{q}\|D^2 u\|_q^{2-\frac{n}{q}}+ C(\mu)\|D^2 u\|_q+\|\nabla f\|_q\,, $$
 that, with a further application of Lemma\,\ref{LL3}, gives estimate \eqref{tre1}, uniformly in $\eta>0$ and $\rho>0$.
  \vskip0.2cm\noindent {\it Step III: \underline{$D^3 u\in L^q(\R^n)$, uniformly in $\eta>0$ and $\rho>0$.}}
 \vskip0.1cm
 Let us multiply
system \eqref{stokesS7} by $a(\mu, u)D \Delta u |D\Delta u|^{q-2}$,  and integrate over $\R^n$.
By H\o lder' s inequality
we get
$$\ba{ll}\dy\vs1  \big\|{D\Delta u}\big\|_q ^{q}\leq \sum_{i=1}^6 \|I_i\|_q \|a(\mu, u)\|_\infty \big\|{D\Delta u}\big\|_q ^{q-1}\ea$$
Hence, dividing by  $\big\|D\Delta
u\big\|_q^{q-1}$ and applying Lemma\,\ref{LL3}
with $p=2$, we find \be\label{M2q}\ba{ll}\dy\vs1
\big\|{D^3 u}\big\|_q\leq C\sum_{i=1}^6
\|I_i\|_q \|a(\mu, u)\|_\infty\,. \ea\ee Each
$I_i$ has been estimated in Step II. Further we
can increase the $L^{2q}$-norm of $D^2u$ in
\eqref{I1}--\eqref{I89} via \eqref{II}, and the
$L^q$-norm of $ \frac{D\Delta u}{a(\mu,u)}$ in
\eqref{I98} and \eqref{I89} via estimate
\eqref{tre1}. Hence, from \eqref{M2q}, a direct
calculation gives
$$\ba{rl}\dy\vs1
 \big\|{D^3 u}\big\|_q\leq\!\!\! &\dy C(\mu)\|D^3 u\|_q^\frac{n}{q}\|D^2 u\|_q^{2-\frac{n}{q}}\|a(\mu, u)\|_\infty\\ &\dy+C(\mu)\|D^2 u\|_q\|a(\mu, u)\|_\infty+\|\nabla f\|_q\ \|a(\mu, u)\|_\infty\,. \ea$$
We apply to the first term on the right-hand side Young's inequality, with exponents $\frac qn>1$ and $\frac{q}{q-n}>1$ . Hence,
 $$\ba{rl}\dy\vs1
 \big\|{D^3 u}\big\|_q\leq\!\!\! &\dy \frac 12\|D^3 u\|_q+ C(\mu)\|D^2 u\|_q^{\frac{2q-n}{q-n}}\|a(\mu, u)\|_\infty^\frac{q}{q-n}\\ &\dy+C(\mu)\|D^2 u\|_q\|a(\mu, u)\|_\infty+\|\nabla f\|_q\ \|a(\mu, u)\|_\infty\,.\ea$$
 Taking into account estimate \eqref{passeqq1noee} on $\|D^2u\|_q$, estimate \eqref{gp2} on $\|\nabla u\|_p$, and the consequent estimate of $\|\nabla u\|_\infty$,  one ends up with estimate \eqref{trev} with a constant $C$ independent of $\eta$ and $\rho$.   The proof is achieved.
\chiu

\begin{prop}\label{Vteo7}{\sl
Let the assumption of Proposition\,\ref{Vteo6} be satisfied. Then, there exists a solution $u$ of system \eqref{stokesS5}, and there exist some constants C, independent of $\mu$, such that
for $r=q_1$ and $r=q$,
 \be\label{tpasse}
     \big\|\frac{\nabla\nabla u}{a(\mu, u)}\big\|_{r}\leq C\, \|f\|_{r}\,,   \ee
     \be\label{spres2}
    \|\nabla \Pi(u,\chi^\rho )\|_r\leq C\, \|f\|_r\,,\ee
  \be\label{tpasseqq1} \big\|{D^2u}\big\|_{q_1}+\big\|{D^2u}\big\|_{q}\leq\dy
C(q_1,q)\, \mu^{\frac{ 2-p}{2} } (\|f\|_{q_1}+\|f\|_q) +C(q_1,q)\, (\|f\|_{q_1}+\|f\|_q) ^\frac{1}{p-1}\,,\ee
with
 $a=\frac{nq}{n(q-p)+pq}$.
Moreover
\be\label{gp3}\|\nabla u\|_p+\|\Pi(u,\chi^\rho)\|_p\leq \dy
C(q_1,q)\, \mu^{\frac{ 2-p}{2} } (\|f\|_{q_1}+\|f\|_q) +C(q_1,q)\, (\|f\|_{q_1}+\|f\|_q) ^\frac{1}{p-1}\,.\ee
Finally there exists a constant $C$, independent of $\eta$ and $\rho$, such that
     \be\label{trevv}
     \big\|D^3 u\big\|_{q}\leq C(\mu, \|f\|_{q_1}, \|f\|_{q}, \|\nabla f\|_{q})\,.
      \ee
}
\end{prop}
 \Pr  For $\eta>0$, let us denote by $\{u^\eta\}$ the sequence of solution of system \eqref{stokesS4}
 obtained in Proposition\,\ref{Vteo3}. These are also solutions of system \eqref{stokesS6},
 obtained by multiplying system \eqref{stokesS4} by $a(\mu,u)$. The bounds
  \eqref{passeqq1noee}  and \eqref{gp2} imply the existence of a field $u$
and of a subsequence, that for the sake of simplicity we do not relabel, such that,
$$\nabla u^\eta\rightharpoonup \nabla u \textrm{  in } L^p(\R^n),$$
$$D^2u^\eta\rightharpoonup D^2u \textrm{ in } L^{q_1}(\R^n)\cap L^{q}(\R^n), $$
and the limit function $u$ satisfies, from the lower
semi-continuity, estimate
\eqref{tpasseqq1} and \eqref{gp3}.
Moreover, from  Lemma \ref{converg1},
for any $\vp\in L^{q'}(\R^n)$,
\be\label{weakteta}
\lim_{\eta\to 0} (\frac{D^2 u^\eta}{a(\mu, u^\eta)}- \frac{D^2 u}{a(\mu, u)}, \vp)=0\,.\ee
Let us
 show
that the limit $u$ is actually a solution of system \eqref{stokesS5} in $\R^n$.
Below we argue as in the proof of Lemma \ref{converg2}. Let us consider the corresponding
subsequence of $\{{\mathcal A}_1(\mu, u^\eta)\}$. Recall that $J_\eta(\frac{\nabla\nabla u^\eta\, \chi^\rho}{a(\mu, u^\eta)})$ has compact support in $\R^n$, included in $\overline{B_{2\rho}}$. By using Rellich compactness
 theorem, and then that strong convergence implies almost everywhere convergence along a subsequence, there exists a subsequence such that
 \be\label{extr7}A^1(\mu, u^\eta)\to A^1(\mu, u) , \ \mbox{ a.e. in } \ B_{2\rho}\,.\ee
 We fix the last subsequence
 of $\{ A^1(\mu, u^\eta)J_\eta(\frac{\nabla\nabla u^\eta\chi^\rho}{a(\mu, v^\nu)})\}$.
 Let us set
\be\label{forpeta}\ba{ll}\vs1\dy
\!\!\big((A^1(\mu,\!u^\eta)\!-\!
A^1(\mu,\!u))\!\cdot\!J_\eta(\frac{\nabla\nabla u^\eta\chi^\rho}{a(\mu,\!u^\eta)}),\!\vp\big)\!+\!
\big(A^1(\mu,\!u)\!\cdot\!J_\eta\big(\frac{\nabla\nabla u^\eta\chi^\rho}{a(\mu,\!u^\eta)}\!-\!\frac{\nabla\nabla u\chi^\rho}{a(\mu,\!u)}\big),\!\vp\big)
\\\dy \hfill +\big(A^1(\mu, u) \cdot \big(J_\eta\big(\frac{\nabla\nabla u\,\chi^\rho}{a(\mu, u^)})-\frac{\nabla\nabla u\,\chi^\rho}{a(\mu, u)}\big), \vp\big) :=J_1^\eta+J_2^\eta+
J_3^\eta\,,
\ea\ee
with an arbitrary $\vp\in C_0^\infty(\R^n)$. One easily recognizes that $J_1^\eta\to 0$. Indeed, firstly note that
$|A^1(\mu, u^\eta) )-A^1(\mu, u) ) |\leq 2$. Then, using \eqref{extr7} and recalling that
  $J_\eta(\frac{\nabla\nabla u^\eta\,\chi^\rho}{a(\mu, u^\eta)})$ is bounded uniformly in $\eta$ thanks to  \eqref{passeunoee}, from the Lebesgue dominated convergence theorem $J_1^\eta$ goes to zero.
As far as the sequence of integrals $J_2^\eta$ is concerned, using Fubini's theorem we set
$$\ba{rl}\dy\vs1J_2^\eta=&\dy \big(\big(\frac{\nabla\nabla u^\eta\,\chi^\rho}{a(\mu, u^\eta)}-\frac{\nabla\nabla u\,\chi^\rho}{a(\mu, u)}\big), J_\eta(A^1(\mu, u)\cdot \vp\big)\big)\\\vs1 =&\dy\big(\big(\frac{\nabla\nabla u^\eta\,\chi^\rho}{a(\mu, u^\eta)}-\frac{\nabla\nabla u\,\chi^\rho}{a(\mu, u)}\big), J_\eta(A^1(\mu, u)
\cdot \vp\big)-A^1(\mu, u)
\cdot \vp\big)\\
&\dy +\big(\big(\frac{\nabla\nabla u^\eta\,\chi^\rho}{a(\mu, u^\eta)}-\frac{\nabla\nabla u\,\chi^\rho}{a(\mu, u)}\big), A^1(\mu, u)
\cdot \vp\big)= J_{21}^\eta+J_{22}^\eta\,.
\ea $$
Since $A^1(\mu, u)
\cdot \vp \in L^{q'}(\R^n) $, then $J_\eta(A^1(\mu, u)
\cdot \vp\big)\to A^1(\mu, u)
\cdot \vp$ strongly in $L^{q'}(\R^n)$, and
$J_{21}^\eta$ tends to zero
using the uniform bound \eqref{passeunoee}, or \eqref{passeqq1noee}.
$J_{22}^\eta$ tends to zero
by using the $L^q$ weak convergence
$\frac{D^2 u^\eta}{a(\mu, u^\eta)} $ to $\frac{D^2 u}{a(\mu, u)}$ given in \eqref{weakteta}.
The last sequence $J_3^\eta$ tends to zero from the strong convergence of $J_\eta\big(\frac{\nabla\nabla u\,\chi^\rho}{a(\mu, u^)})$ to
$\frac{\nabla\nabla u\,\chi^\rho}{a(\mu, u)}$, ensured by the validity of \eqref{tpasse} and the properties of the mollifier.
Therefore we have found a further subsequence such that
$$\lim_{\eta\to 0}\big(A^1(\mu, u^\eta)\cdot  J_\eta\big(\frac{\nabla\nabla u^\eta\,\chi^\rho}{a(\mu, u^\eta)}), \vp\big)=
 \big(A^1(\mu, u)\cdot  \frac{\nabla\nabla u\,\chi^\rho}{a(\mu, u)}, \vp\big).$$
Finally, let us consider the corresponding subsequence of $\{\nabla \Pi^\eta\}$, where $\Pi^\eta:=\Pi(u^\eta, \chi^\rho)$ is given by \eqref{pi}, with $u$ replaced by $u^\eta$.
 From estimate \eqref{gp2}, the sequence
 $\{\Pi^\eta\}$ is bounded in $L^p(\R^n)$.
 Moreover,  from \eqref{spres1}, $\nabla\Pi^\eta$  is uniformly bounded $L^{q_1}(\R^n)\cap L^q(\R^n)$.
Hence there exist a field $\Pi$
and a subsequence, that we still denote by $\Pi^\eta$, such that
$$\Pi^\eta\rightharpoonup  \Pi \textrm{ in } L^{p}(\R^n)\,,$$
$$\nabla\Pi^\eta\rightharpoonup \nabla\Pi \textrm{  in } L^{q_1}(\R^n)\cap L^q(\R^n), $$
and, from the lower
semi-continuity, estimate \eqref{spres2} and \eqref{gp3} are satisfied.
It remains to show that
$$\Pi=\Pi(u, \chi^\rho)
=(2-p)\int_{\R^n}\nabla_y  {\mathcal E} (x-y)
\, \frac{D_{y_i}u_j\,\nabla u}{a_1(\mu, u)}\cdot (\frac{D_{y_j}\nabla u\,\chi^\rho}{a(\mu, u)})dy\,. $$ Let $\vp\in C_0^\infty(\R^n)$
and set
$$g(y)= \int_{\R^n}\nabla_x {\mathcal E} (x-y)\, \vp(x)\,dx\,.$$
Then
$$\ba{rl}\vs1\dy
(\Pi^\eta,\vp)\dy=&\dy (p-2) \int_{\R^n}\frac{\nabla u_i^\eta   (y)D_{y_j}u_h^\eta(y) }{a_1(\mu, u^\eta(y))} J_\eta(\frac{D^2_{y_iy_j} u_h^\eta(y)\,\chi^\rho(y)}{a(\mu, u^\eta(y))}\big)\cdot g(y)\,dy\,,
\ea$$
and the convergence of a subsequence to $(\Pi,\vp)$ can be obtained repeating the arguments used in \eqref{forpeta}.
Therefore
$$(\frac{\Delta u}{a(\mu,u)}
+(p-2) A^1(\mu,u)\cdot \frac{\nabla\nabla u\,\chi^\rho}{a(\mu, u)}-\nabla \Pi(u, \chi^\rho)- f,\vp)=0\,,$$
for all $\vp\in C_0^\infty(\R^n)$, which ensures
that $u$ satisfies \eqref{stokesS5} a.e. in $\R^n$. Finally, the same subsequence satisfies the uniform bound \eqref{trev}. Hence, there exists a subsequence weakly converging to $D^3 u$ in $L^q(\R^n)$. By lower semicontinuity, $D^3 u$ satisfies the bound \eqref{trevv}.\chiu \vskip0.1cm
\section{\large The solution of problem \eqref{stokesS8} ($\rho\to \infty$)}\label{secrho}
We show that the solution of system \eqref{stokesS5}, in the limit as $\rho$ tends to infinity, converges to a solution of system
\eqref{stokesS8}. Subsequently we also prove that such a solution satisfies the equation $\nabla \cdot u=0$ in $\R^n$. This task is a consequence of a maximum principle on a suitable elliptic system.
\begin{prop}\label{Vteo8}{\sl Let $q\in (n,+\infty)$, $q_1=\frac{np}{n+p}$,
 and assume that $\ov M(2)$, $M(q_1)$ and $M(q)$ are positive constants.
  Let $f\in L^{q_1}(\R^n)\cap L^{q}(\R^n)$, with $\nabla f\in L^q(\R^n)$.
 Then, there exists a solution $u$ of system \eqref{stokesS8}, and there exist some
 constants C, independent of $\mu$, such that, for $r=q_1$ and $r=q$,
 \be\label{rhopasse}
     \big\|\frac{\nabla\nabla u}{a(\mu, u)}\big\|_{r}\leq C\, \|f\|_{r}\,, \ee
      \be\label{spres3} \|\nabla \Pi(u)\|_r\leq C(r)\, \|f\|_{r} \,,\ee
 \be\label{rhopasseqq1} \big\|{D^2u}\big\|_{q_1}+\big\|{D^2u}\big\|_{q}\leq\dy
C(q_1,q)\, \mu^{\frac{ 2-p}{2} } (\|f\|_{q_1}+\|f\|_q) +C(q_1,q)\, (\|f\|_{q_1}+\|f\|_q) ^\frac{1}{p-1}\,,\ee
where $\Pi(u)$ is given by \eqref{piu}.
Moreover,
\be\label{gp4}\|\nabla u\|_p+\|\Pi(u)\|_p\leq \dy
C(q_1,q)\, \mu^{\frac{ 2-p}{2} } (\|f\|_{q_1}+\|f\|_q) +C(q_1,q)\, (\|f\|_{q_1}+\|f\|_q) ^\frac{1}{p-1}\,.\ee
Finally there exists a constant $C$ such that
     \be\label{trevvv}
     \big\|D^3 u\big\|_{q}\leq C(\mu, \|f\|_{q_1}, \|f\|_{q}, \|\nabla f\|_{q})\,.
      \ee}
\end{prop}
\Pr
 For any $\rho>0$, let us denote by  $u^\rho$ the solution of system \eqref{stokesS5}, satisfying the properties of Proposition\,\ref{Vteo7}.
  Since, from \eqref{tpasseqq1}, $D^2 u^\rho$ is bounded in $L^{q_1}(\R^n)\cap L^q(\R^n)$, and,   from \eqref{gp3}, $\nabla u^\rho$ is bounded in $L^p(\R^n)$, there exist a
 a field $u$
and a subsequence, that we do not relabel, such that
 $$\nabla u^\rho\rightharpoonup \nabla u \textrm{  in } L^p(\R^n),$$
 $$D^2u^\rho\rightharpoonup D^2u \textrm{  in } L^{q_1}(\R^n)\cap L^{q}(\R^n),$$
 and estimate \eqref{rhopasseqq1}  and \eqref{gp4} hold.
 Further,   Lemma \ref{converg1} and Lemma \ref{converg3} in $E=\R^n$, applied either with $r=q$, $\ov r=q_1$ and  with $r=\ov r=q_1$, imply that, along a subsequence
\be\label{weakrho}
\frac{D^2 u^\rho}{a(\mu, u^\rho)}\rightharpoonup
\frac{D^2 u}{a(\mu, u)}\,,\textrm{  in } L^{q_1}(\R^n)\cap L^{q}(\R^n)\,,\ee
\be\label{weakrho1}
A^{\frac{4-p}{2}}_{ijhk}(\mu, u^\rho)D_{hk}^2 u_j^\rho
\rightharpoonup A^{\frac{4-p}{2}}_{ijhk}(\mu, u)D_{hk}^2 u_j\,,\textrm{ in } L^{q_1}(\R^n)\cap L^{q}(\R^n)\,,\ee
\be\label{weakrho2}
A^{\frac{4-p}{2}}_{ijhk}(\mu, u^\rho)D_{hk}^2 u_i^\rho
\rightharpoonup A^{\frac{4-p}{2}}_{ijhk}(\mu, u)D_{hk}^2 u_i\,,\textrm{  in } L^{q_1}(\R^n)\cap L^{q}(\R^n)\,,\ee
and estimate \eqref{rhopasse}
holds. Let us
 show
that the limit $u$ is actually a solution of system \eqref{stokesS8}.
Let  $\vp$ be an arbitrary function in
$C_0^\infty(\R^n)$. 
Clearly, from \eqref{weakrho}, $$ \lim_{\rho\to \infty}
 (\frac{\Delta u^\rho}{a(\mu, u^\rho)}- \frac{\Delta u}{a(\mu, u)}, \vp)=0\,.$$
 Let $R<\rho$. Then $\chi^\rho=1$ in $B_R$ and we set \be\label{forpra}\ba{ll}\dy\vs1
\big(\,A^\frac{4-p}{2}(\mu, u^\rho) \cdot \nabla\nabla u^\rho\,\chi^\rho-
A^\frac{4-p}{2}(\mu, u)\cdot \nabla\nabla u, \vp\big)\\\vs1 \dy =
\int_{B_R}\big(A^\frac{4-p}{2}(\mu, u^\rho) \cdot\nabla\nabla u^\rho-
A^\frac{4-p}{2}(\mu, u)\cdot \nabla\nabla u\big) \cdot \vp\, dx\\
\dy +\int_{\R^n-B_R}\big(A^\frac{4-p}{2}(\mu, u^\rho) \cdot\nabla\nabla u^\rho\,\chi^\rho-
A^\frac{4-p}{2}(\mu, u)\cdot \nabla\nabla u\big)\cdot \vp\, dx
\ea\ee
 The first term on the right-hand side tends to zero, uniformly in $R<\rho$,
  thanks to \eqref{weakrho1}. Then the second term on the right-hand side
  in \eqref{forpra} tends to zero, since the terms in the brackets are bounded in $L^{q_1}(\R^n)$, uniformly in $\rho$, and
  $\vp\in L^{q_1'}(\R^n)$ implies that $\|\vp\|_{L^{q_1'}(\R^n-B_R)}\to 0$ as $R\to \infty$.\par
 Finally, let us consider the corresponding subsequence of $\{\nabla \Pi^\rho\}$, where
$\Pi^\rho:=\Pi(u^\rho, \chi^\rho)$. From estimate \eqref{gp3}, $\{ \Pi^\rho\}$ is bounded in $L^{p}(\R^n)$, and, from estimate \eqref{spres2},  $\{\nabla  \Pi^\rho\}$ is bounded in $L^{q_1}(\R^n)\cap L^q(\R^n)$. Hence
there exist a field $\Pi$
and a subsequence such that
$$\Pi^\rho\rightharpoonup\Pi \textrm{  in } L^{p}(\R^n),$$
$$\nabla\Pi^\rho\rightharpoonup \nabla\Pi \textrm{  in } L^{q_1}(\R^n)\cap L^q(\R^n), $$
and the limit satisfies estimates \eqref{spres3} and \eqref{gp4}. We show that
$\Pi=\Pi(u)\,, $ with $\Pi(u)$ defined in \eqref{piu}.
It is enough to show this convergence on a subsequence, the convergence of the whole
sequence follows by uniqueness. Let $\vp\in C_0^\infty(\R^n)$. Setting
$$g(y)= \int\nabla_x\, {\mathcal E} (x-y)\, \vp(x)\,dx\,,$$
from the Hardy-Littlewood-Sobolev theorem $g\in L^s(\R^n)$, for $s\in (\frac{n}{n-1},\infty)$. Then, writing $(\Pi^\rho,\vp)$ as
$$
(\Pi^\rho,\vp)\dy=(p-2) \int_{\R^n}\frac{\nabla u_i^\rho (y)D_{y_j}u_h^\rho(y)}{a_1(\mu, u^\rho(y))} J_\eta(\frac{D^2_{y_iy_j} u_h^\rho(y)\,\chi^\rho(y)}{a(\mu, u^\rho(y))}\big)\cdot g(y)\,dy\,,$$
the convergence of a subsequence to $(\Pi(u),\vp)$ can be obtained reasoning as in \eqref{forpra}.
Therefore
\be\label{umusol1}(\frac{\Delta u}{a(\mu,u)}
+(p-2) A^{\frac{4-p}{2}}(\mu, u)\cdot {\nabla\nabla u}-
\nabla \Pi(u)-f,\vp)=0\,,\ee
for all $\vp\in C_0^\infty(\R^n)$. The existence of the third derivatives and estimate \eqref{trevvv} follow from lower semicontinuity  arguments.\chiu
\vskip0.1cm
\begin{prop}\label{Vteo9}{\sl Let the assumptions of Proposition\,\ref{Vteo8} be satisfied, and assume that $\nabla\cdot f=0$. Then $\nabla \cdot u=0$.}
\end{prop}
\Pr Let us multiply system \eqref{stokesS8} by $\nabla \vp$, for some $\vp\in C_0^\infty(\R^n)$, and integrate on $\R^n$:
\be\label{SS8}\big(\frac{\Delta u}{a(\mu,u)},\nabla\vp\big)
+\frac{(p-2)}{2} \big(\frac{\nabla u_j\cdot \nabla|\nabla u|^2}{a_{\frac{4-p}{2}}(\mu, u)}, D_{y_j}\vp)=
(\nabla \Pi(u),\nabla\vp)+( f,\nabla\vp).
\ee
Below we rewrite each term in \eqref{SS8} in a suitable way. Integrating by parts we find
\be\label{SS9}\big(\frac{\Delta u}{a(\mu,u)},\nabla\vp\big)=-
\big(\frac{\Delta(\nabla\cdot u)}{a(\mu,u)},\vp\big)+\frac{2-p}{2}\big(\frac{\Delta u\cdot \nabla|\nabla u|^2}{a_{\frac{4-p}{2}}(\mu,u)},\vp\big)\,.\ee
Integrating by parts,  
 we have
\be\label{SS10}\ba{ll}\vs1\dy\frac{(p-2)}{2} \big(\frac{\nabla u_j\cdot \nabla|\nabla u|^2}{a_{\frac{4-p}{2}}(\mu, u)}, D_{y_j}\vp)=
\frac{(2-p)}{2}\big(D_{y_j}\frac{\nabla u_j\cdot \nabla|\nabla u|^2}{a_{\frac{4-p}{2}}(\mu, u)}, \vp)\\
\vs1\dy=\frac{(2-p)}{2} \big(\frac{\nabla(\nabla\cdot u)\cdot \nabla|\nabla u|^2}{a_{\frac{4-p}{2}}(\mu, u)}, \vp)
+\frac{(2-p)}{2} \big(\frac{\nabla u_j\cdot D_{y_j}\nabla|\nabla u|^2}{a_{\frac{4-p}{2}}(\mu, u)}, \vp)\\
\dy \hfill+
\frac{(p-2)}{2}\frac{(4-p)}{2} \big(\frac{D_{y_i} u_j\cdot (D_{y_i}|\nabla u|^2D_{y_j}|\nabla u|^2)}{a_{\frac{6-p}{2}}(\mu, u)}, \vp).
\ea\ee
As far as the term with $\nabla\Pi$, we argue as follows. Recall
that $$
\Pi(u):=\frac{(2-p)}{2}\int_{\R^n}D_{y_i}  {\mathcal E} (x-y)
\frac{(D_{y_i}u_j\,D_{y_j} |\nabla u|^2)(y)}{a_{\frac{4-p}{2}}(\mu, u(y))}\,dy\,.$$
Then
\be\label{SS11}\ba{ll} \vs1\dy
(\nabla\Pi(u),\nabla \vp)=-(\Delta \Pi(u),\vp)=
\frac{(2-p)}{2}(D_{y_i} \big(\frac{D_{y_i}u_j\, D_{y_j}|\nabla u|^2}{a_{\frac{4-p}{2}}(\mu, u)} \big),\vp)\\
\vs1\dy=\frac{(2-p)}{2}(\frac{\Delta u_j D_{y_j}|\nabla u|^2}{a_{\frac{4-p}{2}}(\mu, u)},\vp)+
\frac{(2-p)}{2}(\frac{D_{y_i} u_j D^2_{y_iy_j}|\nabla u|^2}{a_{\frac{4-p}{2}}(\mu, u)} ,\vp)\\
\dy \hfill
+\frac{(p-2)}{2}\frac{(4-p)}{2}(\frac{D_{y_i}u_j\, D_{y_i}|\nabla u|^2\, D_{y_j}|\nabla u|^2}{a_{\frac{6-p}{2}}(\mu, u)} ,\vp)\,.
\ea\ee
Therefore, using \eqref{SS9}, \eqref{SS10} and \eqref{SS11} in \eqref{SS8},
we obtain that, for any $\vp\in C_0^\infty(\R^n)$,
$$ \big(\frac{\Delta(\nabla\cdot u)}{a(\mu,u)}+\frac{(p-2)}{2} \frac{\nabla(\nabla\cdot u)\cdot \nabla|\nabla u|^2}{a_{\frac{4-p}{2}}(\mu, u)}, \vp\big)=(f,\nabla\vp)\,,$$
which, from the assumption on $f$ and the integrability of $u$, implies that
\be\label{modulus}
\frac{\Delta(\nabla\cdot u)}{a(\mu,u)}+ \frac{(p-2)}{2}\frac{\nabla(\nabla\cdot u)\cdot \nabla|\nabla u|^2}{a_{\frac{4-p}{2}}(\mu, u)}=0\,.\ee
Let us set $U:=\nabla\cdot u$. Then $U$ is solution of the following elliptic system
\be\label{modulus1}
\Delta U+ \frac{(p-2)}{2} \frac{\nabla U\cdot \nabla|\nabla u|^2}{a_1(\mu, u)}=0\,.\ee
We recall that $\frac{\nabla|\nabla u|^2}{a_{1}(\mu, u)}\in L^\infty(\R^n)$
  and $D^2 U \in L^{q}(\R^n)$. Therefore we can apply the maximum principle for elliptic equations (see \cite{GT}, chap. 9) on each $B_{\sigma}\subset \R^n$. Since
  $U\to 0$ as $\sigma\to \infty$, we get $U\equiv 0$.\chiu

\section{\large Proof of Theorem\,\ref{mainTT}}\label{proof}
  \vskip0.2cm{\it Step I: \underline{The case $\mu>0$ and $f\in J^q(\R^n)\cap J^{q_1}(\R^n)$}} -
 Let us consider the sequence  $\{f^\nu\}\in {\mathscr C}_0(\R^n)$ converging to $f$ in $J^{q_1}(\R^n)\cap J^{q}(\R^n)$,
whose existence is ensured by Lemma\,\ref{density}.
 This sequence satisfies the assumptions of Proposition\,\ref{Vteo8}  and Proposition\,\ref{Vteo9}.
Hence, there exists a sequence $\{u^\nu\}\in J^\frac{np}{n-p}(\R^n)$ of solutions
of system \eqref{stokesS8}, with $\Pi(u^\nu)$ given by \eqref{piu}, with $\nabla\cdot u^\nu=0$, and satisfying estimates \eqref{rhopasse}--\eqref{trevvv}, uniformly in $\nu>0$.
  Since $\{D^2 u^\nu\}$ is bounded in $L^{q_1}(\R^n)\cap L^q(\R^n)$, and $\{\nabla u^\nu\}$ is bounded in $L^p(\R^n)$
  there exist a field $u$ and a (not relabeled) subsequence such that
$$\nabla u^\nu\rightharpoonup \nabla u \textrm{ in } L^p(\R^n),$$
  $$D^2u^\nu\rightharpoonup D^2u \textrm{ in } L^{q_1}(\R^n)\cap L^{q}(\R^n), $$
and \eqref{mainETT1} holds for the limit function $u$.
  Moreover
  Lemma \ref{converg1} and Lemma \ref{converg3} in $E=\R^n$,
applied either   with $r=q$, $\ov r=q_1$ and   with $r=\ov r=q_1$, imply the existence of a subsequence, that we do not relabel, such that
\be\label{weaknu} \frac{D^2
u^\nu}{a(\mu, u^\nu)}
\rightharpoonup \frac{D^2 u}{a(\mu,
u)}, \textrm{ in }
L^{q_1}(\R^n)\cap L^{q}(\R^n)\,,\ee
\be\label{weaknu1}
A^{\frac{4-p}{2}}_{ijhk}(\mu,
u^\nu)D_{hk}^2 u_j^\nu
\rightharpoonup
A^{\frac{4-p}{2}}_{ijhk}(\mu,
u)D_{hk}^2 u_j\,,\textrm{  in }
L^{q_1}(\R^n)\cap L^{q}(\R^n)
\,,\ee \be\label{weaknu2}
A^{\frac{4-p}{2}}_{ijhk}(\mu,
u^\nu)D_{hk}^2 u_i^\nu
\rightharpoonup
A^{\frac{4-p}{2}}_{ijhk}(\mu,
u)D_{hk}^2 u_i\,, \textrm{ in }
L^{q_1}(\R^n)\cap L^{q}(\R^n)\,.\ee
Let us consider the corresponding
subsequence of $\{\nabla
\Pi^\nu\}$, where $ \Pi^\nu:=
\Pi(u^\nu)$, with $\Pi(u^\nu)$
defined in \eqref{piu}. From
\eqref{gp4} $\{\Pi^\nu\}$ is
bounded in $L^p(\R^n)$, while from
\eqref{spres3}
$\{\nabla\Pi^\nu\}$ is uniformly bounded in $L^{q_1}(\R^n)\cap L^q(\R^n)$. Hence
there exist a field $\Pi$
and a subsequence, that we still denote by $\Pi^\nu$, such that
$$\Pi^\nu\rightharpoonup\Pi \textrm{  in } L^{p}(\R^n),$$
$$\nabla\Pi^\nu\rightharpoonup \nabla\Pi \textrm{  in } L^{q_1}(\R^n)\cap L^q(\R^n), $$
with $\nabla\Pi$ satisfying estimate \eqref{mainETT2} and, from \eqref{gp4} and the lower semi-continuity,
\be\label{piuip}
\|\Pi\|_p\leq \dy
C(q_1,q)\, \mu^{\frac{ 2-p}{2} } (\|f\|_{q_1}+\|f\|_q) +C(q_1,q)\, (\|f\|_{q_1}+\|f\|_q) ^\frac{1}{p-1}\,.\ee
We want to show that
$\Pi=\Pi(u)\,, $ where $\Pi(u)$ is defined in \eqref{piu}.
It is enough to show this convergence on a subsequence, and the convergence of the whole
sequence follows by uniqueness. Let $\vp\in C_0^\infty(\R^n)$ and set
$$g(y)= \int\nabla_x\, {\mathcal E} (x-y)\, \vp(x)\,dx\,.$$
We can write
$$(\Pi^\nu,\vp)\dy=(p-2) \int_{\R^n}\frac{\nabla u_i^\nu (y)D_{y_j}u_h^\nu(y)}{a_1(\mu, u^\nu(y))} (\frac{D^2_{y_iy_j} u_h^\nu(y)}
{a(\mu, u^\nu(y))}\big)\cdot g(y)\,dy\,.$$
Now the convergence of a subsequence to $(\Pi(u),\vp)$ follows from  \eqref{weaknu2}, observing that, from the Hardy-Littlewood-Sobolev theorem, $g\in L^{q_1'}(\R^n)$.
Therefore
\be\label{umusol2}(\frac{\Delta u}{a(\mu,u)}
+(p-2) A^{\frac{4-p}{2}}(\mu, u)\cdot {\nabla\nabla u}-
\nabla \Pi(u)-f,\vp)=0\,,\ee
for all $\vp\in C_0^\infty(\R^n)$. The existence of the third derivatives and estimate \eqref{trevvv} follow from lower semicontinuity  arguments.
By employing the divergence free condition satisfied by the sequence $\{u^\nu\}$,
and the convergence of $(\nabla\cdot u^\nu)$ to $(\nabla\cdot u)$ in $\R^n$,
we obtain that $u$ is a high regular solution of \eqref{stokes}, in the sense of Definition\,\ref{DS}.
  \vskip0.4cm\noindent{\it Step II: \underline{The case $\mu=0$ and $f\in J^q(\R^n)\cap J^{q_1}(\R^n)$}} -
For
$\mu>0$, let $\{u^\mu\}$ be the sequence of
solutions of \eqref{stokesS8} obtained in Step I.
We have to
show that, as $\mu\to 0$,  the sequence $\{u^\mu\}$ converges, in suitable norms, to a limit function  $u$ which is solution of
problem \eqref{stokes} with $\mu=0$. Actually we show that, \be\label{lob3}\nabla\cdot u(x)=0\,,\ \mbox{ in } \ \R^n\,,\ee and, for any $\vp\in W^{1,p}(\R^n)$,
\be\label{lob4}
(|\nabla u|^{p-2} \nabla u,\nabla \vp)=-(\nabla\Pi(u)+f, \vp)\,,\ee
with
\be\label{pif}
\Pi(u):=\int_{\R^n}^*D^2_{y_iy_j}{\mathcal E}(x-y) \,\frac{D_{y_i} u_j}{|\nabla u|}{\atop^{\!\!\!2-p}}\,\,dy\,,\ee
where $\int_{\R^n}^*G(y)dy$
 denotes the principal value singular integral in the Cauchy sense.
\par Firstly note that, from the  $\mu$-uniform bound \eqref{mainETT1} obtained in the first part of the proof,
\be\label{loba}
\big\|\frac{\nabla u^\mu}{a(\mu,u^\mu)}\big\|_{p'}\leq \|\nabla u^\mu\|_p<c\|D^2 u^\mu\|_{q_1}<+\infty\,.\ee
Therefore, there exist two fields, $u$ and $\Psi$, and a subsequence such that, in the limit as $\mu\to 0$,
\be\label{loba1}\ba{ll}
\vs1\dy u^\mu  \rightharpoonup  u \mbox{ in } L^{\frac{np}{n-p}}(\R^n)\,,\\\vs1\dy
\dy\nabla u^\mu  \rightharpoonup \nabla u \mbox{ in } L^{p}(\R^n)\,, \\
\dy\vs1\,D^2 u^\mu\rightharpoonup D^2u \mbox{ in } L^{q_1}(\R^n)\cap L^q(\R^n)\,,\\\dy
\frac{\nabla u^\mu}{a(\mu,u^\mu)} \rightharpoonup  \Psi \ \mbox{ in \ } L^{p'}(\R^n)\,.
\ea\ee
Moreover, $D^2u$ satisfies estimate
 \be\label{mupasseqq1} \big\|{D^2u}\big\|_{q_1}+\big\|{D^2u}\big\|_{q}\leq\dy
C(q_1,q)\, (\|f\|_{q_1}+\|f\|_q) ^\frac{1}{p-1}\,.\ee
By Rellich's theorem, for any compact set $K\subset \R^n$ there exists a subsequence, depending on $K$,  such that
$$\nabla u^\mu  \rightarrow \nabla u \mbox{  in } L^{p}(K)\,.$$
This last convergence implies the almost everywhere convergence in $K$. Therefore, we also have
$$\frac{\nabla u^\mu}{a(\mu,u^\mu)} \rightarrow  \frac{\nabla u}{|\nabla u|}{\atop^{2-p}} \ \mbox{ a.e. in \ } K\,. $$
Since, from \eqref{loba1}$_4$, this subsequence weakly converges  to $\Psi$ in $L^{p'}(K)$,
we find that $\Psi=\frac{\nabla u}{|\nabla u|}{\atop^{2-p}}$, on each compact $K\subset \R^n$, which ensures that
\be\label{lob2}
\frac{\nabla u^\mu}{a(\mu,u^\mu)} \rightharpoonup  \frac{\nabla u}{|\nabla u|}{\atop^{2-p}}
 \ \mbox{ in \ } L^{p'}(\R^n) \mbox{ as } \mu\to 0\,.
\ee
Finally, let us consider the corresponding subsequence of $\{\nabla \Pi^\mu\}$, where $\Pi^\mu:= \Pi(u^\mu)$. From the validity of \eqref{mainETT2} for $\mu>0$,  $\{\nabla \Pi^\mu\}$ is uniformly bounded in $L^{q_1}(\R^n)\cap L^q(\R^n)$. Moreover, from \eqref{piuip},  $\{\Pi^\mu\}$ is uniformly bounded in $L^{p}(\R^n)$.
Hence
there exist a field $\Pi$
and a subsequence such that
\be\label{sn1}\ba{ll}\vs1\dy
\hskip0.5cm \Pi^\mu\rightharpoonup\Pi \textrm{ in } L^{p}(\R^n),\\
\dy\nabla\Pi^\mu\rightharpoonup \nabla\Pi \textrm{ in } L^{q_1}(\R^n)\cap L^q(\R^n), \ea\ee
with the limit satisfying \eqref{mainETT2}.
These arguments ensure that identity \eqref{lob4} holds, with some function $\Pi$ on the right-hand side.
 We want to show that
$\Pi=\Pi(u)\,, $ where $\Pi(u)$ is given by \eqref{pif}.
We prove that in any compact $\ov{B_{R_\circ}}$ the family $\{\Pi^\mu\}$ weakly converges to $\Pi(u)$. By uniqueness of the weak limit, we get $\Pi=\Pi(u)$.
To this end we remark that, by virtue of \eqref{loba} and \eqref{mupasseqq1}, we can deduce that, for all $R>0$, $ \frac{\nabla u^\mu}{a(\mu,u^\mu)}$ strongly converges to $ \frac{\nabla u}{|\nabla u|}{\atop^{2-p}}
$ in $L^{p'}(B_{R_\circ})$ and  $ |\frac{\nabla u^\mu}{a(\mu,u^\mu)}|+| \frac{\nabla u}{|\nabla u|}{\atop^{2-p}}|\leq C$, for all $\mu>0$. Hence, for all $\vp\in L^{p'}(B_{R_\circ})$ we easily obtain, for all $R>R_\circ$
\be\label{pM}\ba{ll}\vs1\dy (\Pi^\mu-\Pi(u), \vp)=\dy\int_{\R^n}\int_{\R^n}\!D^2_{y_iy_j} {\mathcal E}(x-y)\, (\frac{D_{y_i}u_j^\mu}{a(\mu, u^\mu)}-\frac{D_{y_i}u_j}{|\nabla u|}{\atop^{\!2-p}})\,dy\, \vp(x)\, dx\\
\hfill= \dy \int_{B_{R_\circ}}\int_{B_{2R}}\!\!D^2_{y_iy_j} {\mathcal E}(x-y)\, (\frac{D_{y_i}u_j^\mu}{a(\mu, u^\mu)}-\frac{D_{y_i}u_j}{|\nabla u|}{\atop^{\!2-p}})\,dy\, \vp(x)\, dx\\\hfill\dy + \int_{B_{R_\circ}}\int_{\R^n-B_{2R}}\!\!\!\!D^2_{y_iy_j} {\mathcal E}(x-y)\, (\frac{D_{y_i}u_j^\mu}{a(\mu, u^\mu)}-\frac{D_{y_i}u_j}{|\nabla u|}{\atop^{\!2-p}})\,dy\, \vp(x)\, dx:=I_1+I_2\,.
\ea\ee
Applying H\o lder's inequality, we get
$$|I_1|\leq \|\vp\|_{L^p(B_{R_\circ})} \| \frac{D_{y_i}u_j^\mu}{a(\mu, u^\mu)}-\frac{D_{y_i}u_j}{|\nabla u|}{\atop^{\!2-p}}\|_{L^{p'}(B_{2R})}\,,$$
$$\ba{ll}\dy |I_2|\leq \|\vp\|_{L^1(B_{R_\circ})} (\int_{|x-y|>R}\frac{dy}{|x-y|}{\atop^{\!np}})^\frac 1p(\| \frac{D_{y_i}u_j^\mu}{a(\mu, u^\mu)}\|_{p'}+\|\frac{D_{y_i}u_j}{|\nabla u|}{\atop^{\!2-p}}\|_{p'})\\
\dy \hfill \leq \|\vp\|_{L^1(B_{R_\circ})} R^{-\frac{n}{p'}}C\,.\ea$$
Hence in the limit as $\mu\to 0$ we get
$$\lim_{\mu\to 0} | (\Pi^\mu-\Pi(u), \vp)|\leq  \|\vp\|_{L^1(B_{R_\circ})} R^{-\frac{n}{p'}}\,.$$ Since $R$ is arbitrary, we get that $$\lim_{\mu\to 0}  (\Pi^\mu-\Pi(u), \vp)=0\,, \ \forall \vp \in L^{p'}(B_{R_\circ})\,,$$
which implies the thesis.
\par Finally, by employing the divergence free condition satisfied by the sequence $\{u^\mu\}$,
and the convergence of $(\nabla\cdot u^\mu)$ to $(\nabla\cdot u)$ in $\R^n$,
the validity of \eqref{lob3} also follows.
\vskip0.3cm\noindent{\it Step III: \underline{The case $\mu\geq 0$ and  $f\in L^q(\R^n)\cap L^{q_1}(\R^n)$}} - Let us consider the Helmholtz decomposition of $f$:
\be\label{pression1}
f=F+\nabla \psi\,,\ee
with $\psi\in \widehat W^{1,q}(\R^n)\cap \widehat W^{1,q_1}(\R^n)$ solution of the Poisson equation
\be\label{pression2}
\Delta\psi=\nabla\cdot f\,, \ee
and $F=P_q f=P_{q_1}f\in J^q(\R^n)\cap J^{q_1}(\R^n)$, where $P_q$ and $P_{q_1}$ are the projection operators in $J^q(\R^n)$ and $J^{q_1}(\R^n)$, respectively.
 From the results of the previous steps, with the right-hand side $f$ now replaced by  $F$, we find a solution $(u,\widetilde\pi)$ of
\eqref{stokes}
with \be\label{pression3}
\widetilde\pi=(2-p)\int_{\R^n}D_{y_i}  {\mathcal E} (x-y)
\frac{(D_{y_i}u_j\, \nabla u)(y)}{a_{\frac{4-p}{2}}(\mu, u(y))}\cdot D_{y_j}\nabla u(y)\,dy\,, \mbox{ for } \mu>0\,,\ee
\be\label{pression3M}
\widetilde\pi=\int_{\R^n}^*D^2_{y_iy_j}{\mathcal E}(x-y) \,\frac{D_{y_i} u_j}{|\nabla u|}{\atop^{\!\!\!2-p}}\,\,dy\,,\ \mbox{ for } \mu=0\,.\ee
Clearly the pair $(u, \pi)$, with
$\pi=\widetilde\pi+\psi$ is a solution of problem \eqref{stokes} in the sense of Definition\,\ref{DS}, with $f\in L^q(\R^n)\cap L^{q_1}(\R^n)$, and satisfies estimates \eqref{mainETT1} and \eqref{mainETT2}. In this connection, we observe that, from Step II, $\widetilde \pi\in L^{p'}(\R^n)$ by interpolation, while $\psi\in L^{p'}(\R^n)$, by using the representation of solutions to the Poisson equation and then applying the Hardy-Littlewood-Sobolev theorem with $f\in L^{q_1}(\R^n)$.  \chiu
\vskip0.2cm \noindent{\it Step IV: \underline{Uniqueness}} -
The solution $u$, obtained in the first part of the proof, is also a weak solution of \eqref{stokes}, in the sense of
Definition\,\ref{DWS}.
Let $v$ be another weak solution. Then  $u-v\in \widehat J^{1,p}(\R^n)$ can be approximated by a sequence $\{\vp^k\}\in {\mathscr C}_0(\R^n)$. Hence
\be\label{ngiu}\int ({\mathbb S}(\nabla u)-{\mathbb S}(\nabla v))
\cdot \nabla \varphi^k
\,dx=0\,.\ee
Passing to the limit as $k\to\infty$, and then applying Lemma\,\ref{giu} we find
$$\ba{rl}\dy 0\!\!&\dy=\lim_{k\to\infty}\int ({\mathbb S}(\nabla u)-{\mathbb S}(\nabla v))
\cdot \nabla \varphi^k
\,dx\\ \dy &\dy =\int ({\mathbb S}(\nabla u)-{\mathbb S}(\nabla v))
\cdot \nabla (u-v)
\,dx\geq \int C (\mu+ |\nabla u|^2 + |\nabla v|^2)^\frac{p-2}{2}|\nabla (u-v)|^2\,dx,\ea $$
which ensures that $u=v$ a.e. in $\R^n$, and completes the proof.  \chiu

\section{\large Proof of Theorem\,\ref{maindual}}\label{proof1}
We start by proving the following lemma, which is an easy consequence of Theorem\,\ref{mainTT}.
\begin{lemma}\label{CmainT}{\sl Let $\mu>0$, and let the assumptions of Theorem\,\ref{mainTT} be satisfied.  Then the following estimate holds
\be\label{d2lq}
\|D^2 u\|_q\leq c\|f\|_q\|a(\mu,u)\|_\infty^{2-p}.\ee }
\end{lemma}
\Pr Since $\mu>0$, we can write system \eqref{stokes} as follows
\be\label{stokesSa}\frac{\Delta u}{a(\mu,u)}
+(p-2) \frac{(\nabla u\otimes \nabla u)\cdot \nabla\nabla u}{a_{\frac{4-p}{2}}(\mu, u)}=
\nabla \pi+ f\,,\ \ \nabla\cdot u=0,\ \mbox{ in } \R^n\,.
\ee
Let us multiply system \eqref{stokesSa} by $a(\mu, u) \Delta u |\Delta u|^{q-2}$,  and integrate over $\R^n$.
By H\o lder' s inequality
 we get
$$\ba{rl}\dy\vs1  \big\|{\Delta u}\big\|_q ^{q}\leq&\dy \!
(2-p) \big\|\nabla\nabla u\big\|_{q}
 \big\|\Delta u\big\|_q^{q-1} \\&\dy +\|\nabla \pi\|_{q}
  \big\|\Delta u\big\|_q^{q-1}\|a(\mu, u)\big\|_{\infty} +
 \|f\|_{q}  \big\|\Delta u\|_q^{q-1}\|a(\mu, u)\|_{\infty},
    \ea$$
    which, by using estimates
   $$
 \|D^2 u\|_{q} \leq H (q)\|\Delta u\|_{q}\,,$$
and \eqref{mainETT2} on the pressure gradient,  easily gives
    $$\big\|{D^2 u}\big\|_q \leq\dy  c\|f\|_{q} \|a(\mu, u)\|_{\infty}.$$ \chiu
    \vskip0.1cm
    In the proof of Theorem\,\ref{mainTT}, either for the regularity results and for
    the uniqueness result, it is crucial to have $\nabla u\in L^p(\R^n)$. In Theorem\,\ref{mainTT} it was not possible to
    deduce this result from an {\it energy estimate}, and it was deduced from $D^2u\in L^{q_1}(\R^n)$.
    In Theorem\,\ref{maindual} we want to relax the assumptions on $f$, by requiring only
    $f\in L^{q_2}(\R^n)\cap L^q(\R^n)$. Then we have not anymore at disposal $D^2 u\in L^{q_1}(\R^n)$ to get $\nabla u\in L^p(\R^n)$, and
    we use a {\it weak form of the energy estimate} (see \eqref{c1}). Hence we have to develop a suitable construction, which makes the proof involved.
 \vskip0.1cm\noindent
 {\it Proof of Theorem\,\ref{maindual}} -
For some fixed $\zeta> 0$, let $\chi^\zeta(x)$ be a smooth nonnegative cut-off function, with $|\nabla \chi^\zeta|\leq c\, \zeta^{-1}$. Let us consider the following system
\be\label{stokeschi}
\nabla \cdot( (\mu+|\nabla u|^2)^\frac{p-2}{2} \,\nabla u )-\nabla\pi= f\chi^\zeta\,,\quad \nabla \cdot u=0\ \mbox{ in } \R^n\,,\ n\geq 3\,,\ee
with $\mu> 0$.  As $f\chi^\zeta\in L^q(\R^n)\cap L^{q_1}(\R^n)$, from Theorem\,\ref{mainTT}, there exists a unique high regular solution $(u^{\mu,\zeta}, \pi^{\mu,\zeta})$.    Let us set
$$\R^n_<=\{x\in \R^n: |\nabla u^{\mu,\zeta}(x)|^2<\mu\}\,, $$
$$\R^n_>=\{x\in \R^n: |\nabla u^{\mu,\zeta}(x)|^2\geq \mu\}\,.$$
These sets are well defined, due to the continuity of $\nabla u^{\mu,\zeta}$ ensured by Theorem\,\ref{mainTT}.
Then
\be\label{c0a}
\int_{\R^n_<} a(\mu, u^{\mu,\zeta})|\nabla  u^{\mu,\zeta}|^2dx\geq (2\mu)^{\frac{p-2}{2}}\int_{\R^n_<} |\nabla  u^{\mu,\zeta}|^2dx\,,\ee
and
\be\label{c0b}
\int_{\R^n_>} a(\mu, u^{\mu,\zeta})|\nabla  u^{\mu,\zeta}|^2dx\geq 2^{\frac{p-2}{2}}\int_{\R^n_>} |\nabla  u^{\mu,\zeta}|^pdx\,.\ee
Multiplying both sides of \eqref{stokeschi} by $u^{\mu,\zeta}$,  integrating by parts,
 using that $u^{\mu,\zeta}$ is divergence free, and Sobolev's embedding we have
\be\label{c1}\int_{\R^n}\!\!a(\mu, u^{\mu,\zeta})|\nabla  u^{\mu,\zeta}|^2dx\leq \|f\|_{q_2} \|\chi^\zeta u^{\mu,\zeta}\|_\frac{np}{n-p}\leq c
 \|f\|_{q_2} \|\nabla (\chi^\zeta u^{\mu,\zeta})\|_p\,,\ee
 where $q_2= \frac{np}{np-n-p}$.
 Let us estimate the $L^p$-norm of the right-hand side:
\be\label{c2}\ba{rl}\vs1\dy\| \nabla (\chi^\zeta u^{\mu,\zeta})\|_p^p\leq &\!\!\!\!  \dy c \int_{K(\zeta)} \frac{|u^{\mu,\zeta}|^p}{\zeta^p}\,dx+\int (\chi^\zeta)^p|\nabla u^{\mu,\zeta}|^p\,dx
 \\ \leq&\dy\!\!\!\!  c \int_{K(\zeta)} \frac{|u^{\mu,\zeta}|^p}{\zeta^p}+\int_{\R^n_>} (\chi^\zeta)^p|\nabla u^{\mu,\zeta}|^pdx+
 \int_{\R^n_<} (\chi^\zeta)^p|\nabla u^{\mu,\zeta}|^pdx
  \\ \dy\leq &\!\!\!\! \dy c \int_{K(\zeta)} \frac{|u^{\mu,\zeta}|^p}{\zeta^p}+\int_{\R^n_>} |\nabla u^{\mu,\zeta}|^pdx+
 \mu^\frac p2 (2\zeta)^{n}\,.
  \ea\ee
  with $K(\zeta)=\{  \zeta\leq |x|\leq 2\zeta  \}$.
 On the other hand, observing that \eqref{c0a} and \eqref{c0b} imply
 $$\ba{rl}\vs1
 \dy \int_{\R^n}\!\!a(\mu, u^{\mu,\zeta})|\nabla  u^{\mu,\zeta}|^2dx\geq &\dy\frac 12 \int_{\R^n} \!\!a(\mu, u^{\mu,\zeta})|\nabla  u^{\mu,\zeta}|^2dx \\ &\dy+
2^{\frac{p-4}{2}}\mu^{\frac{p-2}{2}}\int_{\R^n_<} |\nabla  u^{\mu,\zeta}|^2dx+
2^{\frac{p-4}{2}}\int_{\R^n_>} |\nabla  u^{\mu,\zeta}|^pdx\,,\ea$$
from \eqref{c1} and \eqref{c2}, by applying Young's inequality  we find
$$\ba{ll}\dy
\frac 12 \int_{\R^n} \!\!a(\mu, u^{\mu,\zeta})|\nabla  u^{\mu,\zeta}|^2dx+
2^{\frac{p-4}{2}}\mu^{\frac{p-2}{2}}\int_{\R^n_<} |\nabla  u^{\mu,\zeta}|^2dx+
2^{\frac{p-4}{2}}\int_{\R^n_>} |\nabla  u^{\mu,\zeta}|^pdx\\
\dy\hfill \leq
c
 \|f\|_{q_2}\left(\|\frac{u^{\mu,\zeta}}{\zeta}\|_{L^p(K(\zeta))}+ \mu^\frac 12 (2\zeta)^\frac np \right)
+ \ve \int_{\R^n_>} |\nabla u^{\mu,\zeta}|^pdx+ c(\ve) \|f\|_{q_2}^{p'}\,,
\ea$$
which finally gives
\be\label{c3}\ba{rl}\vs1\dy
\int_{\R^n} \!\!a(\mu, u^{\mu,\zeta})|\nabla  u^{\mu,\zeta}|^2dx\leq &\dy \!\!\!c\|f\|_{q_2}
\left(\|\frac{u^{\mu,\zeta}}{\zeta}\|_{L^p(K(\zeta))}+ \mu^\frac 12 (2\zeta)^\frac np\! \right)+c\|f\|_{q_2}^{p'}\\ :=&\dy \!\!
c\,\|f\|_{q_2} B(\mu, \zeta)+c\, \|f\|_{q_2}^{p'}\,.\ea\ee
A straightforward calculation, together with Gagliardo-Nirenberg's inequality give
\be\label{c4}\ba{ll}\dy \vs1
|\nabla u^{\mu,\zeta}(x)|^p=\frac{|\nabla u^{\mu,\zeta}(x)|^p\,|\nabla u^{\mu,\zeta}(x)|^{2-p}}
{(\frac 12|\nabla u^{\mu,\zeta}(x)|^2+\frac 12 |\nabla u^{\mu,\zeta}(x)|^2)}{\atop^{\frac{2-p}{2}}}
\dy \leq \frac{2^{\frac{2-p}{2}}|\nabla u^{\mu,\zeta}(x)|^2}{(\mu+|\nabla u^{\mu,\zeta}(x)|^2)}{\atop^{\frac{2-p}{2}}}
\\ \dy \leq c\,2^{\frac{2-p}{2}}
\|\nabla(\frac{|\nabla u^{\mu,\zeta}|^2}{(\mu+|\nabla u^{\mu,\zeta}|^2)}{\atop^{\frac{2-p}{2}}})\|_q^a
\|\frac{|\nabla u^{\mu,\zeta}|^2}{(\mu+|\nabla u^{\mu,\zeta}|^2)}{\atop^{\frac{2-p}{2}}}\|_1^{1-a}\,,\ \forall x\in \R^n_>\,,\ea\ee
with $a=\frac{nq}{nq+q-n}\,$.
For the $L^q$-norm in \eqref{c4}, by using Lemma\,\ref{CmainT} for the second derivatives of $u^{\mu,\zeta}$ we get
\be\label{c5}\ba{ll}\dy\vs1
\|\nabla(\frac{|\nabla u^{\mu,\zeta}|^2}{(\mu+|\nabla u^{\mu,\zeta}|^2)}{\atop^{\frac{2-p}{2}}})\|_q\leq c (\int |D^2 u^{\mu,\zeta}(x)|^q
|\nabla u^{\mu,\zeta}(x)|^{(p-1)q}dx)^\frac 1q\\
\dy\hfill \leq
 c\|\nabla u^{\mu,\zeta}\|_\infty^{p-1} \|f\chi^\zeta\|_{q} (\mu^\frac{2-p}{2}+\|\nabla u^{\mu,\zeta}\|_{\infty}^{2-p}).  \ea\ee
Therefore, taking into account that the $L^1$-norm in \eqref{c4} can be estimated via \eqref{c3},
and recalling that $\|\nabla u^{\mu, \zeta}\|_{
\infty}\leq \|\nabla u^{\mu, \zeta}\|_{L^{\infty}(\R^n_{>})}+\mu^\frac 12$, estimate \eqref{c4} gives
\be\label{c6}\ba{ll}\dy \vs1
|\nabla u^{\mu,\zeta}(x)|^p\leq c \, (\mu^\frac{(2-p)a}{2}\|\nabla u^{\mu,\zeta}\|_\infty^{(p-1)a}  +
 \|\nabla u^{\mu,\zeta}\|_\infty^{a}) \|f\|_{q}^a
\|f\|_{q_2}^{1-a}
B(\mu, \zeta)^{1-a}\\\hfill
\vs1\dy +c (\mu^\frac{(2-p)a}{2}\|\nabla u^{\mu,\zeta}\|_\infty^{(p-1)a}  +
 \|\nabla u^{\mu,\zeta}\|_\infty^{a}) \|f\|_{q}^a  \|f\|_{q_2}^{p'(1-a)}\\\vs1
 \hfill \dy \leq c \, (\mu^\frac{a}{2}+\mu^\frac{(2-p)a}{2}\|\nabla u^{\mu,\zeta}\|_{L^\infty(\R^n_>)}^{(p-1)a}  +
\|\nabla u^{\mu,\zeta}\|_{L^\infty(\R^n_>)}^{a}
) \|f\|_{q}^a
\|f\|_{q_2}^{1-a}
B(\mu, \zeta)^{1-a}\\\dy
+ c \,(\mu^\frac{a}{2}\!+\mu^\frac{(2-p)a}{2}\|\nabla u^{\mu,\zeta}\|_{L^\infty(\R^n_>)}^{(p-1)a}\!+
\|\nabla u^{\mu,\zeta}\|_{L^\infty(\R^n_>)}^{a}
) \|f\|_{q}^a  \|f\|_{q_2}^{p'(1-a)}, \forall x\in \R^n_>\,.\ea\ee
Observing that either $(p-1)a<1$ and $a<1$
we can apply Young's inequality, and find
\be\label{c7}\ba{ll}\dy \vs1
|\nabla u^{\mu,\zeta}(x)|^p\leq c \, \mu^\frac{a}{2}\|f\|_{q}^a (
\|f\|_{q_2}^{1-a}
B(\mu, \zeta)^{1-a}+
 \|f\|_{q_2}^{p'(1-a)})
+ 2\ve
\|\nabla u^{\mu,\zeta}\|_{L^\infty(\R^n_>)}^{p}\\\vs1 \dy +
c(\ve)(\mu^\frac{(2-p)a}{2}\|f\|_{q}^a
\|f\|_{q_2}^{1\!-\!a}
B(\mu,\!\zeta)^{1\!-\!a})^\frac{p}{p-(p\!-\!1)a}\!+
c(\ve)(\|f\|_{q}^a
\|f\|_{q_2}^{1\!-\!a}
B(\mu,\!\zeta)^{1\!-\!a})^\frac{p}{p-a}\\
\dy \hfill+
c(\ve)(\mu^\frac{(2-p)a}{2}\|f\|_{q}^a
\|f\|_{q_2}^{p'(1-a)})^\frac{p}{p-(p-1)a}+
c(\ve)(\|f\|_{q}^a
\|f\|_{q_2}^{p'(1-a)})^\frac{p}{p-a}\,,\ \forall x\in \R^n_>\,.\ea\ee
Raising both sides to $\frac 1p$, taking the supremum, and choosing a suitably small $\ve$, we easily find
\be\label{c8}\ba{ll}\dy \vs1
\|\nabla u^{\mu,\zeta}\|_{L^\infty(\R^n_>)}\leq c \, \mu^\frac{a}{2p}\|f\|_{q}^\frac{a}{p} (
\|f\|_{q_2}^{\frac{1-a}{p}}
B(\mu, \zeta)^\frac{1-a}{p}+ \|f\|_{q_2}^{\frac{1-a}{p-1}})\\\vs1 \dy +
c(\mu^\frac{(2-p)a}{2}\|f\|_{q}^a
\|f\|_{q_2}^{1-a}
B(\mu, \zeta)^{1-a})^\frac{1}{p-(p-1)a}+
c(\|f\|_{q}^a
\|f\|_{q_2}^{1-a}
B(\mu, \zeta)^{1-a})^\frac{1}{p-a}
\\
\dy \hfill+
c(\mu^\frac{(2-p)a}{2}\|f\|_{q}^a
\|f\|_{q_2}^{p'(1-a)})^\frac{1}{p-(p-1)a}+
c(\|f\|_{q}^a
\|f\|_{q_2}^{p'(1-a)})^\frac{1}{p-a}\,.\ea\ee
As $\|\nabla u^{\mu,\zeta}\|_{L^\infty(\R^n_<)}\leq \mu^\frac 12$,
we finally get
\be\label{c9}\ba{ll}\dy \vs1
\|\nabla u^{\mu,\zeta}\|_{\infty}\leq \mu^\frac 12+ c \, \mu^\frac{a}{2p}\|f\|_{q}^\frac{a}{p}
(\|f\|_{q_2}^{\frac{1-a}{p}}
B(\mu, \zeta)^\frac{1-a}{p}+ \|f\|_{q_2}^{\frac{1-a}{p-1}})\\ \vs1\dy +
c(\mu^\frac{(2-p)a}{2}\|f\|_{q}^a
\|f\|_{q_2}^{1-a}
B(\mu, \zeta)^{1-a})^\frac{1}{p-(p-1)a}+
c(\|f\|_{q}^a
\|f\|_{q_2}^{1-a}
B(\mu, \zeta)^{1-a})^\frac{1}{p-a}
\\
\dy \hfill+
c(\mu^\frac{(2-p)a}{2}\|f\|_{q}^a
\|f\|_{q_2}^{p'(1-a)})^\frac{1}{p-(p-1)a}+
c(\|f\|_{q}^a
\|f\|_{q_2}^{p'(1-a)})^\frac{1}{p-a}
\,.\ea\ee
Applying once again Lemma\,\ref{CmainT} and then estimate \eqref{c9} we have
\be\label{c10}\ba{ll}\vs1\dy
\|D^2 u^{\mu,\zeta}\|_q\leq c\mu^{\frac{2-p}{2}}\|f\|_q+ c\|f\|_q
 (\mu^\frac{a}{2p}\|f\|_{q}^\frac{a}{p}
(\|f\|_{q_2}^{\frac{1-a}{p}}
B(\mu, \zeta)^\frac{1-a}{p}+
 \|f\|_{q_2}^{\frac{1-a}{p-1}})
)^{2-p}\\\vs1 \dy +
c\|f\|_q(\mu^\frac{(2\!-\!p)a}{2}\!\|f\|_{q}^a
\|f\|_{q_2}^{1\!-\!a}\!
B(\mu,\!\zeta)^{1\!-\!a})^{\!\frac{2-p}{p-(p\!-\!1)a}}\!\!+
c\|f\|_q(\|f\|_{q}^a
\|f\|_{q_2}^{1\!-\!a}\!
B(\mu,\!\zeta)^{1\!-\!a})^\frac{2\!-\!p}{p\!-\!a}\\
\dy \hfill+
c\|f\|_q(\mu^\frac{(2-p)a}{2}\|f\|_{q}^a
\|f\|_{q_2}^{p'(1-a)})^\frac{2-p}{p-(p-1)a}+
c\|f\|_q(\|f\|_{q}^a
\|f\|_{q_2}^{p'(1-a)})^\frac{2-p}{p-a}
\,.
 \ea\ee
\vskip0.1cm\noindent
{\it \underline{The limit as $\mu\to 0$}} -
Following the arguments used in Step II and Step III of the  proof of Theorem\,\ref{mainTT}, we find that for any fixed $\zeta$ there exists a subsequence
of $\{(u^{\mu,\zeta}, \pi^{\mu,\zeta})\}$
 converging in suitable norms (see \eqref{loba1} and \eqref{sn1}) to a limit function $(u^\zeta,\pi^\zeta)$. For any $\zeta>0$, $u^\zeta$ satisfies estimate \eqref{mupasseqq1}. The expression of $\pi^{\zeta}$ is
 \be\label{pression3Mz}
\pi^\zeta=\widetilde\pi^\zeta+\psi^\zeta=\int_{\R^n}^*D^2_{y_iy_j}{\mathcal E}(x-y) \,\frac{D_{y_i} u_j^\zeta}{|\nabla u^\zeta|}{\atop^{\!\!\!2-p}}\,\,dy+
\int_{\R^n}\nabla_y {\mathcal E}(x-y)\cdot \,f\chi^\zeta\,dy\,.\ee
For any $\zeta>0$, the sequence $\{\pi^\zeta\}$ is bounded in $L^p(\R^n)$ (see \eqref{piuip}), and $\nabla\pi^\zeta$ satisfies estimate
 \eqref{Edual2}. Our next aim is  to get a bound, on $D^2 u^\zeta$ in $L^q(\R^n)$ and on $\nabla u^\zeta$ in $L^p(\R^n)$, uniformly
with respect to $\zeta$, which
implies a bound on $\pi^\zeta$in
$L^{p'}(\R^n)$, and then pass to
the limit as $\zeta\to \infty$. By
applying Hardy's inequality, we
write $B(\mu,\zeta)$, defined in
\eqref{c3} as follows
$$\ba{ll}\vs1\dy B(\mu, \zeta)=\|\frac{u^{\mu,\zeta}}{\zeta}\|_{L^p(K(\zeta))}+ \mu^\frac 12 (2\zeta)^\frac np\leq
\|\frac{u^{\mu,\zeta}-u^\zeta}{\zeta}\|_{L^p(K(\zeta))}+\|\frac{u^{\zeta}}{\zeta}\|_{L^p(K(\zeta))} +\mu^\frac 12 (2\zeta)^\frac np
\\\hfill  \vs1\dy
\leq \|\frac{u^{\mu,\zeta}-u^\zeta}{\zeta}\|_{L^p(K(\zeta))}+c\|\nabla u^{\zeta}\|_{L^p(|x|\geq \zeta))}+ \mu^\frac 12 (2\zeta)^\frac np\,.\ea$$
Therefore, using the strong convergence, as $\mu\to 0$, of $\{u^{\mu,\zeta}\}$ to $\{u^\zeta\}$ in $L^p(K(\zeta))$, from
\eqref{c10} we find the following bound
\be\label{c11}\ba{ll}\dy\vs1
\|D^2 u^\zeta\|_q=\liminf_{\mu\to 0} \|D^2 u^{\mu,\zeta}\|_q\leq
\|f\|_q\!+c\|f\|_q \|f\|_{q_2}^\frac{(1-a)(2-p)}{p-1}\!\!\\\vs1
\dy \hfill  +
c\|f\|_{q}^\frac{p(1\!-\!a)+a}{p-a}
\|f\|_{q_2}^\frac{p'(2-p)(1-a)}{p-a}\!\!\!+ c\|f\|_q^\frac{p(1-a)+a}{p-a}
\|f\|_{q_2}^\frac{(2-p)(1-a)}{p-a}\!
\liminf_{\mu\to 0}\!B(\mu,\!\zeta)^\frac{(2-p)(1-a)}{p-a}\\\vs1
\dy \hfill\leq
c\|f\|_q+c\|f\|_q \|f\|_{q_2}^\frac{(1-a)(2-p)}{p-1}
+
c\|f\|_{q}^\frac{p(1-a)+a}{p-a}
\|f\|_{q_2}^\frac{p'(2-p)(1-a)}{p-a}\\\hfill  \dy + c\|f\|_{q}^\frac{p(1-a)+a}{p-a}
\|f\|_{q_2}^\frac{(2-p)(1-a)}{p-a}
\|\nabla u^{\zeta}\|_{L^p(|x|\geq \zeta))}^\frac{(2-p)(1-a)}{p-a}
\ea\ee
On the other hand, by uniqueness, we know that for any $\zeta>0$ $u^\zeta$ coincides with
 the unique weak solution of \eqref{stokeschi}, which, from Lemma \ref{exists}
 satisfies
 \be\label{gradz}
 \|\nabla u^{\zeta}\|_p\leq c\| f\|_{q_2}^\frac{1}{p-1}\,.\ee
 Hence \eqref{c11} becomes
 \be\label{c12}
\|D^2 u^\zeta\|_q\leq
 c\|f\|_q(1+ \|f\|_{q_2}^\frac{(1-a)(2-p)}{p-1}+\|f\|_{q}^\frac{a(2-p)}{p-a}
\|f\|_{q_2}^{\frac{(2-p)(1-a)}{p-a}\frac{p}{p-1}})\,,\ee
with $a=\frac{nq}{nq+q-n}\,$,
 which gives the wanted uniform bound on the sequence $\{D^2 u^\zeta\}$ in $L^q(\R^n)$.
 \vskip0.2cm\noindent{\it\underline{The limit as $\zeta\to \infty$}} -
We want to pass to the limits as $\zeta\to\infty$ and
show that, as $\zeta\to\infty$, a subsequence of $\{(u^{\zeta}, \pi^{\zeta})\}$
  converges, in suitable norms,
 to a pair $(u, \pi)$, satisfying  Definition\,\ref{DS}
 Firstly note that, from estimate \eqref{gradz},
\be\label{pczp}
\big\||\nabla u^{\zeta}|^{p-2}\nabla u^{\zeta}\big\|_{p'}= \|\nabla u^{\zeta}\|_p\leq c\, \|f\|_{q_2}^\frac{1}{p-1}\,.\ee
Therefore, there exist a subsequence and two limit functions $\Psi$ and $u$, such that, in the limit as $\zeta\to\infty$
$$ u^{\zeta}  \rightharpoonup  u \mbox{ in } L^{\frac{np}{n-p}}(\R^n),$$
$$\nabla u^{\zeta} \rightharpoonup \nabla u \mbox{ in } L^{p}(\R^n),$$
\be\label{loba11}
|\nabla u^{\zeta}|^{p-2}\nabla u^{\zeta}\rightharpoonup  \Psi \ \mbox{ in \ } L^{p'}(\R^n)\,.
\ee
Moreover, from estimates \eqref{c12} there exists a
subsequence,  such that
$$\,D^2 u^{\zeta}\rightharpoonup D^2u \mbox{ in } L^q(\R^n)\,, \mbox{ as } \zeta\to \infty\,,$$
 with $D^2u$ satisfying estimate
  \be\label{c13}
\|D^2 u\|_q\leq
 c\|f\|_q(1+ \|f\|_{q_2}^\frac{(1-a)(2-p)}{p-1}+\|f\|_{q}^\frac{a(2-p)}{p-a}
\|f\|_{q_2}^{\frac{(2-p)(1-a)}{p-a}\frac{p}{p-1}})\,,\ee
where  $a=\frac{nq}{nq+q-n}$. Further, for any compact set $K\subset \R^n$ there exists a subsequence, depending on $K$,  such that
$$\nabla u^{\zeta}\rightarrow \nabla u \mbox{  in } L^{p}(K)\,, \mbox{ as } \zeta\to \infty\,.$$
This last convergence implies the almost everywhere convergence in $K$. Therefore, we also have
$$|\nabla u^{\zeta}|^{p-2}\nabla u^{\zeta} \rightarrow  |\nabla u|^{p-2}\nabla u \ \mbox{ a.e. in \ } K\,. $$
From  the weak convergence to
$\Psi$ in $L^{p'}(K)$ given in
\eqref{loba11}, we find that
$\Psi=|\nabla u|^{p-2}\nabla u$, on
each compact $K\subset \R^n$, which
ensures that \be\label{lob21}
|\nabla u^\zeta|^{p-2}\nabla
u^{\zeta} \rightharpoonup  |\nabla
u|^{p-2}\nabla u
 \ \mbox{  in \ } L^{p'}(\R^n).
\ee Finally, let us consider the
corresponding subsequence of
$\{\nabla \pi^{\zeta}\}$, with
$\pi^{\zeta}$ given by
\eqref{pression3Mz}. The sequence
$\{\pi^\zeta\}$ is uniformly
bounded in $L^{p'}(\R^n)$ and
satisfies \eqref{EDUALP}. Indeed
$\widetilde\pi$ is uniformly
bounded in $L^{p'}(\R^n)$ by using
the uniform bound \eqref{pczp} and
then applying the
Calder\'on-Zygmund theorem on
singular integrals, while
$\psi^\zeta$ is uniformly bounded
in $L^{p'}(\R^n)$ by using the
assumption $f\in L^{q_2}(\R^n)$ and
then applying the
Hardy-Littlewood-Sobolev theorem.
Further $\nabla \pi^\zeta$ is
uniformly bounded in $L^q(\R^n)$.
Hence there exist two fields
$\widetilde\pi$ and $\psi$ and a
subsequence such that, in the limit
as $\zeta\to\infty$, \be\label{cpz}
\widetilde\pi^\zeta\rightharpoonup
\widetilde\pi \textrm{ in }
L^{p'}(\R^n), \ee \be\label{cpz1}
\psi^\zeta\rightharpoonup \psi
\textrm{ in } L^{p'}(\R^n), \ee
\be\label{c14}
\nabla\widetilde\pi^\zeta\rightharpoonup
\nabla\widetilde\pi  \textrm{ in }
L^q(\R^n), \ee
\be\label{c15}\nabla\psi^{\zeta}\rightharpoonup
\nabla\psi \textrm{  in }
L^q(\R^n).\ee It remains to show
that \be\label{pression3Mzz}
\pi:=\widetilde\pi+\psi=\int_{\R^n}^*D^2_{y_iy_j}{\mathcal
E}(x-y) \,\frac{D_{y_i}
u_j}{|\nabla
u|}{\atop^{\!\!\!2-p}}\,\,dy+
\int_{\R^n}\nabla_y {\mathcal
E}(x-y)\cdot \,f\,dy\,.\ee This can
be obtained by repeating the
arguments used in the proof of
Theorem\,\ref{mainTT} (see Step II,
\eqref{pM}), since we have at
disposal the same convergences. For
brevity we omit the details. From
\eqref{lob21}, \eqref{c14},
\eqref{c15} and \eqref{stokeschi}
with $\mu=0$, it is routine to find
that $u$ is solution of
\eqref{stokes} in the sense of
Definition\,\ref{DS}. \chiu

 \vskip0.1cm
 {\bf Acknowledgment} - The authors are grateful
 to C.R. Grisanti for his valuable comments to
 the   paper.\par
The paper is
 performed under the
auspices of GNFM-INdAM.

\end{document}